# MULTIVARIATE ANALYSIS AND JACOBI ENSEMBLES: LARGEST EIGENVALUE, TRACY–WIDOM LIMITS AND RATES OF CONVERGENCE[1]

By Iain M. Johnstone

*Stanford University*

Let $A$ and $B$ be independent, central Wishart matrices in $p$ variables with common covariance and having $m$ and $n$ degrees of freedom, respectively. The distribution of the largest eigenvalue of $(A+B)^{-1}B$ has numerous applications in multivariate statistics, but is difficult to calculate exactly. Suppose that $m$ and $n$ grow in proportion to $p$. We show that after centering and scaling, the distribution is approximated to second-order, $O(p^{-2/3})$, by the Tracy–Widom law. The results are obtained for both complex and then real-valued data by using methods of random matrix theory to study the largest eigenvalue of the Jacobi unitary and orthogonal ensembles. Asymptotic approximations of Jacobi polynomials near the largest zero play a central role.

**1. Introduction.** It is a striking feature of the classical theory of multivariate statistical analysis that most of the standard techniques—principal components, canonical correlations, multivariate analysis of variance (MANOVA), discriminant analysis and so forth—are founded on the eigenanalysis of covariance matrices.

If, as is traditional, one assumes that the observed data follow a multivariate Gaussian distribution, then that theory builds on the eigenvalues and eigenvectors of one or two matrices following the Wishart distribution. Since the "single Wishart" problem can be viewed as a limiting case of the "double Wishart" one, the fundamental setting is that of the generalized eigenproblem $\det[B - \theta(A+B)] = 0$. In the idioms of MANOVA, $A$ repre-

Received November 2007.

[1]Supported in part by NSF Grants DMS-00-72661 and DMS-05-05303 and NIH RO1 CA 72028 and EB 001988.

sents the "within groups" or "error" covariance matrix, and $B$ the "between groups" or "hypothesis" covariance.

In each of the standard techniques, there is a conventional "null hypothesis"—independence, zero regression, etc. Corresponding test statistics may use either the full set of eigenvalues, as, for example, in the likelihood ratio test, or simply the extreme eigenvalues, as in the approach to inference advanced by S. N. Roy.

This paper focuses on the largest eigenvalue, or "latent root," and in particular on its distribution under the null hypothesis, in other words when the two Wishart matrices $A$ and $B$ are independent, central, and have common covariance matrix.

Even under the assumption of Gaussian data, the null distribution of the largest root is difficult to work with. It is expressed in terms of a hypergeometric function of matrix argument, with no general and simple closed form. It depends on three parameters—the common dimension of the two Wishart matrices and their respective degrees of freedom. Traditional textbooks have included tables of critical points which, due to the three parameters, can run up to twenty-five pages [Morrison (2005), Timm (1975)]. Traditional software packages have often used a one-dimensional $F$ distribution approximation that can be astonishingly inaccurate for dimensions greater than two or three. Recently, as will be reviewed below, some exact algorithms have been made available, but they are not yet in wide use. One can speculate that the use of largest root tests has been limited in part by the lack of a simple, serviceable approximation.

The goal of this paper is to provide such an approximation, which turns out to be expressed in terms of the Tracy–Widom distribution $F_1$ of random matrix theory. This distribution is free of parameters, and can be tabulated or calculated on the fly; it plays here a role analogous to that of the standard normal distribution $\Phi$ in central limit approximations. The three Wishart parameters appear in the centering and scaling constants for the largest eigenvalue, for which we give readily computable formulas.

The approximation is an asymptotic one, developed using the models and techniques of random matrix theory in which the dimension $p$ increases to infinity, and the degrees of freedom parameters grow in proportion to $p$. A pleasant surprise is that the approximation has a "second-order" accuracy; in that sense loosely reminiscent of the $t$-approximation to normal. The traditional percentage points in the upper tail of the null distribution— 90%, 95%, etc.—are reasonably well approximated for $p$ as small as 5. In a companion paper Johnstone (2009), it is argued that over the entire range of the parameters (i.e., $p$ as small as 2), the Tracy–Widom approximation can yield a first screening of significance level for the largest root test that may be adequate in many, and perhaps most, applied settings.

Some words about the organization of the paper. The remainder of this Introduction develops the double Wishart setting and states the approximation result, first for real-valued data, and then for complex-valued data, the latter involving the Tracy–Widom $F_2$ distribution. Section 2 collects some of the statistical settings to which the Tracy–Widom approximation applies, along with a "dictionary" that translates the result into each setting.

The remainder of the paper develops the proofs, using methods of random matrix theory (RMT). Section 3 reformulates the results in the language of RMT, to say that the scaled Jacobi unitary and orthogonal ensembles converge to Tracy–Widom at the soft upper edge. Section 4 gives a detailed outline of the proof, noting points of novelty. As is conventional, the unitary (complex) case is treated first (Section 7), and then used as a foundation for the orthogonal (real) setting of primary interest in statistics (Section 8). Everything is based on Plancherel–Rotach asymptotics of Jacobi polynomials near their largest zero; this is developed in Sections 5 and 6 using the Liouville–Green approach to the corresponding differential equation. Some of the results of this paper were announced in Johnstone (2007).

1.1. *Statement of results.* Let $X$ be an $m \times p$ normal data matrix: each row is an independent observation from $N_p(0, \Sigma)$. A $p \times p$ matrix $A = X'X$ is then said to have a Wishart distribution $A \sim W_p(\Sigma, m)$. Let $B \sim W_p(\Sigma, n)$ be independent of $A \sim W_p(\Sigma, m)$. Assume that $m \geq p$; then $A^{-1}$ exists and the nonzero eigenvalues of $A^{-1}B$ generalize the univariate $F$ ratio. The scale matrix $\Sigma$ has no effect on the distribution of these eigenvalues, and so without loss of generality suppose that $\Sigma = I$.

The matrix analog of a Beta variate is based on the eigenvalues of $(A + B)^{-1}B$, and leads to

DEFINITION 1 [Mardia, Kent and Bibby (1979), page 84]. Let $A \sim W_p(I, m)$ be independent of $B \sim W_p(I, n)$, where $m \geq p$. Then the largest eigenvalue $\theta$ of $(A + B)^{-1}B$ is called the *greatest root statistic* and a random variate having this distribution is denoted $\theta_1(p, m, n)$, or $\theta_{1,p}$ for short.

Since $A$ is positive definite, $0 < \theta < 1$. Equivalently $\theta_1(p, m, n)$ is the largest root of the determinantal equation

$$\det[B - \theta(A + B)] = 0. \tag{1}$$

Specific examples will be given below, but in general the parameter $p$ refers to dimension, $m$ to the "error" degrees of freedom and $n$ to the "hypothesis" degrees of freedom. Thus $m + n$ represents the "total" degrees of freedom.

The greatest root distribution has the property

$$\theta_1(p, m, n) \stackrel{\mathcal{D}}{=} \theta_1(n, m + n - p, p), \tag{2}$$

useful in particular in the case when $n < p$ [e.g., Mardia, Kent and Bibby (1979), page 84].

Assume $p$ is even and that $p, m = m(p)$ and $n = n(p) \to \infty$ together in such a way that

$$\text{(3)} \qquad \lim_{p \to \infty} \frac{\min(p, n)}{m + n} > 0, \qquad \lim_{p \to \infty} \frac{p}{m} < 1.$$

A consequence of our main result, stated more completely below, is that with appropriate centering and scaling, the logit transform $W_p = \operatorname{logit} \theta_{1,p} = \log(\theta_{1,p}/(1 - \theta_{1,p}))$ is approximately Tracy–Widom distributed:

$$\text{(4)} \qquad \frac{W_p - \mu_p}{\sigma_p} \overset{\mathcal{D}}{\Rightarrow} Z_1 \sim F_1.$$

The distribution $F_1$ was found by Tracy and Widom (1996) as the limiting law of the largest eigenvalue of a $p \times p$ Gaussian symmetric matrix; further information on $F_1$ is reviewed, for example, in Johnstone (2001). Its appearance here is an instance of the universality properties expected for largest eigenvalue distributions in random matrix theory [e.g., Deift (2007), Deift and Gioev (2007), Deift et al. (2007)].

The centering and scaling parameters are given by

$$\text{(5)} \qquad \begin{aligned} \mu_p &= 2 \log \tan\left(\frac{\varphi + \gamma}{2}\right), \\ \sigma_p^3 &= \frac{16}{(m+n-1)^2} \frac{1}{\sin^2(\varphi + \gamma) \sin \varphi \sin \gamma}, \end{aligned}$$

where the angle parameters $\gamma, \varphi$ are defined by

$$\text{(6)} \qquad \begin{aligned} \sin^2\left(\frac{\gamma}{2}\right) &= \frac{\min(p, n) - 1/2}{m + n - 1}, \\ \sin^2\left(\frac{\varphi}{2}\right) &= \frac{\max(p, n) - 1/2}{m + n - 1}. \end{aligned}$$

As will be discussed later, the "correction factors" of $-\frac{1}{2}$ and $-1$ yield a second-order rate of convergence that has important consequences for the utility of the approximation in practice. Indeed, our main result can be formulated as follows.

THEOREM 1. *Assume that $m(p), n(p) \to \infty$ as $p \to \infty$ through even values of $p$ according to (3). For each $s_0 \in \mathbb{R}$, there exists $C > 0$ such that for $s \geq s_0$,*

$$|P\{W_p \leq \mu_p + \sigma_p s\} - F_1(s)| \leq C p^{-2/3} e^{-s/2}.$$

*Here $C$ depends on $(\gamma, \varphi)$ and also on $s_0$ if $s_0 < 0$.*

1.2. *Exact expressions.* Assume that $m, n \geq p$ and that $A \sim W_p(I, m)$ independently of $B \sim W_p(I, n)$. The joint density of the eigenvalues $1 \geq \theta_1 \geq \theta_2 \geq \cdots \geq \theta_p \geq 0$ of $(A+B)^{-1}B$, or equivalently, of the roots of $\det[B - \theta(A+B)] = 0$, simultaneously derived in 1939 by Fisher, Girshick, Hsu, Mood and Roy, is given by Muirhead (1982), page 112:

$$(7) \qquad f(\theta) = c_1 \prod_{i=1}^{p} (1-\theta_i)^{(m-p-1)/2} \theta_i^{(n-p-1)/2} \prod_{i<j}^{p} |\theta_i - \theta_j|.$$

The normalizing constant $c_1 = c_1(p, m, n)$ involves the multivariate gamma function; we shall not need it here.

Exact evaluation of the marginal distribution of the largest root $\theta_1$ is not a simple matter. Constantine (1963) showed that the marginal distribution could be expressed in terms of a hypergeometric function of matrix argument. Let $t = (n-p-1)/2$. Then

$$(8) \qquad P\{\theta_{1,p} \leq x\} = c_2 x^{pm/2} {}_2F_1\left(\frac{m}{2}, -t; \frac{m+p+1}{2}; xI\right).$$

When $t$ is a nonnegative integer, there is a terminating series [Koev (n.d.), Muirhead (1982), page 483, Khatri (1972)] in terms of zonal polynomials $C_\kappa$:

$$(9) \qquad P\{\theta_{1,p} \leq x\} = x^{pm/2} \sum_{k=0}^{pt} \sum_{\kappa \vdash k, \kappa_1 \leq t} \frac{(m/2)_\kappa C_\kappa((1-x)I)}{k!},$$

where $\kappa \vdash k$ signifies that $\kappa = (\kappa_1, \ldots, \kappa_n)$ is a partition of $\kappa$, and $(\frac{m}{2})_\kappa$ is a generalized hypergeometric coefficient. Further details and definitions may be found in Koev (n.d.) and Muirhead (1982). Johnstone (2009) lists further references in the literature developing tables of the distribution of $\theta_1(p, m, n)$.

Recently, Koev and Edelman (2006) have exploited recursion relations among Jack functions to develop efficient evaluations of hypergeometric functions of matrix argument. Current MATLAB software implementations (Koev, private communication) allow convenient—up to 1-sec computation time—evaluation of (8) for $m, n, p \leq 17$, and (9) for $m, n, p \leq 40$ when $t$ is integer.

1.3. *Numerical illustrations.* Table 1 and Figure 1 show results of some simulations to test the Tracy–Widom approximation. A companion paper [Johnstone (2009)] has further information on the quality of the distributional approximation.

To the left of the vertical line are three situations in which $m = 8p$ and $n = 2p$, that is, where the error and hypothesis degrees of freedom are comfortably larger than dimension. In the second setting, to the right of the line, this is no longer true: $m = 2p$ and $n = p$.

TABLE 1
*First column shows the percentiles of the $F_1$ limit distribution corresponding to fractions in second column. Next three columns show estimated cumulative probabilities for $\theta_1$ in $R = 10{,}000$ repeated draws from the two Wishart setting of Definition 1, with indicated values of $(p, m, n)$. The following three columns show estimated cumulative probabilities for $w = \log \theta/(1-\theta)$ again in $R = 10{,}000$ draws with the indicated values of $(p, m, n)$. Final column gives approximate standard errors based on binomial sampling. Bold font highlights some conventional significance levels. The Tracy–Widom distribution $F_1$ was evaluated on a grid of 121 points $-6(0.1)6$ using the Mathematica package p2Num written by Craig Tracy. Remaining computations were done in MATLAB, with percentiles obtained by inverse interpolation, and using randn() for normal variates and norm() to evaluate the largest eigenvalue of the matrices appearing in Definition 1*

| Percentile | TW | | $p, n$ | 20,40 | 5,10 | 2,4 | | 50,50 | 5,5 | 2,2 | |
|---|---|---|---|---|---|---|---|---|---|---|---|
| | | | $m$ | 160 | 40 | 16 | | 100 | 10 | 4 | |
| | | $\mu_\theta$ | | 0.49 | 0.48 | 0.44 | $\mu$ | 2.06 | 1.93 | 1.69 | |
| | | $\sigma_\theta$ | | 0.02 | 0.06 | 0.18 | $\sigma$ | 0.127 | 0.594 | 1.11 | 2 * SE |
| $-3.90$ | 0.01 | | | 0.010 | 0.008 | 0.000 | | 0.007 | 0.002 | 0.010 | (0.002) |
| $-3.18$ | 0.05 | | | 0.052 | 0.049 | 0.009 | | 0.042 | 0.023 | 0.037 | (0.004) |
| $-2.78$ | 0.10 | | | 0.104 | 0.099 | 0.046 | | 0.084 | 0.062 | 0.074 | (0.006) |
| $-1.91$ | 0.30 | | | 0.311 | 0.304 | 0.267 | | 0.289 | 0.262 | 0.264 | (0.009) |
| $-1.27$ | 0.50 | | | 0.507 | 0.506 | 0.498 | | 0.499 | 0.495 | 0.500 | (0.010) |
| $-0.59$ | 0.70 | | | 0.706 | 0.705 | 0.711 | | 0.708 | 0.725 | 0.730 | (0.009) |
| 0.45 | **0.90** | | | **0.904** | **0.910** | **0.911** | | **0.905** | **0.919** | **0.931** | (0.006) |
| 0.98 | **0.95** | | | **0.950** | **0.955** | **0.958** | | **0.953** | **0.959** | **0.966** | (0.004) |
| 2.02 | **0.99** | | | **0.990** | **0.992** | **0.995** | | **0.990** | **0.991** | **0.993** | (0.002) |

In the first setting, the largest eigenvalue distribution is concentrated around $\mu_\theta \approx 0.5$ and with scale $\sigma_\theta$ small enough that the effect of the boundary at 1 is hardly felt. In this setting, the logit transform $w = \log \theta/(1-\theta)$ is less important for improving the quality of the approximation. Indeed, the three first columns show the result of using

$$\mu_\theta = \frac{e^{\mu_p}}{1 + e^{\mu_p}}, \qquad \sigma_\theta = \mu_\theta(1 - \mu_\theta)\sigma_p.$$

1.4. *Complex-valued data.* Data matrices $X$ based on complex-valued data arise frequently, for example, in signal processing applications [e.g., Tulino and Verdu (2004)]. If the rows of $X$ are drawn independently from a complex normal distribution $\mathbb{C}N(\mu, \Sigma)$ [see, e.g., James (1964), Section 7], then we say $A = \bar{X}'X \sim \mathbb{C}W_p(\Sigma, n)$. Here $\bar{X}'$ denotes the conjugate transpose of $X$.

In parallel with the real case definition, if $A \sim \mathbb{C}W_p(I, m)$ and $B \sim \mathbb{C}W_p(I, n)$ are independent, then the joint density of the eigenvalues $1 \geq \theta_1 \geq \theta_2 \geq \cdots \geq$

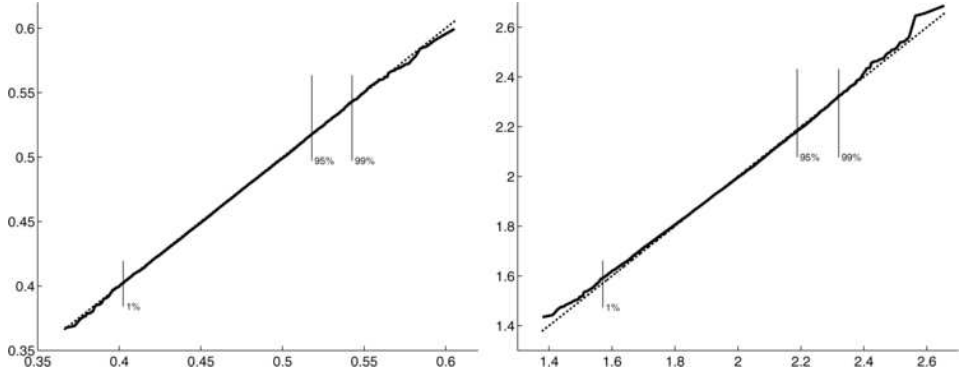

FIG. 1. *First panel: probability plots of $R = 10{,}000$ observed replications of $\theta_1$, largest root of (1), for $p = 20$, $n = 40$, $m = 160$. That is, the $10{,}000$ ordered observed values of $\theta_1$ are plotted against $F_1^{-1}((i - 0.5)/R)$; $i = 1, \ldots, R$. The vertical lines show 1st, 95th and 99th percentiles. The dotted line is the $45°$ line of perfect agreement of empirical law with asymptotic limit. Second panel: same plots for $p = n = 50$, $m = 100$, but now on logit scale, plotting $W_p = \log(\theta_1/(1 - \theta_1))$.*

$\theta_p \geq 0$ of $(A + B)^{-1}B$, or equivalently, the roots of $\det[B - \theta(A + B)] = 0$, is given, for example, by James (1964),

(10) $$f(\theta) = c \prod_{i=1}^{p}(1 - \theta_i)^{m-p}\theta_i^{n-p}\prod_{i<j}(\theta_i - \theta_j)^2.$$

The largest eigenvalue $\theta^C(p, m, n)$ of $(A + B)^{-1}B$ is called the greatest root statistic, with distribution $\theta^C(p, m, n)$. The property (2) carries over to the complex case.

Again let $W^C = \text{logit}\,\theta_p^C = \log(\theta_p^C/(1 - \theta_p^C))$.

THEOREM 2. *Assume that $m(p), n(p) \to \infty$ as $p \to \infty$ according to (3). For each $s_0 \in \mathbb{R}$, there exists $C > 0$ such that for $s \geq s_0$,*

$$|P\{W_p^C \leq \mu_p^C + \sigma_p^C s\} - F_2(s)| \leq C p^{-2/3} e^{-s/2}.$$

*Here $C$ depends on $(\gamma, \varphi)$ and also on $s_0$ if $s_0 < 0$.*

The limiting distribution is now the *unitary* Tracy–Widom distribution [Tracy and Widom (1994)]. To describe the complex centering and scaling constants, we introduce a parameterization basic to the paper:

$$N = \min(n, p), \qquad \alpha = m - p, \qquad \beta = |n - p|.$$

Then $\mu^C, \sigma^C$ use weighted averages based on the parameter sets $(N, \alpha, \beta)$ and $(N - 1, \alpha, \beta)$:

$$\mu^C = \frac{\tau_N^{-1} u_N + \tau_{N-1}^{-1} u_{N-1}}{\tau_N^{-1} + \tau_{N-1}^{-1}}, \qquad (\sigma^C)^{-1} = \tfrac{1}{2}(\tau_N^{-1} + \tau_{N-1}^{-1}),$$

where
$$w_N = 2\log\tan\left(\frac{\varphi_N + \gamma_N}{2}\right),$$
$$\omega_N^3 = \frac{16}{(2N + \alpha + \beta + 1)^2}\frac{1}{\sin^2(\varphi_N + \gamma_N)\sin\varphi_N \sin\gamma_N}$$

and
$$\sin^2\left(\frac{\gamma_N}{2}\right) = \frac{N + 1/2}{2N + \alpha + \beta + 1},$$
$$\sin^2\left(\frac{\varphi_N}{2}\right) = \frac{N + \beta + 1/2}{2N + \alpha + \beta + 1}.$$

Quantities $w_{N-1}, \omega_{N-1}$ are based on $\varphi_{N-1}, \gamma_{N-1}$ with $N-1$ substituted everywhere for $N$, but with $\alpha, \beta$ unchanged.[2]

The remarks made earlier about exact expressions for largest eigenvalue distributions have analogs in the complex case—see the references cited earlier, Dumitriu and Koev (2008) and determinantal identities (2.10) and (2.11) in Koev (n.d.).

1.5. *Remark on software.* Given a routine to compute the Tracy–Widom distribution, it is a simple matter to code and use the formulas given in this paper. Some references to then-extant software were given in Johnstone (2007). Further detail is planned for the completed version of Johnstone (2009).

**2. Related statistical settings and implications.** In the first part of this section, we list five common settings in multivariate statistics to which the largest eigenvalue convergence result applies, along with the parameterizations appropriate to each.

2.1. *Double Wishart models.*

2.1.1. *Canonical correlation analysis.* Suppose that there are $n$ observations on each of $p + q$ variables. For definiteness, assume that $p \leq q$. The first $p$ variables are grouped into an $n \times p$ data matrix $X = [\mathbf{x}_1\ \mathbf{x}_2\ \cdots\ \mathbf{x}_p]$ and the last $q$ into $n \times q$ matrix $Y = [\mathbf{y}_1\ \mathbf{y}_2\ \cdots\ \mathbf{y}_q]$. Write $S_{XX} = X^T X$, $S_{XY} = X^T Y$ and $S_{YY} = Y^T Y$ for the cross-product matrices. Canonical correlation analysis (CCA), or more precisely, the zero-mean version of CCA, seeks the

---
[2] This use of the notation $w_N$, $\omega_N$ is local to this Introduction and the Remark concluding Section 7.1 and not to be confused with other uses in the detailed proofs.

linear combinations $a^T x$ and $b^T y$ that are most highly correlated, that is, to maximize

$$r = \text{Corr}(a^T x, b^T y) = \frac{a^T S_{XY} b}{\sqrt{a^T S_{XX} a} \sqrt{b^T S_{YY} b}}. \tag{11}$$

This leads to a maximal correlation $r_1$ and associated canonical vectors $a_1$ and $b_1$, usually each taken to have unit length. The procedure may be iterated, restricting the search to vectors orthogonal to those already found:

$$r_k = \max\{a^T S_{XY} b : a^T S_{XX} a = b^T S_{YY} b = 1, \text{ and}$$
$$a^T S_{XX} a_j = b^T S_{YY} b_j = 0, \text{ for } 1 \leq j < k\}.$$

The successive canonical correlations $r_1 \geq r_2 \geq \cdots \geq r_p \geq 0$ may be found as the roots of the determinantal equation

$$\det(S_{XY} S_{YY}^{-1} S_{YX} - r^2 S_{XX}) = 0 \tag{12}$$

[see, e.g., Mardia, Kent and Bibby (1979), page 284]. A typical question in application is then how many of the $r_k$ are significantly different from zero.

Substitute the cross-product matrix definitions into the CCA determinantal equation (12) to obtain $\det(X^T Y (Y^T Y)^{-1} Y^T X - r^2 X^T X) = 0$. Let $P$ denote the $n \times n$ orthogonal projection matrix $Y(Y^T Y)^{-1} Y^T$ and $P^\perp = I - P$ its complement. Then with $B = X^T P X$ and $A = X^T P^\perp X$, (12) becomes

$$\det(B - r^2(A + B)) = 0. \tag{13}$$

Now assume that $Z = [X \ Y]$ is an $n \times (p + q)$ normal data matrix with mean zero. The covariance matrix is partitioned

$$\Sigma = \begin{bmatrix} \Sigma_{XX} & \Sigma_{XY} \\ \Sigma_{YX} & \Sigma_{YY} \end{bmatrix}.$$

Under these Gaussian assumptions, the $X$ and $Y$ variable sets will be independent if and only if $\Sigma_{XY} = 0$. This is equivalent to asserting that the population canonical correlations all vanish: $\rho_1 = \cdots = \rho_p = 0$.

The canonical correlations $(\rho_1, \ldots, \rho_p)$ are invariant under block diagonal transformations $(x_i, y_i) \to (Bx_i, Cy_i)$ of the data (for $B$ and $C$ nonsingular $p \times p$ and $q \times q$ matrices, resp.). It follows that under the null hypothesis $H_0 : \Sigma_{XY} = 0$, the distribution of the canonical correlations can be found (without loss of generality) by assuming that $\Sigma_{XX} = I_p$ and $\Sigma_{YY} = I_q$. In this case, the matrices $A$ and $B$ of (13) are independent with $B \sim W_p(q, I)$ and $A \sim W_p(n - q, I)$.

From the definition, the largest squared canonical correlation $\theta_1 = r_1^2$ has the $\theta(p, n - q, q)$ distribution under the null hypothesis $\Sigma_{XY} = 0$.

*Mean correction.* In practice, it is more common to allow each variable to have a separate, unknown mean. One forms the variable means $\bar{x}_i = n^{-1}\sum_{k=1}^{n} x_{i,k}$ and replaces $\mathbf{x}_i$ by $\mathbf{x}_i - \bar{x}_i\mathbf{1}$, and similarly for the second set of variables $\mathbf{y}_j$. The entries $S_{XY}$, etc. in (11) are now blocks in the partitioned sample covariance matrix: if $P_c = I_n - n^{-1}\mathbf{1}\mathbf{1}^T$, then

$$S_{XY} = (P_c Y)^T (P_c X) = Y^T P_c X, \qquad S_{XX} = X^T P_c X, \text{ etc.}$$

For the distribution theory, suppose that $Z = [X\ Y]$ is an $n \times (p+q)$ normal data matrix with mean $(\mu_X\ \mu_Y)$ and covariance $\Sigma$. Applying mean-corrected CCA as above, then under $\Sigma_{XY} = 0$, the largest squared canonical correlation $\theta_1 = r_1^2$ has distribution $\theta_1(p, n'-q, q)$, where $n' = n - 1$. Indeed, let $H'$ be the upper $(n-1) \times n$ block of an orthogonal matrix with $n$th row equal to $n^{-1/2}\mathbf{1}^T$. Then $Z' = H'Z$ turns out [e.g., Mardia, Kent and Bibby (1979), page 65] to be a normal data matrix with mean 0, covariance $\Sigma$ and sample size $n - 1$, to which our mean-zero discussion above applies. Since $H'^T H' = P_c$, the mean-zero prescription (11) applied to $[X'\ Y']$ yields the same canonical correlations as does the usual mean centered approach applied to $[X\ Y]$.

2.1.2. *Angles and distances between subspaces.* The cosine of the angle between two vectors $u, v \in \mathbb{R}^n$ is given by

$$\sigma(u, v) = |u^T v|/(\|u\|_2 \|v\|_2).$$

Consequently, (11) becomes $r = \sigma(Xa, Yb)$. Writing $\mathcal{X}$ and $\mathcal{Y}$ for the subspaces spanned by the columns of $X$ and $Y$, then the canonical correlations $r_k = \cos\vartheta_k$ are just the cosines of the principal angles between $\mathcal{X}$ and $\mathcal{Y}$ [e.g., Golub and Van Loan (1996), page 603].

The closeness of two equidimensional subspaces can be measured by the largest angle between vectors in the two spaces:

$$d(\mathcal{X}, \mathcal{Y}) = \min_{a,b} \sigma(Xa, Yb) = r_p,$$

the *smallest* canonical correlation. This is equivalent to the 2-norm of the distance between orthoprojections on $\mathcal{X}$ and $\mathcal{Y}$: $\|P_\mathcal{X} - P_\mathcal{Y}\|_2 = \sin\theta_p = \sqrt{1 - r_p^2}$.

*Random subspaces.* A standard way to realize a draw from the uniform (Haar) distribution on the Grassmann manifold of $p$-dimensional subspaces of $\mathbb{R}^n$ is to let $\mathcal{X} = \text{span}(X)$, with the entries of the $n \times p$ matrix $X$ being i.i.d. standard Gaussian. If $X_{n \times p}$ and $Y_{n \times q}$ are two such independent Gaussian matrices, then the squared cosines of the principal angles between $\mathcal{X}$ and $\mathcal{Y}$ are given by the roots of (13), with $A \sim W_p(n-q, I)$ independently of $B \sim W_p(q, I)$. In the language of the next section, the Jacobi orthogonal ensemble thus arises as the distribution of the squared principal cosines between two random subspaces. Similar statements hold for complex-valued Gaussian data matrices $X$, $Y$ and the Jacobi unitary ensemble [cf. Collins (2005), Theorem 2.2, and Absil, Edelman and Koev (2006)].

2.1.3. *Multivariate linear model.* In the standard generalization of the linear regression model to allow for multivariate responses, it is assumed that

$$Y = XB + U,$$

where $Y(n \times p)$ is an observed matrix of $p$ response variables on each of $n$ individuals, $X(n \times q)$ is a known design matrix, $B(q \times p)$ is a matrix of unknown regression parameters and $U$ is a matrix of unobserved random disturbances. For distribution theory it is assumed that $U$ is a normal data matrix, so that the rows are independent Gaussian, each with mean 0 and common covariance $\Sigma$.

Consider a null hypothesis of the form $CBM = 0$. Here it is assumed that $C(g \times q)$ has rank $g$. The rows of $C$ make assertions about the effect of linear combinations of the "independent" variables on the regression: the number of hypothesis degrees of freedom is $g$. The matrix $M(p \times r)$ is taken to have rank $r$. The columns of $M$ focus attention of particular linear combinations of the dependent variables: the "dimension" of the null hypothesis equals $r$. The union-intersection test of this null hypothesis is based on the greatest root of $\theta$ of $H(H + E)^{-1}$ for the independent Wishart matrices $H$ and $E$ described, for example, in Mardia, Kent and Bibby (1979), page 162. Under the null hypothesis, $\theta \sim \theta(r, n-q, g)$. The companion paper Johnstone (2009) focuses in greater detail on the application of Theorem 1 in the multivariate linear model.

2.1.4. *Equality of covariance matrices.* Suppose that independent samples from two normal distributions $N_p(\mu_1, \Sigma_1)$ and $N_p(\mu_2, \Sigma_2)$ lead to covariance estimates $\hat{\Sigma}_i$ which are independent and Wishart distributed on $n_i$ degrees of freedom: $A_i = n_i \hat{\Sigma}_i \sim W_p(n_i, \Sigma_i)$ for $i = 1, 2$. Then the largest root test of the null hypothesis $H_0 : \Sigma_1 = \Sigma_2$ is based on the largest eigenvalue $\theta$ of $(A_1 + A_2)^{-1} A_2$, which under $H_0$ has the $\theta(p, n_1, n_2)$ distribution [Muirhead (1982), page 332].

2.1.5. *Multiple discriminant analysis.* Suppose that there are $g$ populations, the $i$th population being assumed to follow a $p$-variate normal distribution $N_p(\mu_i, \Sigma)$, with the covariance matrix assumed to be unknown, but common to all populations. A sample of size $n_i$ is available from the $i$th population, yielding a total $n = \sum n_i$ observations. Multiple discriminant analysis uses the "within groups" and "between groups" sums of squares and products matrices $W$ and $B$ to construct linear discriminant functions based on eigenvectors of $W^{-1} B$. A test of the null hypothesis that discrimination is not worthwhile ($\mu_1 = \cdots = \mu_g$) can be based, for example, on the largest root of $W^{-1} B$, which leads to use of the $\theta(p, n-g, g-1)$ distribution [Mardia, Kent and Bibby (1979), pages 318 and 138].

Table 2

| Setting | | $p$ | $m$ | $n$ |
|---|---|---|---|---|
| CCA | $[X\ Y] \sim N_{p+q}(0, I_n \otimes \Sigma)$<br>$H_0: \Sigma_{XY} = 0$ | $p$ | $n-q$ | $q$ |
| Multivariate | $\underset{n \times p}{Y} = \underset{q \times p}{X\ B} + \underset{n \times p}{U}$ | $r$ | $n-q$ | $g$ |
| Linear<br>　model | $H_0: \underset{g \times q}{C}\ \underset{q \times p}{\beta}\ \underset{p \times r}{M} = 0$ | ↑<br>dimen | ↑<br>error d.f. | ↑<br>hypoth. d.f. |
| Equality<br>　of covariance | $n_i \hat{\Sigma}_i \sim W_p(n_i, \Sigma_i)$<br>$H_0: \Sigma_1 = \Sigma_2$ | $p$ | $n_1$ | $n_2$ |
| Mult.<br>　discrim. | $n_i$ obs on $g$ pops $N_p(\mu_i, \Sigma)$<br>$i = 1, \ldots, g$ | $p$ | $n-g$ | $g-1$ |

Table 2 summarizes the correspondences between the parameters in these various models and those used in Theorem 1.

2.2. *Discussion and implications.*

*Limiting empirical spectrum.* The empirical distribution of eigenvalues $\theta_i$ of $(A+B)^{-1}B$ is defined to be

$$F_p(\theta) = p^{-1} \#\{i : \theta_i \leq \theta\}.$$

Wachter (1980) obtained the limiting distribution of $F_p$ in an asymptotic regime (3) in which $m$ and $n$ grow proportionally with $p$. We recall Wachter's result, in the new parameterization given by (6). Suppose, for convenience, that $p \leq n$, and let

$$\theta_{\pm} = \sin^2\left(\frac{\varphi \pm \gamma}{2}\right),$$

or, more precisely, the limit as $p \to \infty$ under assumption (3). Then for each $\theta \in [0, 1]$, $F_p(\theta) \to \int_0^\theta f(\theta')\,d\theta'$, where the limiting density has the form

$$f(\theta) = \frac{c\sqrt{(\theta_+ - \theta)(\theta - \theta_-)}}{\theta(1 - \theta)}, \qquad c = 2\pi \sin^2(\gamma/2).$$

This is the analog for two Wishart matrices of the celebrated semicircle law for square symmetric matrices, and the Marčenko–Pastur quarter-circle law for a single Wishart matrix [for references, see, e.g., Johnstone (2001)].

In the canonical correlation setting—$p$ and $q$ variables and $n$ samples—the parameters $\theta_{\pm}$ represent the limiting maximum and minimum squared correlation. They are expressed in terms of the half-angles $\gamma/2$ and $\varphi/2$, which for $p/n$ and $q/n$ fairly small are roughly

$$\gamma/2 \doteq \sqrt{p/n}, \qquad \varphi/2 \doteq \sqrt{q/n}.$$

*Smallest eigenvalue.* Assume that $A \sim W_p(I, m)$ independently of $B \sim W_p(I, n)$ and that both $m, n \geq p$. In this case, all $p$ eigenvalues of $A^{-1}B$ are positive a.s., and we have the identity

(14) $$\theta_1\{(A+B)^{-1}B\} = 1 - \theta_p\{(A+B)^{-1}A\},$$

where $\theta_k(C)$ denotes the $k$th ordered eigenvalue of $C$, with $\theta_1$ being smallest.

Let $\theta^-(p, m, n)$ denote a random variable having the distribution of the *smallest* eigenvalue of $(A+B)^{-1}B$. Clearly

$$\theta^-(p, m, n) \stackrel{\mathcal{D}}{=} 1 - \theta(p, n, m)$$

and if we set $\theta' = \theta(p, n, m)$, then

$$W_p^- = \log \frac{\theta^-}{1 - \theta^-} = -\log \frac{\theta'}{1 - \theta'}.$$

If we therefore set

$$\mu^-(p, m, n) = -\mu(p, n, m), \qquad \sigma^-(p, m, n) = \sigma(p, n, m),$$

where $\mu(p, n, m)$ and $\sigma(p, n, m)$ are given by (5)–(6) (with $n$ and $m$ interchanged), we have convergence to a Tracy–Widom distribution reflected at 0:

$$(W_p^- - \mu_p^-)/\sigma_p^- \stackrel{\mathcal{D}}{\Rightarrow} -F_1,$$

or, writing $\bar{F}(t) = 1 - F_1(t)$ for the complementary Tracy–Widom distribution function,

$$|P\{W_p^- \leq \mu_p^- - \sigma_p^- t\} - \bar{F}_1(t)| \leq Cp^{-2/3}e^{-ct}.$$

This form highlights the phenomenon that convergence for the distribution of $\theta^-(p, m, n)$ is best in the *left* tail. As with $\theta_1(p, m, n)$, the approximation is best in the part of the distribution furthest from the bulk of the eigenvalues.

*Analogy with* $t$. Here is an admittedly loose analogy between the null distribution of the largest eigenvalue and that of the $t$-statistic. Both cases assume Gaussian data, but in the $t$ case, the test is on the mean $\mu$, while for $\theta_1$ it concerns the covariance structure. In both settings, the exact null distribution is known, but one is interested in the rate of convergence to the limiting distribution which is used for approximation. Table 3 compares our result—in the canonical correlations version—with a standard fact about the Gaussian approximation to the $t$ distribution: if the parent distribution is Gaussian, then the convergence is second-order, in that the error term is of order $1/n$ rather than the first-order error $1/\sqrt{n}$ associated with central limit theorem convergence.

TABLE 3

|  | *t*-statistic $\sqrt{n}\bar{x}/s$ | largest root $u_1$ of $A, B$ |
|---|---|---|
| Model: | $X_i \stackrel{\text{ind}}{\sim} N(\mu, \sigma^2)$ | $\binom{X_i}{Y_i} \sim N(0, \Sigma)$ |
|  | $H_0 : \mu = 0$ | $H_0 : \Sigma_{XY} = 0$ |
| Exact law: | $t \sim t_{n-1}$ | $u_1 \sim JOE_p(n-q-p, q-p)$ |
| Approx. law: | $\Phi(x) = \int_{-\infty}^x \phi(s)\, ds$ | $F_1(x) = \exp\{-\frac{1}{2}\int_x^\infty q(s) + (x-s)^2 q(s)\, ds\}$ |
| Convergence: | $O(n^{-1})$, not $O(n^{-1/2})$ | $O(p^{-2/3})$, not $O(p^{-1/3})$ |

*Convergence of quantiles.* Let $F_{p,1}$ denote the distribution function of $(W_p - \mu_p)/\sigma_p$—Theorem 1 asserts the convergence of $F_{p,1}(s)$ to $F_1(s)$ at rate $p^{-2/3}$. For given $\alpha \in (0,1)$, let $s_p(\alpha) = F_{p,1}^{-1}(\alpha)$ and $s(\alpha) = F_1^{-1}(\alpha)$ denote the $\alpha$th quantiles of $F_{p,1}$ and $F_1$, respectively. Under the assumptions of Theorem 1, we have convergence of the quantiles at rate $p^{-2/3}$: this follows from

$$f_{\min}(\alpha)|s_p(\alpha) - s(\alpha)| \leq |F_1(s_p(\alpha)) - F_1(s(\alpha))|$$
$$= |F_1(s_p(\alpha)) - F_{p,1}(s_p(\alpha))| \leq C(\alpha) p^{-2/3},$$

where $f_{\min}(\alpha)$ denotes the minimum value of $f_1(s)$ for values of $s$ between $s_p(\alpha)$ and $s(\alpha)$.

*Convergence of $\theta_{1,p}$.* An informal way of writing the conclusion of Theorem 1 is

$$W_p = \mu_p + \sigma_p Z_1 + O(p^{-4/3}),$$

as may be seen noting that $W_p$ and $Z_1$ can be defined on a common space using a $U(0,1)$ variate $U$, setting $W_p = \mu_p + \sigma_p F_{p,1}^{-1}(U)$ and $Z_1 = F_1^{-1}(U)$ and using the remark of the previous paragraph.

A straightforward delta-method argument now shows, for smooth functions $g(w)$, that

$$g(W_p) = g(\mu_p) + \sigma_p g'(\mu_p) Z_1 + O(p^{-4/3}).$$

In particular, with the logistic transformation $g(w) = e^w/(1 + e^w)$, we get, on the original scale,

$$\theta_{1,p} = \mu_\theta + \sigma_\theta Z_1 + O(p^{-4/3}),$$

with

$$\mu_\theta = \sin^2\left(\frac{\varphi + \gamma}{2}\right), \qquad \sigma_\theta^3 = \frac{\sin^4(\varphi + \gamma)}{4(m+n-1)^2 \sin\varphi \sin\gamma}.$$

If the $p^{-2/3}$ convergence rate is invariant to smooth transformations $g$, what is the special role of the logit transform in Theorem 1? Several related observations may be offered. First, empirical data analysis of quantities between 0 and 1, such as $\theta_i$, is often aided by the logit transform—indeed the idea to use it in this setting first emerged from efforts to improve the approximation in probability plots such as Figure 1. [Noting that $\theta_i$ may be thought of as squared canonical correlations, one also recalls that Fisher's $z = \tanh^{-1} r$ transformation improves the normal approximation in the case of a single coefficient.]

Second, at a technical level, as explained in Section 4, our operator convergence argument uses an integral representation of the Jacobi correlation kernel whose form is most similar to the Airy kernel when expressed using the hyperbolic tangent.

Finally, a geometric perspective: a natural metric on the cone of positive definite symmetric matrices is given by $d^2(A, B) = \sum \log^2 w_i$, where $\{w_i\}$ are the eigenvalues of $A^{-1}B$. Expressed in terms of the eigenvalues $\theta_i = w_i/(1 + w_i)$ of $(A+B)^{-1}B$, we get $d^2(A, B) = \sum \log^2[\theta_i/(1-\theta_i)]$, which is just Euclidean distance on the logit scale.

*Single Wishart limit.* If $\theta$ is an eigenvalue of $(A+B)^{-1}B$, then $m\theta/(1-\theta)$ is an eigenvalue of $mA^{-1}B$. Since $mA^{-1} \to I_p$ as $m \to \infty$, information about the largest eigenvalue of a single Wishart matrix is encoded in the double Wishart setting. This is spelled out in terms of hypergeometric functions by Koev (n.d.). However, as regards asymptotics, the preliminary limit $m \to \infty$ takes us out of the domain (3), and a separate treatment of the Tracy–Widom approximation is needed. Jiang (2008) does this in the Jacobi setting for $m^2/n \to \infty$ [and assuming $p/n \to c \in (0, \infty)$]. For the single Wishart matrix problem, in the complex case, see El Karoui (2006), and for the real setting, Ma (n.d.). The last two references both focus on second-order accuracy results.

*On the assumption that p is even.* The method of proof of Theorem 1 relies on a determinantal representation (48), that is, valid in the real case only for $p$ (written as $N + 1$ in the notation there) even. There is no such concern in the complex case.

Numerical investigation, both for this paper (Table 1) and its companion [Johnstone (2009)], suggests that the centering and scaling formulas (5) and Tracy–Widom approximation work as well for $p$ odd as for the $p$ even cases considered in the proofs. However, theoretical support for this observation remains incomplete. On the one hand, interlacing results would allow the largest eigenvalue for $p$ odd to be bracketed between settings with $p \pm 1$ [e.g., Chen (1971), Golub and Van Loan (1996), Corollary 8.6.3]. On the other hand, attempts to translate this directly to an $O(p^{-2/3})$ bound for the approximate distribution of $W_p$ in Theorem 1 encounter the following obstacle. Writing $\mu_p + \sigma_p s = \mu_{p+1} + \sigma_{p+1} s'$ so as to exploit the convergence result for $p+1$ leads to

$$s' - s = (\mu_{p+1} - \mu_p)/\sigma_p + (\sigma_{p+1}/\sigma_p - 1)s,$$

and calculations similar to those for Lemma 4 below show that $(\mu_{p+1} - \mu_p)/\sigma_p$ is generally $O(p^{-1/3})$.

**3. Jacobi ensembles.** We turn to a formulation of our results in the notation of random matrix theory (RMT), which provides tools and results needed for our proofs. A probability distribution on matrices is called an ensemble; that ensemble is termed unitary (resp., orthogonal) if the matrix elements are complex (real) and it is invariant under the action of the unitary (orthogonal) group. Because of the invariance, interest focuses on the joint density of the eigenvalues. A class of such ensembles of special interest in statistics has joint eigenvalue densities of the form

$$f_{N,\beta} = c_{N,\beta} \prod_{j<k} |x_j - x_k|^\beta \prod_{j=1}^{N} w_\beta(x_j).$$

The index $\beta = 1$ for real (orthogonal) ensembles and $\beta = 2$ for complex (unitary) ones. Here $w_\beta$ is one of the classical weight functions from the theory of orthogonal polynomials [Szegö (1967)], for which $\log w$ is a rational function with denominator degree $d \leq 2$. Most studied are the Gaussian ensembles ($d = 0$), with $w(x) = e^{-x^2/2}$, leading to Hermite polynomials, and corresponding to the eigenvalues of $N \times N$ Hermitian matrices with complex or real entries.

Next ($d = 1$) is the weight function $w(x) = e^{-x} x^\alpha$ of the Laguerre polynomials, corresponding to the eigenvalues of Gaussian covariance matrices, or equivalently to singular values of $N \times (N + \alpha)$ matrices with independent Gaussian real or complex entries.

Our interest in this paper lies with the final classical case ($d = 2$), with weight $w(x) = (1-x)^\alpha (1+x)^\beta$ leading to the Jacobi polynomials $P_N^{\alpha,\beta}(x)$. While the associated Jacobi unitary and orthogonal ensembles may have received relatively less attention in RMT, they may be seen as fundamental to the classical null hypothesis problems of multivariate statistical analysis.

REMARK. Second-order convergence results, with centering and scaling constants, are developed for the Laguerre case by El Karoui (2006) and Ma (n.d.) for the complex and real cases, respectively. A forthcoming manuscript will describe $N^{-2/3}$ convergence in the simplest Gaussian ensemble settings.

*Unitary case.* Along with the Gaussian and Laguerre ensembles, the eigenvalue density (17) of the Jacobi ensemble has the form $c_N \prod_i w(x_i) \cdot \Delta_N^2(x)$ with $w(x) = (1-x)^\alpha (1+x)^\beta$ being one of the classical weight functions of the theory of orthogonal polynomials, and

$$\Delta_N(x) = \prod_{i<j}(x_i - x_j) = \det[x_j^{k-1}]$$

being the Vandermonde determinant. Let $\phi_k(x)$ be the functions obtained by orthonormalizing the sequence $x^k w^{1/2}(x)$ in $L^2(-1,1)$. In fact,

(15) $$\phi_k(x) = h_k^{-1/2} w^{1/2}(x) P_k^{\alpha,\beta}(x),$$

where $P_k^{\alpha,\beta}(x)$ are the Jacobi polynomials, defined as in Szegö (1967). By a standard manipulation of the squared Vandermonde determinant, the eigenvalue density has a determinantal representation

$$f_{N,2}(x) = \frac{1}{N!} \det[S_{N,2}(x_j, x_k)]$$

with the correlation kernel having a Mercer expansion

(16) $$S_{N,2}(x,y) = \sum_{k=0}^{N-1} \phi_k(x)\phi_k(y).$$

The joint density of the eigenvalues is assumed to be

(17) $$f_{N,2}(x) = c \prod_{i=1}^{N}(1-x_i)^\alpha (1+x_i)^\beta \prod_{i<j}(x_i - x_j)^2.$$

With the identifications

(18) $$\begin{pmatrix} N \\ \alpha \\ \beta \end{pmatrix} = \begin{pmatrix} p \\ m-p \\ n-p \end{pmatrix}, \qquad \frac{1+x}{2} = \theta,$$

we recover the joint density of the roots of the double Wishart setting given at (10).

Our asymptotic model, which is equivalent to (3), assumes that $\alpha = \alpha(N)$ and $\beta = \beta(N)$ increase with $N$ in such a way that

(19) $$\frac{\alpha(N)}{N} \to a_\infty \in (0,\infty), \qquad \frac{\beta(N)}{N} \to b_\infty \in [0,\infty).$$

The dependence on $N$ will not always be shown explicitly. Introduce parameters:

(20) $$\kappa_N = \alpha + \beta + 2N + 1, \qquad \cos\varphi = \frac{\alpha-\beta}{\kappa_N}, \qquad \cos\gamma = \frac{\alpha+\beta}{\kappa_N}$$

and centering and scaling constants (on the $x$-scale):

$$(21) \qquad x_N = -\cos(\varphi + \gamma), \qquad \sigma_N^3 = \frac{2\sin^4(\varphi + \gamma)}{\kappa_N^2 \sin\varphi \sin\gamma}.$$

It will turn out that a better approximation is obtained working on the $u$-scale defined through the transformation $x = \tanh u$. The centering and scaling parameters become

$$(22) \qquad u_N = \tanh^{-1} x_N, \qquad \tau_N = \sigma_N/(1 - x_N^2).$$

The final centering and scaling parameters are suitable averages of those required for approximation at polynomial degree $N$ and $N-1$:

$$(23) \qquad \mu = \frac{\tau_N^{-1} u_N + \tau_{N-1}^{-1} u_{N-1}}{\tau_N^{-1} + \tau_{N-1}^{-1}}, \qquad \sigma^{-1} = \tfrac{1}{2}(\tau_N^{-1} + \tau_{N-1}^{-1}).$$

THEOREM 3. *There exist positive finite constants $c$ and $C(s_0)$ so that for $s \geq s_0$*

$$|P\{(\tanh^{-1} x_{(1)} - \mu)/\sigma \leq s\} - F_2(s)| \leq C N^{-2/3} e^{-cs}.$$

A consequence of our approach is a convergence result for the two-point correlation kernel, rescaled by $\tau(s) = \tanh(\mu + \sigma s)$, to the Airy kernel

$$(24) \qquad S_A(s, t) = \frac{\mathrm{Ai}(s)\,\mathrm{Ai}'(t) - \mathrm{Ai}(t)\,\mathrm{Ai}'(s)}{s - t}.$$

Indeed, uniformly on half intervals $[s_0, \infty)$, we show (in Section 7) that

$$(25) \qquad \sqrt{\tau'(s)\tau'(t)} S_{N,2}(\tau(s), \tau(t)) = S_A(s, t) + O(N^{-2/3} e^{-(s+t)/4}).$$

*Orthogonal case.* Suppose that $N+1$ is even. The joint density of the eigenvalues is assumed to be

$$(26) \qquad f(x) = c \prod_{i=1}^{N+1} (1 - x_i)^{(\alpha-1)/2} (1 + x_i)^{(\beta-1)/2} \prod_{i<j}^{N+1} |x_i - x_j|.$$

With the identifications

$$(27) \qquad \begin{pmatrix} N+1 \\ \alpha \\ \beta \end{pmatrix} = \begin{pmatrix} p \\ m-p \\ n-p \end{pmatrix}, \qquad \frac{1+x}{2} = \theta,$$

we recover the joint density of the roots of the real double Wishart setting given at (7). We match $N+1$ (rather than $N$) to $p$ because of a key formula relating the Jacobi orthogonal ensemble to the Jacobi unitary ensemble, (50) below.

The asymptotic model is the same as in the unitary case, that is, (19), as are the definitions of $(\kappa_N, \varphi, \gamma)$ in (20) and $(x_N, \sigma_N)$ in (21).

The final centering and scaling parameters are given by

$$(28) \qquad \mu = u_N, \qquad \sigma = \tau_N,$$

where $(u_N, \tau_N)$ are as in (22). Thus, after inserting (21) into (22),

$$(29) \qquad \mu = \log \tan\left(\frac{\varphi + \gamma}{2}\right), \qquad \sigma^3 = \frac{2}{\kappa_N^2 \sin^2(\varphi + \gamma) \sin\varphi \sin\gamma}.$$

THEOREM 4. *With $\mu, \sigma$ defined by (28), there exist positive finite constants $c$ and $C$ so that for $s \geq s_L$*

$$|P\{(\tanh^{-1} x_{(1)} - \mu)/\sigma \leq s\} - F_1(s)| \leq CN^{-2/3} e^{-cs}.$$

*Related work.* The first asymptotic analyses [e.g., Nagao and Wadati (1993), Nagao and Forrester (1995)] of the Jacobi correlation kernel assumed $\alpha, \beta$ fixed as $N \to \infty$. As a result, the upper limit of the $N$ eigenvalues equalled the upper limit of the base interval $[-1, 1]$, a "hard" edge.

The "double scaling limit" (19), natural for statistical purposes, has also arisen recently in RMT. Baik et al. (2006) [see also Deift (2007)] develop a probabilistic model leading to JUE for the celebrated observations of Krbalek and Seba (2000) that the bus spacing distribution in Cuernavaca, Mexico is well modeled by GUE. Baik et al. (2006) consider the double scaling limit (19) in the bulk.

Turning to the edge, Collins (2005) has shown that the centered and scaled distribution of the largest eigenvalue (in fact eigenvalues) of JUE converge to the Tracy–Widom distribution $F_2$ under asymptotic model (19) in the "ultraspherical" case in which $\alpha(N) = \beta(N)$. Our Theorem 3 applies also when $\alpha(N) \neq \beta(N)$ and provides an $O(N^{-2/3})$ rate bound. Collins uses a somewhat different centering and scaling, and with those proves convergence of the two-point correlation kernel with error $O(N^{-1/3+\varepsilon})$.

We remark that the universal Airy scaling limit arises in the double scaling limit because (19) forces the upper edge to be "soft," converging to $x_\infty < 1$.

**4. Strategy of proof.** A kernel $A(x, y)$ defines an operator $A$ on functions $g$ as usual via $(Ag)(y) = \int A(x, y)g(y)\,dy$. For suitable functions $f$, denote by $Sf$ the operator with kernel $S(x, y)f(y)$. Let $E_N$ denote expectation with respect to the density function (17). A key formula for unitary ensembles [e.g., Tracy and Widom (1998)], valid in particular for (17), states that

$$(30) \qquad E_N \prod_{j=1}^{N}[1 + f(x_j)] = \det(I + S_{N,2}f),$$

where the right-hand side is a Fredholm determinant of the operator $S_{N,2}f$ [Riesz and Sz.-Nagy (1955), Gohberg and Krein (1969), Chapter 4]. The choice $f = -\chi_0$, where $\chi_0(x) = I_{(x_0,1]}(x)$, yields the determinantal expression for the distribution of $x_{(1)}$:

$$F_{N2}(x_0) = P\left\{\max_{1 \leq j \leq N} x_j \leq x_0\right\} = \det(I - S_{N,2}\chi_0). \tag{31}$$

Tracy and Widom (1994) showed that the distribution $F_2$ has a determinantal representation

$$F_2(s_0) = \det(I - S_A),$$

where $S_A$ denotes the Airy operator on $L^2(s_0, \infty)$ with kernel (24). We introduce a rescaling $x = \tau(s)$, with $x_0 = \tau(s_0)$. To derive bounds on the convergence of $F_{N,2}(x_0)$ to $F_2(s_0)$, we use a bound due to Seiler and Simon (1975):

$$|\det(I - S_\tau) - \det(I - S_A)| \leq \|S_\tau - S_A\|_1 \exp(\|S_\tau\|_1 + \|S_A\|_1 + 1). \tag{32}$$

Here, operator $S_\tau$ has kernel

$$S_\tau(s,t) = \sqrt{\tau'(s)\tau'(t)} S_N(\tau(s), \tau(t)) \tag{33}$$

and is a suitably transformed, centered and scaled version of $S_{N,2}$ and $\|\cdot\|_1$ denotes trace class norm on operators on $L^2(s_0, \infty)$. The role of the nonlinear transformation contained within $\tau$ will be discussed further below. This bound reduces the convergence question to study of convergence of the kernel $S_\tau(x,y)$ to $S_A(x,y)$. For this, we use integral representations of both kernels. For the Airy kernel [Tracy and Widom (1994)]

$$S_A(s,t) = \int_0^\infty \operatorname{Ai}(s+z) \operatorname{Ai}(t+z) \, dz, \tag{34}$$

while for the Jacobi kernel, we use a formula to be found in Forrester (2004), Chapter 4. To state it, define

$$\hat{\phi}_N(u) = \frac{\phi_N(\tanh u)}{\cosh u}, \qquad \hat{S}_{N,2}(u,v) = \frac{S_N(\tanh u, \tanh v)}{\cosh u \cosh v}. \tag{35}$$

Then, from the final display in the proof of Forrester's Proposition 4.11,

$$\hat{S}_{N,2}(u,v) = \frac{(\kappa_N - 1)a_N}{2} \int_0^\infty [\hat{\phi}_N(u+w)\hat{\phi}_{N-1}(v+w) \\ + \hat{\phi}_{N-1}(u+w)\hat{\phi}_N(v+w)] \, dw. \tag{36}$$

The convergence argument will therefore be based on bounding the convergence of a suitably transformed, centered and scaled version of the weighted Jacobi polynomials $\phi_N(x)$ and $\phi_{N-1}(x)$ to the Airy function $\operatorname{Ai}(s)$.

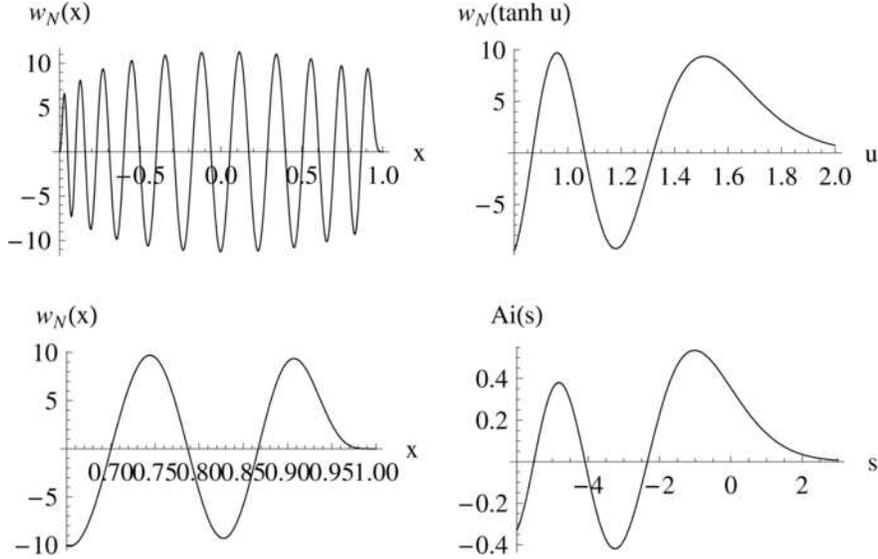

FIG. 2. *Top left: weighted Jacobi polynomial $w_N(x) = (1-x^2)^{1/2} w^{1/2}(x) P_N^{\alpha,\beta}(x)$ for $\alpha = 10, \beta = 5, N = 20$. Bottom left: focus on $w_N(x)$ in neighborhood of largest zero. Top right: nonlinear transformation of abscissa: $w_N(\tanh u)$. Bottom right: limiting Airy function* Ai(s). *Note improvement in approximation due to stretching of abscissa by hyperbolic tangent.*

The strategy for approximation by Ai($s$) is shown in Figure 2. The weighted Jacobi polynomial

$$(37) \qquad w_N(x) = (1-x)^{(\alpha+1)/2}(1+x)^{(\beta+1)/2} P_N^{\alpha,\beta}(x)$$

has $N$ zeros in $(-1, 1)$, shown in the top left panel. Zooming into a neighborhood of the largest zero (bottom left panel) shows a similarity with the graph of the Airy function. The nonlinear transformation $x = \tanh u$ of the abscissa is suggested by the form of the integral representation (35)–(36); in particular of course, it stretches $x \in (-1, 1)$ to $u \in (-\infty, \infty)$. The top right panel shows $w_N(\tanh u) = h_N^{1/2} \hat{\phi}_N(u)$: the stretching of the abscissa has improved the visual approximation to the Airy function, especially in the right tail.

To carry out the Jacobi polynomial asymptotics, several approaches are available, including saddle-point methods based on a contour integral representation [e.g., Wong and Zhao (2004)] and Riemann–Hilbert methods [e.g., Kuijlaars et al. (2004)]. Our situation is nonstandard because our model supposes that the parameters $\alpha(N), \beta(N)$ increase proportionally with $N$. We use the Liouville–Green approach set out in Olver (1974), since it comes with ready-made bounds for the error of approximation which are of great use in this paper. The Liouville–Green approximation relies on the fact that Jacobi polynomials and hence the function $w_N$ satisfy a second-order differential equation, (71) below, which may be put into the form

$$(38) \qquad w''(x) = \{\kappa^2 f(x) + g(x)\} w(x),$$

where $\kappa = 2N + \alpha + \beta + 1$ is the large parameter, and

$$(39) \qquad f(x) = \frac{(x - x_{N-})(x - x_{N+})}{4(1 - x^2)^2}, \qquad g(x) = -\frac{3 + x^2}{4(1 - x^2)^2}.$$

The values $x_{N-}$ and $x_{N+}$, given precisely at (75) below, are approximately the locations of the smallest and largest zeros of $P_N^{\alpha,\beta}$, respectively. They are the turning points of the differential equation; for example, $w_N(x)$ passes from oscillation to exponentially fast decay as $x$ moves through $x_{N+}$.

The Liouville–Green transformation is defined by ignoring $g(x)$ in (38) and transforming the independent variable $x$ into $\zeta$ via the equation $\zeta^{1/2} d\zeta = f^{1/2}(x) dx$, or equivalently

$$(40) \qquad (2/3)\zeta^{3/2} = \int_{x_+}^{x} f^{1/2}(x') dx'.$$

Then $W = (d\zeta/dx)^{1/2} w$ is close to satisfying the equation $d^2\bar{W}/d\zeta^2 = \kappa^2 \zeta \bar{W}$, which is a scaled form of the Airy differential equation, and so has linearly independent solutions in terms of Airy functions, traditionally denoted by $\text{Ai}(\kappa^{2/3}\zeta)$ and $\text{Bi}(\kappa^{2/3}\zeta)$. In fact, it turns out that

$$w_N(x) \doteq c_N (\dot{\zeta}(x))^{-1/2} \text{Ai}(\kappa^{2/3}\zeta(x)).$$

The value of the constant $c_N$ is fixed by matching the behavior of both sides as $x \to 1$ (Section A.4). For an approximation near $x_N = x_{N+}$, we introduce a new scaling $x = x_N + \sigma_N s_N$. To fix the scale $\sigma_N$, we linearize $\zeta(x_N + \sigma_N s_N)$ about its zero at $x_N$ and choose $\sigma_N$ so that $\kappa^{2/3}\zeta(x) \doteq s_N$. The resulting $\sigma_N$ is of order $N^{-2/3}$. To summarize the results of this local approximation and matching, define a particular multiple of $w_N(x)$, namely

$$(41) \qquad \check{\phi}_N(x) = (1 - x^2)^{1/2} \phi_N(x) / \sqrt{\kappa_N \sigma_N}.$$

Use of the Liouville–Green error bounds (Section 6.1) establishes that for $s_L \le s_N \le C N^{1/6}$,

$$(42) \qquad \check{\phi}_N(x_N + s_N \sigma_N) = \text{Ai}(s_N) + O(N^{-2/3} e^{-s_N/2}).$$

*u-scale.* Consistent with the top right panel of Figure 2, we need a translation of this approximation to $\bar{\phi}_N(u) = \check{\phi}_N(\tanh u)$. The $u$-scale centering and scaling $(u_N, \tau_N)$ are found by matching $\tanh(u_N + \tau_N t) \doteq x_N + \sigma_N t$ to first-order, and yield, for $t_0 \le t_N \le C N^{1/6}$

$$(43) \qquad \bar{\phi}_N(u_N + \tau_N t_N) = \text{Ai}(t_N) + O(N^{-2/3} e^{-t_N/2}).$$

In an entirely parallel development, there is an analogous result at degree $N-1$, again for $t_0 \leq t_{N-1} \leq CN^{1/6}$,

$$(44) \qquad \bar{\phi}_{N-1}(u_{N-1} + \tau_{N-1} t_{N-1}) = \mathrm{Ai}(t_{N-1}) + O(N^{-2/3} e^{-t_{N-1}/2}).$$

*s-scale.* A final calibration of the variables $t_N$ and $t_{N-1}$ is needed to match with the variable $s$ in the Airy function scale. Letting $\mu$ and $\sigma$ denote the centering and scaling for this calibration, yet to be determined, we set

$$(45) \qquad \phi_\tau(s) = \bar{\phi}_N(\mu + \sigma s), \qquad \psi_\tau(s) = \bar{\phi}_{N-1}(\mu + \sigma s).$$

The change of variables $u = \mu + \sigma s, v = \mu + \sigma t, w = \sigma z$ in (36) leads to

$$(46) \qquad S_\tau(s,t) = \frac{e_N}{2} \int_0^\infty [\phi_\tau(s+z)\psi_\tau(t+z) + \psi_\tau(s+z)\phi_\tau(t+z)].$$

The coefficient $e_N = 1 + O(N^{-1})$, as is shown at (138)–(140) below.

The expressions (43), (44) indicate that $O(N^{-2/3})$ error is only attainable for both $\phi_\tau$ and $\psi_\tau$ using separate scalings $t_N$ and $t_{N-1}$. The choice of $\mu$ and $\sigma$ can be made to transfer that $N^{-2/3}$ rate to a particular linear combination of $\phi_\tau$ and $\psi_\tau$: in the complex case, the bound

$$|\phi_\tau(s) + \psi_\tau(s)| \leq CN^{-2/3} e^{-s/4}$$

is convenient for achieving an $N^{-2/3}$ approximation of (46), and indeed, this forces the particular choices of $\mu$ and $\sigma$ in (144) and hence in Theorem 2. It is important for convergence of the integral in (46) that the above bound be global—valid on the right half line—and thus extending the "local" results of (43) and (44). This argument is set out in detail in Section 7.

REMARKS. (a) The function $\sqrt{f}$ appearing in the Liouville–Green transform

$$\sqrt{f}(x) = \frac{\sqrt{(x-x_-)(x-x_+)}}{2(1-x^2)}$$

is the same as the limiting bulk density of the eigenvalues found by Wachter (1980). The same phenomenon occurs in the single Wishart case.

(b) Our approximations are centered around the turning point of differential equation (38), which occurs at $s = 0$ in its Airy limit. The quantile $s = 0$ occurs beyond the upper quartile of the Tracy–Widom $F_1$ distribution, and so it is perhaps not surprising that the numerical quality of the Tracy–Widom approximation in Table 1 is better in the right tail of the distribution. It is a fortunate coincidence that precisely the right tail is the one of primary interest in statistical application.

*Orthogonal case.* A determinant representation for $\beta = 1$ analogous to (30) was developed by Dyson (1970). Tracy and Widom (1998) give a self-contained derivation of the formula

$$\text{(47)} \qquad E \prod_{j=1}^{N+1} [1 + f(x_j)] = \sqrt{\det(I + K_{N+1} f)},$$

with its immediate consequence, for $f = \chi_0 = -I_{(x_0, 1)}$,

$$\text{(48)} \qquad P\left\{ \max_{1 \leq k \leq N+1} x_k \leq x_0 \right\} = \sqrt{\det(I - K_{N+1} \chi_0)}.$$

Here we must assume that $N + 1$ is even, and then $K_{N+1}$ is a $2 \times 2$ matrix-valued operator whose kernel has the structure

$$\text{(49)} \qquad K_{N+1}(x, y) = \begin{pmatrix} I & -\partial_2 \\ \varepsilon_1 & T \end{pmatrix} S_{N+1,1}(x, y) - \begin{pmatrix} 0 & 0 \\ \varepsilon(x - y) & 0 \end{pmatrix}$$

Here $\partial_2$ denotes the operator of partial differentiation with respect to the second variable, and $\varepsilon_1$ the operator of convolution in the first variable with the function $\varepsilon(x) = \frac{1}{2} \text{sgn}(x)$. Thus $(\varepsilon S)(x, y) = \int \varepsilon(x - u) S(u, y) \, du$. Finally $T$ denotes transposition of variables $TS(x, y) = S(y, x)$.

The derivation of Tracy and Widom (1998) does not completely determine $S_{N+1,1}(x, y)$. More explicit expressions, developed by Dyson (1970) and Mahoux and Mehta (1991), use families of skew-orthogonal polynomials. Adler et al. (2000) relate these skew polynomials to the orthogonal polynomials occurring in the unitary case. A key observation is that the weight function of the orthogonal ensemble should be suitably perturbed from that of the corresponding unitary ensemble. This leads Adler et al. (2000) to a formula that is central for this paper:

$$\text{(50)} \quad S_{N+1,1}(x, y) = \left( \frac{1 - y^2}{1 - x^2} \right)^{1/2} S_{N,2}(x, y) + a_N \tilde{\kappa}_N \tilde{\phi}_N(x) (\varepsilon \tilde{\phi}_{N-1})(y).$$

Here $S_{N,2}$ is the unitary kernel (16) associated with the Jacobi unitary ensemble (17), and

$$\text{(51)} \qquad \tilde{\kappa}_N = (2N + \alpha + \beta)/2,$$

$$\text{(52)} \qquad \tilde{\phi}_N(x) = \phi_N(x)/\sqrt{1 - x^2}.$$

The orthogonal kernel is thus expressed in terms of the unitary kernel and a rank-1 remainder term. The formula allows convergence results from the unitary case to be reused, with relatively minor modification.

As regards the limit, Tracy and Widom (2005) showed that

$$F_1(s_0) = \sqrt{\det(I - K_{GOE})},$$

where, for the purposes of this paper, the GOE kernel may be written

$$K_{GOE}(s,t) = \begin{bmatrix} I & -\partial_2 \\ -\tilde{\varepsilon}_1 & T \end{bmatrix} S(s,t)$$

(53)
$$+ \tfrac{1}{2} \begin{bmatrix} \mathrm{Ai}(s) & 0 \\ \tilde{\varepsilon}(\mathrm{Ai})(t) - \tilde{\varepsilon}(\mathrm{Ai})(s) & \mathrm{Ai}(t) \end{bmatrix} + \begin{bmatrix} 0 & 0 \\ -\varepsilon(s-t) & 0 \end{bmatrix}$$

with

$$S(s,t) = S_A(s,t) - \tfrac{1}{2}\mathrm{Ai}(s)\tilde{\varepsilon}(\mathrm{Ai})(t).$$

Here $(\tilde{\varepsilon}f)(s) = \int_s^\infty f$ and $(\tilde{\varepsilon}_1 S)(s,t) = \int_s^\infty S(u,t)\,du$. We use $\tilde{\varepsilon}_1$ in place of $\varepsilon$ in (53) because convergence to $\mathrm{Ai}(s)$ is stable in the right-hand side, but oscillatory and difficult to handle in the left tail.

We bound the convergence of $F_{N+1,1}(s_0)$ to $F_1(s_0)$ via an analog of (32):

(54)
$$|F_{N+1}(s_0) - F_1(s_0)|$$
$$\leq C(s_0) C(K_\tau, K_{GOE})$$
$$\times \left\{ \sum_i \|K_{\tau,ii} - K_{GOE,ii}\|_1 + \sum_{i \neq j} \|K_{\tau,ij} - K_{GOE,ij}\|_2 \right\}.$$

Here, in analogy with (33),

$$K_\tau(s,t) = \sqrt{\tau'(s)\tau'(t)} K_{N+1}(\tau(s), \tau(t)).$$

The detailed work of representing $K_\tau - K_{GOE}$ in terms of the transformation, centering and scaling implicit in $\tau$ is done in Section 8.3.

As noted by Tracy and Widom (2005), a complication arises in the orthogonal case: as so far described, $K_\tau$ is not a trace class operator, as would be required properly to define the (Fredholm) determinant. This obstacle is evaded by regarding $K_\tau$ as a matrix Hilbert–Schmidt operator on $L^2(\rho) \oplus L^2(\rho^{-1})$ where $L^2(\rho^\pm)$ are weighted Hilbert spaces $L^2([s_0,\infty), \rho^\pm(s)\,ds)$. We assume at least that $\rho^{-1} \in L^1$, and so, in particular, it follows that $\varepsilon \colon L^2(\rho) \to L^2(\rho^{-1})$. Section 8.2 has more detail on this.

A few further remarks on the origin of the $N^{-2/3}$ rate in the orthogonal case. In the unitary case, the $N^{-2/3}$ rate of convergence for the kernel $S_\tau(s,t)$ was obtained by a calculated trade-off of centering and scaling in the approximations for $\phi_\tau$ and $\psi_\tau$, at degrees $N$ and $N-1$, respectively. In the orthogonal case, a cancellation of $N^{-1/3}$ terms, somewhat fortuitous and unexplained, occurs between the integral and rank-1 terms in (50), so that such a calculated trade-off is not required. More specifically, we reuse the unitary case approximations to $\phi_\tau$ and $\psi_\tau$, but now with the straightforward choices $\mu = u_N, \sigma = \tau_N$ in (45). In Section 7.3 it is shown that

(55) $$|\phi_\tau(s) - \mathrm{Ai}(s)| \leq C N^{-2/3} e^{-s/4},$$

(56) $$|\psi_\tau(s) - \mathrm{Ai}(s) - \Delta_N \mathrm{Ai}'(s)| \leq C N^{-2/3} e^{-s/4},$$

where $\Delta_N = (u_N - u_{N-1})/\tau_{N-1} = O(N^{-1/3})$. Thus, to obtain the $N^{-2/3}$ rate for $\psi_\tau$ here, it is necessary to retain the derivative term, itself of order $N^{-1/3}$.

Focus on the $(1,1)$ entries of rescaled $K_{N+1}$ and its limit $K_{GOE}$. In Section 8.3, it is shown that (50) may be written

$$\tau'(s)S_{N+1,1}(\tau(s),\tau(t)) = e_N[\bar{S}_\tau(s,t) + \tfrac{1}{2}\phi_\tau(s)\varepsilon\psi_\tau(t)].$$

In contrast with the unitary case, the use of (55) and (56) leads to an $N^{-1/3}$ term:

$$\bar{S}_\tau(s,t) = S_A(s,t) - \frac{\Delta_N}{2}\operatorname{Ai}(s)\operatorname{Ai}(t) + O(N^{-2/3}).$$

Turning to the rank-1 term and again using (55) and (56),

$$\frac{1}{2}\phi_\tau(s)(\varepsilon\psi_\tau)(t) = \frac{1}{2}\operatorname{Ai}(s)[1 - \tilde{\varepsilon}(\operatorname{Ai})(t)] + \frac{\Delta_N}{2}\operatorname{Ai}(s)\operatorname{Ai}(t) + O(N^{-2/3}),$$

and so, remarkably, the $N^{-1/3}$ terms in the previous two displays cancel, yielding $N^{-2/3}$ convergence, at least for the $(1,1)$ entry. The remainder of the convergence argument—including the operator norm bounds to give the $e^{-s_0/2}$ dependence—may be found in Section 8.4.

*Relation to other work.* There is a large literature on the asymptotic behavior of Jacobi polynomials as $N \to \infty$. For fixed $\alpha$ and $\beta$, classical results are given in Szegö (1967); for more recent results [see, e.g., Kuijlaars et al. (2004) and Wong and Zhao (2004)]. There is a smaller, but growing, literature on results when $\alpha$ and $\beta$ depend on $N$ and tend to infinity with $N$ [see, e.g., Chen and Ismail (1991) and Bosbach and Gawronski (1999)]. Closer to our approach is Dunster (1999), who uses Liouville–Green transformations to study ultraspherical polynomials (a subclass of Jacobi polynomials) with $\alpha = \beta$ proportional to $N$, and provide approximations in terms of Whittaker functions. Carteret, Ismail and Richmond (2003) give Airy approximations to Jacobi polynomials with one of the parameters proportional to $N$. Collins (2005), Lemmas 4.12 and 4.14, provides Airy approximations similar to (42), but with error term $O(N^{-2/3+\varepsilon})$; his proof uses the differential equation satisfied by (37), but not the specific Liouville–Green method adopted here.

**5. Jacobi polynomials; preliminaries.** We collect here some useful facts [Szegö (1967), Chapter 4] about the Jacobi polynomials $P_N = P_N^{\alpha,\beta}(x)$. They are orthogonal with respect to the weight function $w(x) = (1-x)^\alpha(1+x)^\beta$ on $[-1,1]$, and have $L_2$ norms:

$$(57) \quad h_N = \int_{-1}^1 P_N^2(x)w(x)\,dx = \frac{2^{\alpha+\beta+1}}{2N+\alpha+\beta+1}\frac{\Gamma(N+\alpha+1)\Gamma(N+\beta+1)}{\Gamma(N+1)\Gamma(N+\alpha+\beta+1)}.$$

The leading coefficient $P_N(x) = l_n x^N + \cdots$ is

$$l_N = 2^{-N} \binom{2N + \alpha + \beta}{N} \tag{58}$$

and the value at $x = 1$ is

$$P_N^{\alpha,\beta}(1) = \binom{N + \alpha}{N}. \tag{59}$$

The Christoffel–Darboux formula states that

$$S_{N,2}(x, y) = \sum_{k=0}^{N-1} \phi_k(x)\phi_k(y) = a_N \frac{\phi_N(x)\phi_{N-1}(y) - \phi_{N-1}(x)\phi_N(y)}{x - y}, \tag{60}$$

$$a_N = \left(\frac{h_N}{h_{N-1}}\right)^{1/2} \frac{l_{N-1}}{l_N}. \tag{61}$$

5.1. *Parameterizations.* We collect and connect several equivalent parameter sets which are each useful at certain points.

(a) *Statistics parameters* $(p, m, n)$. These describe the parameters of the two Wishart distributions described above. We are adopting the notation of Mardia, Kent and Bibby (1979), who interpret the parameter $p$ as "dimension," $m$ as "error" degrees of freedom and $n$ as "hypothesis" degrees of freedom.

(b) *Jacobi parameters* $(N, \alpha, \beta)$. These are the parameters appearing in the conventional Jacobi polynomials [Szegö (1967), Chapter IV]. The connection to the Wishart matrix parameters is given in the complex case by

$$\begin{pmatrix} N \\ \alpha \\ \beta \end{pmatrix} = \begin{pmatrix} p \\ m - p \\ n - p \end{pmatrix}. \tag{62}$$

[For the real case, see (18) below.] The conditions $m, n \geq p$ correspond to $\alpha, \beta \geq 0$.

(c) *Liouville–Green form* $(\kappa, \lambda, \mu)$. To describe compactly the form (87) of the Jacobi differential equation below, introduce

$$\begin{aligned} \kappa &= 2N + \alpha + \beta + 1, \\ \lambda &= \alpha/\kappa \geq 0, \\ \mu &= \beta/\kappa \geq 0. \end{aligned} \tag{63}$$

(c′) A variant is $(\kappa, a, b)$, where on setting $N_+ = N + \frac{1}{2}$,

$$\alpha = N_+ a, \qquad \beta = N_+ b, \tag{64}$$

which implies

$$\lambda = a/(2 + a + b), \qquad \mu = b/(2 + a + b). \tag{65}$$

(d) *Trigonometric forms* $(\gamma, \varphi)$.

$$\cos \gamma = \lambda + \mu, \qquad \cos \varphi = \lambda - \mu. \tag{66}$$

From the definitions of $(\kappa, \lambda, q)$, we deduce the ranges

$$0 < \gamma \leq \pi/2, \qquad 0 < \varphi < \pi, \qquad \gamma \leq \varphi. \tag{67}$$

The last four systems are related by the equalities

$$\cos \gamma = \lambda + \mu = \frac{a+b}{a+b+2} = \frac{\alpha+\beta}{\alpha+\beta+2N_+}, \tag{68}$$

$$\cos \varphi = \lambda - \mu = \frac{a-b}{a+b+2} = \frac{\alpha-\beta}{\alpha+\beta+2N_+}. \tag{69}$$

(e) *Half angle forms* $(\varphi/2, \gamma/2)$.

$$\sin^2(\gamma/2) = \frac{p+1/2}{m+n+1}, \qquad \sin^2(\varphi/2) = \frac{n+1/2}{m+n+1},$$

as may be seen by using $\cos \gamma = 1 - 2\sin^2(\gamma/2)$ on the left-hand side and relations (18) on the right-hand side of (68) and (69).

5.2. *Differential equation for Jacobi polynomials.* The weighted Jacobi polynomial

$$w_N(x) = (1-x)^{(\alpha+1)/2}(1+x)^{(\beta+1)/2} P_N^{\alpha,\beta}(x) \tag{70}$$

satisfies a second-order equation without first-order term [Szegö (1967), (4.24.1)]

$$w''(x) = q(x)w(x) = \frac{\mathfrak{n}(x)}{4(1-x^2)^2} w(x), \tag{71}$$

where the quadratic polynomial

$$\begin{aligned}\mathfrak{n}(x) &= (\alpha^2-1)(x+1)^2 + (\beta^2-1)(x-1)^2 \\ &\quad + [4N(N+\alpha+\beta+1) + 2(\alpha+1)(\beta+1)](x^2-1).\end{aligned} \tag{72}$$

This equation may be put into a form suitable for asymptotics, namely

$$w''(x) = \{\kappa^2 f(x) + g(x)\} w(x), \tag{73}$$

by using the Liouville–Green parameters (63) to set

$$f(x) = \frac{x^2 + 2(\lambda^2 - \mu^2)x + 2\lambda^2 + 2\mu^2 - 1}{4(1-x^2)^2}, \qquad g(x) = -\frac{3+x^2}{4(1-x^2)^2}. \tag{74}$$

These choices for $f(x)$, $g(x)$ and $\kappa$ [taken from Dunster (1999), (4.1)] are not unique; however, our goal of obtaining approximations to $w_N(x)$ with error bounds $O(1/\kappa)$ imposes constraints which lead naturally to this choice (Section A.2 has some details).

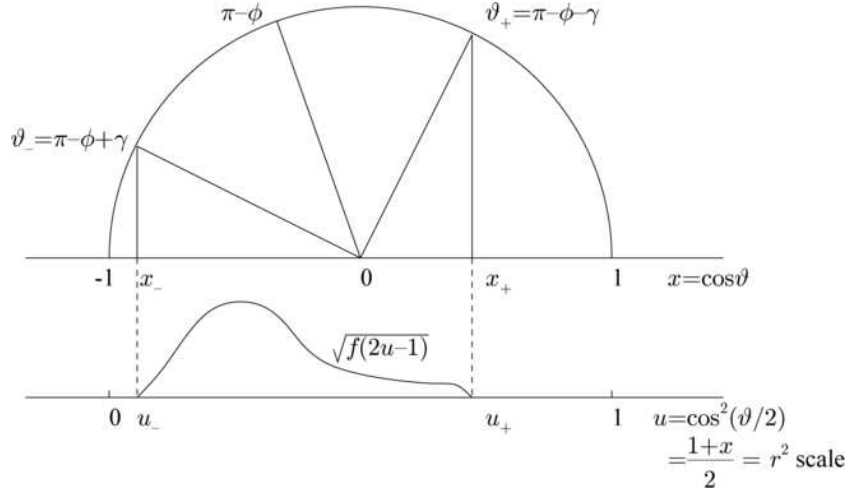

Fig. 3. *Relationship of some key parameters. The turning points $x_+$ and $x_-$ are shown on both the "Jacobi" scale $x \in [-1, 1]$ and the "squared correlation" scale $u = r^2 \in [0, 1]$, along with the respective angle and half-angle interpretations.*

The turning points of the differential equation are given by the zeros of $f$, namely

$$x_\pm = \mu^2 - \lambda^2 \pm \sqrt{\{1-(\lambda+\mu)^2\}\{1-(\lambda-\mu)^2\}} \tag{75}$$

$$= -\cos(\varphi \pm \gamma) = \cos(\pi - \varphi \mp \gamma), \tag{76}$$

where we used (68) and (69). Where necessary to show the dependence on $N$, we write $x_{N\pm}$.

It is easily verified that $\vartheta_\pm = \pi - (\varphi \pm \gamma) \in [0, \pi]$. Indeed (67) entails that $\vartheta_- \leq \pi$, while $\lambda \geq 0$ implies via (68) and (69) that $\cos\varphi \geq -\cos\gamma = \cos(\pi - \gamma)$, and hence that $\varphi \leq \pi - \gamma$ and so $\vartheta_+ \geq 0$. In particular

$$-1 \leq x_- < x_+ \leq 1 \tag{77}$$

and

$$x_+ - x_- = 2\sin\varphi\sin\gamma. \tag{78}$$

The situation is summarized in Figure 3.

Furthermore, a little algebra with (75) shows that both

$$x_+ < 1 \quad \Leftrightarrow \quad \lambda > 0 \quad \Leftrightarrow \quad m > p, \tag{79}$$

$$x_- > -1 \quad \Leftrightarrow \quad \mu > 0 \quad \Leftrightarrow \quad n > p. \tag{80}$$

For later reference, note that we may rewrite $f$ as

$$f(x) = (x - x_{N+})k(x), \qquad k(x) = \frac{x - x_{N-}}{4(1-x^2)^2}. \tag{81}$$

In particular, by combining (78) and (76), we have

$$k_N := \frac{x_{N+} - x_{N-}}{4(1-x_{N+}^2)^2} = \frac{\sin\varphi \sin\gamma}{2\sin^4(\varphi+\gamma)}. \tag{82}$$

5.3. *Asymptotic setting.* The asymptotic model used in this paper supposes that there is a sequence of eigenvalue distribution models, such as (17) and (26), indexed by $N$, the number of variables. The Jacobi parameters $\alpha = \alpha(N), \beta = \beta(N)$ are regarded as functions of $N$.

ASSUMPTION (A). In the following equivalent forms:
For $(N, \alpha, \beta)$:

$$\frac{\alpha(N)}{N} \to a_\infty \in (0,\infty), \qquad \frac{\beta(N)}{N} \to b_\infty \in [0,\infty). \tag{83}$$

For $(\kappa, \lambda, \mu)$:

$$\lambda_N \to \lambda_\infty, \qquad \mu_N \to \mu_\infty \quad \text{s.t.} \quad \lambda_\infty > 0, \qquad \lambda_\infty + \mu_\infty < 1. \tag{84}$$

For $(p, m, n)$:

$$\frac{p}{m+n} \to \rho_\infty > 0, \qquad \frac{m}{p} \to a_\infty + 1 > 1. \tag{85}$$

We emphasize two consequences of these assumptions. First, that $\lim x_{N+} < 1$, so that the right edge is "soft," that is, separated from the upper limit of support of the weight function $w(x)$. Indeed, from (75), $x_+ = 1$ if and only if $\lambda = 0$, and since $\lambda = a/(2+a+b)$, this is prevented because $a_\infty \in (0,\infty)$ and $b_\infty < \infty$.

Second, we have $\lim x_{N+} - x_{N-} > 0$, so that the two turning points are asymptotically separated. Indeed from (78), $x_+ > x_-$ if and only if $\gamma > 0$, and from (68) this occurs if $a_\infty$ and $b_\infty$ are both finite.

The constants in our asymptotic bounds (such as Theorems 3 and 4) will depend on the limiting values $a_\infty, b_\infty$, or their equivalent forms in (84) or (85). These dependencies will not be worked out in detail; instead we will use the following somewhat less precise approach.

Introduce the sets

$$D = \{(\lambda, \mu) : \lambda \geq 0, \mu \geq 0, \lambda + \mu \leq 1\},$$

$$D_\delta = \{(\lambda, \mu) \in D : \lambda \geq \delta, \lambda + \mu \leq 1 - \delta\}.$$

Our results will be valid with $C = C(\delta)$ for all $N$ such that $(\lambda_N, \mu_N) \in D_\delta$.

Under the asymptotic assumptions (A), some straightforward simplifications in formulas occur:

LEMMA 1. *The coefficient $a_N$ defined at (61) satisfies*

(86) $$a_N = (\tfrac{1}{2}\sin\varphi \sin\gamma)(1 + O(N^{-1})).$$

**6. Jacobi polynomial asymptotics near largest zero.** The goal of this section is to establish an Airy function approximation to a version of the weighted Jacobi polynomial $\phi_N(x)$ in a neighborhood of its largest zero, or more correctly, in a shrinking neighborhood about the upper turning point $x_N$ of (75). More precisely, in terms of the weighted version of $\phi_N$ defined at (41) and the scaling parameter $\sigma_N$ defined at (104) below, we develop local approximations of the form

$$\check{\phi}_N(x_N + s_N \sigma_N) = \operatorname{Ai}(s_N) + O(N^{-2/3} e^{-s_N/2}),$$

valid for $s_N \in [s_L, CN^{1/6}]$ and $N > N_0$ sufficiently large.

6.1. *Liouville–Green approach.* We begin with an overview of the Liouville–Green approach to be taken, following [Olver (1974), Chapter 11]. The classical orthogonal polynomials [such as Laguerre $L_N^\alpha(x)$, Jacobi $P_N^{\alpha,\beta}(x)$] satisfy a second-order linear ordinary differential equation which, if the polynomials are multiplied by a suitable weight function, may be put in the form

(87) $$\frac{d^2 w}{dx^2} = q(x)w(x) = \{\kappa^2 f(x) + g(x)\}w(x), \qquad x \in (a,b),$$

where $\kappa = \kappa(N)$ is a parameter, later taken as large. The precise decomposition of $q$ into $\kappa^2 f + g$ is made in order to obtain $O(1/N)$ error bounds below. A zero $x_*$ of $f$ is called a *turning point* because, as will be seen in our example, it separates an interval in which the solution is of exponential type from one in which the solution oscillates. We will assume, for some interval $(a,b)$ containing $x_*$, that $f(x)/(x-x_*)$ is positive and twice continuously differentiable and that $g(x)$ is continuous.

Define new independent and dependent variables $\zeta$ and $W$ via the equations

$$\zeta\left(\frac{d\zeta}{dx}\right)^2 = f(x), \qquad W = \left(\frac{d\zeta}{dx}\right)^{1/2} w.$$

These choices put (87) into the form

(88) $$\frac{d^2 W}{d\zeta^2} = \{\kappa^2 \zeta + \psi(\zeta)\}W,$$

where the perturbation term $\psi(\zeta) = \hat{f}^{-1/4}(d^2/d\zeta^2)(\hat{f}^{1/4}) + g/\hat{f}$. Here $\hat{f}$ is defined by

(89) $$\hat{f}(x) = (d\zeta/dx)^2 = f(x)/\zeta.$$

If the perturbation term $\psi(\zeta)$ in (88) were absent, the equation $d^2W/d\zeta^2 = \kappa^2\zeta W$ would have linearly independent solutions in terms of Airy functions, traditionally denoted by $\text{Ai}(\kappa^{2/3}\zeta)$ and $\text{Bi}(\kappa^{2/3}\zeta)$. Our interest is in approximating the *recessive* solution $\text{Ai}(\kappa^{2/3}\zeta)$, so write the relevant solution of (88) as $W_2(\zeta) = \text{Ai}(\kappa^{2/3}\zeta) + \eta(\zeta)$. In terms of the original dependent and independent variables $\xi$ and $w$, the solution $W_2$ becomes

$$(90) \qquad w_2(x,\kappa) = \hat{f}^{-1/4}(x)\{\text{Ai}(\kappa^{2/3}\zeta) + \varepsilon_2(x,\kappa)\}.$$

Olver [(1974), Theorem 11.3.1] provides an explicit bound for $\eta(\zeta)$ and hence $\varepsilon_2$ and its derivative in terms of the function $\mathcal{V}(\zeta) = \int_\zeta^{\zeta(b)} |\psi(v)v^{-1/2}|\,dv$. To describe these error bounds even in the oscillatory region of $\text{Ai}(x)$, Olver (1974) introduces a positive weight function $E(x) \geq 1$ and positive moduli functions $M(x) \leq 1$ and $N(x)$ such that

$$(91) \qquad \begin{aligned} |\text{Ai}(x)| &\leq M(x)E^{-1}(x) \\ |\text{Ai}'(x)| &\leq N(x)E^{-1}(x) \end{aligned} \qquad \text{for all } x.$$

[Here, $E^{-1}(x)$ denotes $1/E(x)$.] In addition,

$$(92) \qquad \text{Ai}(x) = \frac{1}{\sqrt{2}} M(x)E^{-1}(x), \qquad x \geq c \doteq -0.37$$

and the asymptotics as $x \to \infty$ are given by

$$(93) \quad E(x) \sim \sqrt{2}e^{(2/3)x^{3/2}}, \qquad M(x) \sim \pi^{-1/2}x^{-1/4}, \qquad N(x) \sim \pi^{-1/2}x^{1/4}.$$

The key bounds of Olver [(1974), Theorem 11.3.1] then state, for $x \in (a,b)$,

$$(94) \qquad |\varepsilon_2(x,\kappa)| \leq M(\kappa^{2/3}\zeta)E^{-1}(\kappa^{2/3}\zeta)\left[\exp\left\{\frac{\lambda_0}{\kappa}\mathcal{V}(\zeta)\right\} - 1\right],$$

$$(95) \quad |\partial_x \varepsilon_2(x,\kappa)| \leq \kappa^{2/3}\hat{f}^{1/2}(x)N(\kappa^{2/3}\zeta)E^{-1}(\kappa^{2/3}\zeta)\left[\exp\left\{\frac{\lambda_0}{\kappa}\mathcal{V}(\zeta)\right\} - 1\right],$$

where $\lambda_0 \doteq 1.04$. For $\kappa^{2/3}\zeta \geq c$, (92) shows that the coefficient in (94) is just $\sqrt{2}\,\text{Ai}(\kappa^{2/3}\zeta)$. Here

$$\mathcal{V}(\zeta) = \mathcal{V}(\zeta(x)) = \mathcal{V}_{[x,1]}(H) = \int_x^1 |H'(t)|\,dt$$

is the total variation on $[x,1]$ of the *error control* function

$$H(x) = -\int_0^{\zeta(x)} |v|^{-1/2}\psi(v)\,dv.$$

Section A.3 has more information on $\mathcal{V}(\zeta)$.

*Application to Jacobi polynomials.* In the case of Jacobi polynomials $P_N^{\alpha,\beta}(x)$, the points $x = \pm 1$ are regular singularities and the points $x_\pm$ defined by (75)–(77) are turning points. We are interested in behavior near the upper turning point $x_+$, which is located near the largest zero of $P_N^{\alpha,\beta}$. We apply the foregoing discussion to the interval $(a,b) = (x_0, 1)$, where $x_0 = \frac{1}{2}(x_+ + x_-) = \mu^2 - \lambda^2$. In particular, the independent variable $\zeta(x)$ is given in terms of $f(x)$ by

$$(96) \qquad (2/3)\zeta^{3/2} = \int_{x_+}^{x} f^{1/2}(x')\,dx'.$$

It is easily seen from (74) that $\zeta \to \infty$ as $x \to 1$. More precisely, it is shown at length in Section A.4, Proposition 4 that

$$(2/3)\zeta^{3/2}(x) = \frac{\alpha}{2\kappa}\log(1-x)^{-1} + c_{0N} + o(1),$$

where $c_{0N} = c_{0N}(a,b)$ is given at (230).

Bound (94) is valid only if the integral defining $\mathcal{V}(\zeta)$ converges as $\zeta \to \zeta(b) = \infty$. That this is true for the specific choices of $f$ and $g$ made in (74) follows from arguments in Olver (1974); see Section A.2 for further details.

REMARK A. For behavior near the lower turning point $x_-$, near the smallest zero of $P_N^{\alpha,\beta}$, we would consider instead the interval $(a, a_1) = (-1, -\cos\varphi)$, with a corresponding redefinition of $\zeta(x)$, and we would require convergence of $\int_{-\infty}^{\zeta} |\psi(v)v^{-1/2}|\,dv$. This would be relevant to approximation of the distribution of the smallest canonical correlation, although this can be handled simply through (14) in Section 2.

Bound (94) has a double asymptotic property in $x$ and $\kappa$ which will be useful. First, suppose that $N$ and hence $\kappa$ are held fixed. As $x \to 1$, $\mathcal{V}(\zeta) \to 0$ and so from (94) and its following remarks $\varepsilon_2(x, \kappa) = o(\mathrm{Ai}(\kappa^{2/3}\zeta))$. Consequently, as $x \to 1$

$$(97) \qquad w_2(x, \kappa) \sim \hat{f}^{-1/4}(x)\,\mathrm{Ai}(\kappa^{2/3}\zeta).$$

If the weighted polynomial $w_N(x)$ is a recessive solution of (87), then it must be proportional to $w_2$:

$$(98) \qquad w_2(x, \kappa) = c_N^{-1} w_N(x).$$

The important consequence is that $c_N$ may now be identified by comparing the growth of $w_N(x)$ as $x \to 1$ with that of $w_2(x, \kappa)$. In Section A.5 it is shown that

$$(99) \qquad c_N = e^{\theta''/N}\kappa_N^{1/6}h_N^{1/2},$$

where $\theta'' = O(1)$. Hence

$$(100) \qquad w_N(x) = e^{\theta''/N}\kappa_N^{1/6}h_N^{1/2}w_2(x, \kappa).$$

The second property is that bound (94) holds for all $x \in (a, b)$ and so in any interval $(a_1, b) \subset (a, b)$ upon which $\sup \mathcal{V}(\zeta) < \infty$, we have

$$(101) \qquad |\varepsilon_2(x, \kappa)| = O(1/\kappa) = O(1/N).$$

Comparing the (70) and (15) of $w_N$ and $\phi_N$, we obtain

$$(1-x^2)^{1/2}\phi_N(x) = h_N^{-1/2} w_N(x).$$

Combining (100) with the Airy approximation (90) to $w_2(x,\kappa)$, we obtain

$$(102) \quad (1-x^2)^{1/2}\phi_N(x) = e^{\theta''/N}\kappa_N^{1/6}\hat{f}^{-1/4}(x)\{\text{Ai}(\kappa^{2/3}\zeta) + \varepsilon_2(x,\kappa)\}.$$

*Local scaling parameter $\sigma_N$.* We are chiefly interested in values $x = x_N + \sigma_N s$ near the upper turning point. The scaling constant $\sigma_N$ is chosen so that for fixed $s$, as $N \to \infty$, $\text{Ai}(\kappa^{2/3}\zeta) \to \text{Ai}(s)$. Expand $\zeta(x)$ about the turning point $x_N$:

$$(103) \quad \kappa^{2/3}\zeta(x) = \kappa^{2/3}\zeta(x_N + \sigma_N s) = \kappa^{2/3}\sigma_N s \dot{\zeta}_N + \tfrac{1}{2}\kappa^{2/3}\sigma_N^2 s^2 \ddot{\zeta}_N + \cdots.$$

Setting the coefficient of $s$ to 1 yields

$$(104) \qquad \sigma_N = (\kappa^{2/3}\dot{\zeta}_N)^{-1}.$$

Consider now the coefficient of $\text{Ai}(\kappa^{2/3}\zeta)$ in (102): $e^{\theta''/N}\kappa_N^{1/6}\hat{f}^{-1/4}(x)$. Observe from (89) that $\hat{f}^{-1/4}(x) = \dot{\zeta}(x)^{-1/2}$, and then note that as $x \to x_N$, we have $\dot{\zeta}(x)^{-1/2} \to \dot{\zeta}_N^{-1/2} = \kappa_N^{1/3}\sigma_N^{1/2}$ from (104). Consequently

$$(105) \qquad \kappa_N^{1/6}\hat{f}^{-1/4}(x) = \sqrt{\sigma_N \kappa_N}(\dot{\zeta}(x)/\dot{\zeta}_N)^{-1/2}$$

and so finally we have

$$(106) \quad \check{\phi}_N(x) := (1-x^2)^{1/2}\frac{\phi_N(x)}{\sqrt{\kappa_N \sigma_N}} = \bar{e}_N r_N(x)[\text{Ai}(\kappa^{2/3}\zeta) + \varepsilon_2(x,\kappa)],$$

where $\bar{e}_N = e^{\theta''/N} = 1 + O(N^{-1})$ and

$$(107) \qquad r_N(x) := \frac{\kappa_N^{1/6}}{\sqrt{\kappa_N \sigma_N}}\hat{f}^{-1/4}(x) = \left(\frac{\dot{\zeta}_N(x)}{\dot{\zeta}_N}\right)^{-1/2},$$

where the second equality is just (105).

The goal toward which we are working is a uniform bound on the Airy approximation in a local (but growing) region about $x_N$.

PROPOSITION 1. *For $x = x_N + s_N \sigma_N$ and $s_L \leq s_N \leq CN^{1/6}$, we have*

$$(108) \qquad \check{\phi}_N(x) = \text{Ai}(s_N) + O(N^{-2/3}e^{-s_N/2}),$$

$$(109) \qquad \sigma_N \check{\phi}'_N(x) = \text{Ai}'(s_N) + O(N^{-2/3}e^{-s_N/2}).$$

Before completing the proof, we still require some further preliminaries.

*Properties of the LG transformation.* From (96) it is clear that $\zeta(x_+) = 0$. We exploit a decomposition

$$\dot{\zeta}^2(x) = \hat{f}(x) = k(x)/\xi(x), \tag{110}$$

where, recalling (81), we have

$$k(x) = f(x)/(x - x_N) \quad \text{and} \quad \xi(x) = \zeta(x)/(x - x_N). \tag{111}$$

It follows from Olver (1974), Chapter 11, Lemma 3.1, or directly that $\xi(x)$ is positive and $C^2$ in $(x_0, 1)$, which contains $x_N$. As $x \to x_N$, we have both $\xi(x) \to \dot{\zeta}_N$ and $k(x) \to k_N$, so that $\dot{\zeta}_N^2 = k_N/\dot{\zeta}_N$. Bringing in both (82) and (104), we summarize with

$$\frac{\sin\varphi \sin\gamma}{2\sin^4(\varphi + \gamma)} = k_N = \dot{\zeta}_N^3 = \frac{1}{\kappa_N^2 \sigma_N^3}. \tag{112}$$

Under our asymptotic assumptions (A), we have $\zeta_N(x) \to \zeta_\infty(x)$, along with its first two derivatives, uniformly on compact intervals of $(x_0, 1)$. The dependence on parameters $(\lambda, \mu)$ comes through $x_{N\pm}$ which converge to $x_{\infty\pm}$. From these considerations, we may infer a uniform bound on $D_\delta$ for $\ddot{\zeta}_N$:

$$\sup\{|\ddot{\zeta}_N(x_N + s_N \sigma_N)|, s_N \in [s_L, s_1 N^{1/6}]\} \leq C. \tag{113}$$

Adapting Taylor's expansion (103), and using (104), we find that for some $s^*$ between 0 and $s_N$, and with $x^* = x_N + s^* \sigma_N$,

$$\kappa^{2/3} \zeta(x) = s_N + \tfrac{1}{2} \sigma_N s_N^2 \ddot{\zeta}(x^*)/\dot{\zeta}_N.$$

From (112), it is evident that under assumptions (A), $\dot{\zeta}_N \to \dot{\zeta}_\infty(a_\infty, b_\infty) \in (0, \infty)$. Hence, uniformly for $s_N \in [s_L, s_1 N^{1/6}]$, we have

$$|\kappa^{2/3} \zeta - s_N| \leq C \sigma_N s_N^2. \tag{114}$$

LEMMA 2. *Let $r > 0$ be fixed. For $s_N \geq r^2$, we have $\kappa_N \sigma_N \sqrt{f_N(x)} \geq r$.*

PROOF. Exploiting (111) and then the last two inequalities of (112), we have

$$\sqrt{f(x)} = \sqrt{(x - x_{N+})k(x)} \geq r\sqrt{\sigma_N k_N} = r/(\kappa_N \sigma_N). \qquad \square$$

LEMMA 3. *There exists $C = C(s_L)$ such that for $s_N \geq s_L$,*

$$|\operatorname{Ai}(\kappa_N^{2/3} \zeta(x))| \leq E^{-1}(\kappa_N^{2/3} \zeta(x)) \leq C e^{-s_N}.$$

PROOF. Since $|\operatorname{Ai}(x)| \leq M(x)E^{-1}(x) \leq E^{-1}(x)$, it suffices to use bounds for $E^{-1}$. For $s_N \geq 1$, applying Lemma 2 with $r = 1$, we have

$$\frac{2}{3}\kappa_N \zeta^{3/2} = \kappa_N \int_{x_N}^{x} \sqrt{f} \geq \kappa_N \int_{x_N+\sigma_N}^{x_N+s_N\sigma_N} \sqrt{f} \geq \frac{1}{\sigma_N}(s_N - 1)\sigma_N = s_N - 1.$$

For $x \geq 0$, we have from (93) that $E^{-1}(x) \leq C\exp(-\frac{2}{3}x^{3/2})$, and so

$$E^{-1}(\kappa^{2/3}\zeta) \leq C\exp(-\tfrac{2}{3}\kappa_N \zeta^{3/2}) \leq Ce^{-s_N}.$$

For $s_N \in [s_L, 1]$, it follows from (114) that $|\kappa_N^{2/3}\zeta(x)| \leq C(s_L)$, and hence

$$\sup_{[s_L,1]} |e^{s_N} E^{-1}(\kappa^{2/3}\zeta)| \leq C. \qquad \square$$

PROOF OF PROPOSITION 1.

PROOF OF (108). For error bounds, we use the decomposition suggested by (106):

$$\begin{aligned}\check{\phi}_N(x) &- \operatorname{Ai}(s_N) \\ &= [\bar{e}_N r_N(x) - 1]\operatorname{Ai}(\kappa^{2/3}\zeta) + \operatorname{Ai}(\kappa^{2/3}\zeta) - \operatorname{Ai}(s_N) + \bar{e}_N r_N(x)\varepsilon_2(x, \kappa_N) \\ &= E_{N1} + E_{N2} + E_{N3}.\end{aligned}$$

For the $E_{N1}$ term, first use (113) to conclude that for $s_N \in [s_L, s_1 N^{1/6}]$, we have

$$(115) \qquad |\dot{\zeta}(x)/\dot{\zeta}_N - 1| = \left|\int_{x_N}^{x} \ddot{\zeta}(u)/\dot{\zeta}_N \, du\right| \leq Cs_N \sigma_N.$$

Together with (107), this yields

$$(116) \qquad |\bar{e}_N r_N(x) - 1| \leq C(1 + s_N)\sigma_N.$$

This argument also shows that for $s_N \in [s_L, s_1 N^{1/6}]$,

$$(117) \qquad \tfrac{1}{2} \leq \dot{\zeta}(x)/\dot{\zeta}_N \leq 2.$$

Combined with Lemma 3, we obtain

$$|E_{N1}| \leq C\sigma_N(1 + s_N)e^{-s_N} \leq CN^{-2/3}e^{-s_N/2}.$$

For the $E_{N2}$ term, we first observe, from (114), that for $N > N_0$ we have $|\kappa^{2/3}\zeta - s_N| \leq s_N/4$, and hence that uniformly for $s_N \in (s_L, s_1 N^{1/6})$,

$$|\operatorname{Ai}(\kappa^{2/3}\zeta) - \operatorname{Ai}(s_N)| \leq C\sigma_N s_N^2 \sup\{|\operatorname{Ai}'(t)| : \tfrac{3}{4}s_N \leq t \leq \tfrac{5}{4}s_N\}$$
$$\leq CN^{-2/3}e^{-s_N/2},$$

where we used (118) below.

Finally, for the $E_{N3}$ term, we use (94) and the uniform bound on $\mathcal{V}$ (Section A.3) to get

$$|E_{N3}| \leq C\kappa_N^{-1} r_N(x) M(\kappa^{2/3}\zeta) E^{-1}(\kappa^{2/3}\zeta).$$

For $s_N \in [s_L, 1]$, we observe that $M \leq 1$ and $E \geq 1$, and use (107) together with (117) to conclude that

$$|E_{N3}| \leq C\kappa_N^{-1}\left(\frac{\dot{\zeta}(x)}{\dot{\zeta}_N}\right)^{-1/2} \leq CN^{-2/3} \leq CN^{-2/3}e^{-s_N/2}.$$

For $s_N \in [1, s_1 N^{1/6}]$, we note from (93) that

$$M(\kappa^{2/3}\zeta) \leq c_0 \kappa_N^{-1/6} \zeta^{-1/4}.$$

Since $[\hat{f}(x)]^{-1/4} = [f(x)]^{-1/4}\zeta^{1/4}$, we obtain from the first equality of (107) and then Lemma 2 that

$$r_N(x)M(\kappa^{2/3}\zeta) \leq \frac{c_0}{\sqrt{\kappa_N \sigma_N}} \frac{1}{[f(x)]^{1/4}} \leq \frac{c_0}{\sqrt{r}}.$$

Consequently, from Proposition 3, we arrive at

$$|E_{N3}| \leq C\kappa_N^{-1} E^{-1}(\kappa^{2/3}\zeta) \leq CN^{-2/3}e^{-s_N}. \qquad \square$$

PROPERTIES OF THE AIRY FUNCTION. $\mathrm{Ai}(s)$ satisfies the differential equation $\mathrm{Ai}''(s) = s\,\mathrm{Ai}(s)$. For all $s > 0$, both $\mathrm{Ai}(s) > 0$ and $\mathrm{Ai}'(s) < 0$, and so, from the differential equation, we have $|\mathrm{Ai}'(s)|$ is decreasing for $s > 0$. There are exponential decay bounds: given $s_L$, there exist constants $C_i(s_L)$ such that

(118) $$|\mathrm{Ai}^{(i)}(s)| \leq C_i(s_L)e^{-s}, \qquad s \geq s_L,\ i = 0, 1, 2.$$

6.2. *Approximations at degree $N$ and $N-1$.* The asymptotic model used in this paper supposes that there is a sequence of eigenvalue distribution models, such as (17) and (26), indexed by $N$, the number of variables. The Jacobi parameters $\alpha = \alpha(N)$, $\beta = \beta(N)$ are regarded as functions of $N$. The kernel $S_{N,2}(x,y)$ depends on weighted polynomials $\phi_j = \phi_j(x; \alpha(N), \beta(N))$, $j = 1, \ldots, N$. The Christoffel–Darboux formula and integral representation formulas (36) express $S_{N,2}$ in terms of the two functions

$$\phi_{N-1}(x; \alpha(N), \beta(N)) \quad \text{and} \quad \phi_N(x; \alpha(N), \beta(N)).$$

To construct approximations in the $N$th distribution model, we therefore need separate Liouville–Green asymptotic approximations to both $\phi_{N-1}$ and $\phi_N$. Each of these is defined in turn based on parameters $\lambda, \mu$ and functions $f(x;\lambda,\mu), x_+(\lambda,\mu)$ and $\zeta(x;\lambda,\mu)$. In the case of $\phi_N(x)$, this uses the parameters $(N,\alpha(N),\beta(N))$, while for $\phi_{N-1}(x)$ we use $(N-1,\alpha(N),\beta(N))$. Thus, for example, in comparing the two cases, we have

$$\kappa_N = \kappa(N;\alpha(N),\beta(N)) \qquad \kappa_{N-1} = \kappa(N-1;\alpha(N),\beta(N))$$
$$= 2N + \alpha + \beta + 1, \qquad = 2N + \alpha + \beta - 1,$$
$$x_N = x_+(N;\alpha(N),\beta(N)), \qquad x_{N-1} = x_+(N-1;\alpha(N),\beta(N)),$$
$$\dot\zeta_N = \dot\zeta(x_N;N,\alpha(N),\beta(N)), \qquad \dot\zeta_{N-1} = \dot\zeta(x_{N-1};N-1,\alpha(N),\beta(N)),$$
$$\sigma_N = (\kappa_N^{2/3}\dot\zeta_N)^{-1}, \qquad \sigma_{N-1} = (\kappa_{N-1}^{2/3}\dot\zeta_{N-1})^{-1}$$

and so forth. The analog of (42) for $\phi_{N-1}(x) = \phi_{N-1}(x;\alpha(N),\beta(N))$ states that in terms of the variable $s_{N-1}$:

$$(119) \quad \left(\frac{1-x^2}{\sigma_{N-1}\kappa_{N-1}}\right)^{1/2} \phi_{N-1}(x_{N-1} + \sigma_{N-1}s_{N-1}) = \mathrm{Ai}(s_{N-1}) + \varepsilon_{N-1},$$

with the error bound valid uniformly for compact intervals of $s_{N-1}$.

*Airy approximation to $\phi_{N-1}$.* Define $\check\phi_{N-1}$ correspondingly with every occurrence of $N$ replaced with $N-1$. In parallel with the previous approximations, now with $x = x_{N-1} + s_{N-1}\sigma_{N-1}$ and $s_L \le s_{N-1} \le CN^{1/6}$, we have

$$(120) \quad \check\phi_{N-1}(x) = \mathrm{Ai}(s_{N-1}) + O(N^{-2/3}e^{-s_{N-1}/2}),$$

$$(121) \quad \sigma_{N-1}\check\phi'_{N-1}(x) = \mathrm{Ai}'(s_{N-1}) + O(N^{-2/3}e^{-s_{N-1}/2}).$$

We collect some formulas describing the dependence on $N$ of $a_N, x_N$ and $\sigma_N$. For these formulas we regard $\alpha$ and $\beta$ as constants *not* depending on $N$; in the context of the remark in the previous subsection, we are examining differences between $\phi_N$ and $\phi_{N-1}$ in the $N$th eigenvalue distribution model.

LEMMA 4.

$$(122) \quad \frac{\partial x_{N\pm}}{\partial N} = \frac{\pm 2 \sin^2(\varphi \pm \gamma)}{\kappa_N \sin\varphi \sin\gamma} = \frac{1 - x_{N\pm}^2}{\kappa_N a_N}(1 + O(N^{-1})) = O(N^{-1})$$

and, in particular,

$$(123) \quad \frac{\partial u_N}{\partial N} = \frac{1}{1-x_N^2}\frac{\partial x_N}{\partial N} = \frac{1}{\kappa_N a_N}(1 + O(N^{-1})),$$

$$(124) \quad u_N - u_{N-1} = \partial u_N/\partial N + O(N^{-2}),$$

$$(125) \quad \sigma_N/\sigma_{N-1}, \omega_N/\omega_{N-1}, \tau_N/\tau_{N-1} = 1 + O(N^{-1}),$$

$$(126) \quad \begin{aligned} \Delta_N &:= \frac{u_N - u_{N-1}}{\tau_{N-1}} \\ &= \frac{1}{\tau_N \kappa_N a_N}(1 + \varepsilon_N) = O(N^{-1/3}). \end{aligned}$$

*The constants implicit in the $O(N^{-1/3})$, $O(N^{-1})$ and $O(N^{-2})$ bounds depend on the ratios $\lambda = \alpha/\kappa_N$ and $\mu = \beta/\kappa_N$ defined at* (63).

PROOF. From (63) one obtains

$$\frac{\partial \kappa_N}{\partial N} = 2, \qquad \frac{\partial \lambda}{\partial N} = -\frac{2\alpha}{\kappa^2}, \qquad \frac{\partial \mu}{\partial N} = -\frac{2\beta}{\kappa^2}$$

and then from (68) and (69),

$$(127) \qquad \frac{\partial \gamma}{\partial N} = \frac{2\cos\gamma}{\kappa \sin\gamma}, \qquad \frac{\partial \varphi}{\partial N} = \frac{2\cos\varphi}{\kappa \sin\varphi},$$

and finally from (76),

$$\frac{\partial x_{N\pm}}{\partial N} = \sin(\varphi \pm \gamma)\left(\frac{\partial \varphi}{\partial N} \pm \frac{\partial \gamma}{\partial N}\right) = \frac{2}{\kappa} \frac{\sin(\varphi \pm \gamma)\sin(\gamma \pm \varphi)}{\sin\gamma \sin\varphi}.$$

To obtain the second inequality in (122), use (86) and (76). Since $\partial u/\partial N = (\tanh^{-1})'(x_N)\partial x/\partial N$, (123) follows from (122).

Writing $b_N$ generically for each of $\sigma_N$, $\omega_N$ or $\tau_N$, we obtain (125) by writing

$$\log(b_N/b_{N-1}) = \int_{N-1}^{N} (\log b_t)' \, dt,$$

and verifying that $|(\log b_t)'| = O(N^{-1})$. For example, for $\omega_N = (1 - x_N^2)^{-1} = \sin^{-2}(\varphi + \gamma)$, we have

$$\frac{\partial(\log \omega_N)}{\partial N} = -\frac{2\cos(\gamma + \varphi)}{\sin(\gamma + \varphi)}\left(\frac{\partial \varphi}{\partial N} + \frac{\partial \gamma}{\partial N}\right) = O(N^{-1}),$$

from (127). Similarly, writing $u'_t, u''_t$ for partial derivatives w.r.t. $N$, one verifies that $(\log u'_t)' = u''_t/u'_t = O(N^{-1})$, which shows that $u_N - u_{N-1} = u'_N + O(N^{-2})$, which establishes (124) and allows us to conclude (126) directly from (123) and (125). □

### 7. Unitary case: Theorems 2 and 3.

7.1. *Integral representations for kernel.* In the unitary setting, with eigenvalues $\{x_k\}$ having distribution (17), we have

$$P\left\{\max_k x_k \leq x_0\right\} = \det(I - S_N \chi_0)$$

where $\chi_0(x) = I_{(x_0, 1]}(x)$ and the operator $S_N \chi_0$ is defined via

$$(S_N \chi_0)g(x) = \int_{x_0}^{1} S_N(x, y) g(y) \, dy.$$

Equivalently, we may speak of $S_N$ as an operator on $L_2[x_0, 1)$ with kernel $S_N(x, y)$. On this understanding, we drop further explicit reference to $\chi_0$. Consider now the effect of a change of variables $x = \tau(s)$ with $x_0 = \tau(s_0)$ and $\tau : [s_0, \infty) \to [x_0, 1)$ strictly monotonic. Clearly, with $s_k = \tau^{-1}(x_k)$, we have

$$P\{\max x_k \leq x_0\} = P\{\max s_k \leq s_0\},$$

while we claim (see Section A.6) that

(128) $$\det(I - S_N) = \det(I - S_\tau),$$

where $S_\tau$ is the operator on $L_2(s_0, \infty)$ with kernel

(129) $$S_\tau(s, t) = \sqrt{\tau'(s)\tau'(t)} S_N(\tau(s), \tau(t)).$$

The transformations we consider involve both a nonlinear mapping and a rescaling:

$$\tau(s) = \tau_1 \circ \tau_2(s) = \tanh(\mu + \sigma s).$$

The nonlinear mapping $\tau_1(u) = \tanh u$ has already appeared in (35), giving rise to the integral representation (36) for the kernel $\hat{S}_{N,2}(u, v)$. The rescaling $\tau_2(s) = \mu + \sigma s$ is used for the Airy approximation—it yields the rescaled kernel

(130) $$S_\tau(s, t) = \sigma \hat{S}_N(\mu + \sigma s, \mu + \sigma t).$$

The asymptotic analysis of the edge scaling of the Jacobi kernel has both local and global aspects. The first step is to establish a local Airy approximation to the weighted Jacobi polynomials appearing in (36). The approximation is centered around $x_N = x(N, \alpha, \beta)$, the upper turning point of the differential equation (73) satisfied by $\sqrt{1 - x^2}\phi_N(x)$—this turning point lies within $O(N^{-2/3})$ of the largest zero of $\phi_N$.

The Liouville–Green approximation is made by transforming $x$ to a new independent variable $\zeta$, and as explained at (103)–(104), the scaling $\sigma_N = \sigma(N, \alpha, \beta)$ is defined by

$$\sigma_N = (\kappa_N \dot{\zeta}(x_N))^{-1}.$$

With centering and scaling $(x_N, \sigma_N)$, we establish a local Airy approximation, for $x = x_N + s_N \sigma_N$:

(131) $$\check{\phi}_N(x) = (1 - x^2)^{1/2} \frac{\phi_N(x)}{\sqrt{\kappa_N \sigma_N}} = \text{Ai}(s_N) + O(N^{-2/3} e^{-s_N/2}),$$

uniformly in the variable $s_N$ in the range $s_L \leq s_N \leq C N^{1/6}$. A similar local approximation is shown to hold for $\phi_{N-1}(x)$ with centering $x_{N-1} = x(N - 1, \alpha, \beta)$ and scaling $\sigma_{N-1} = \sigma(N - 1, \alpha, \beta)$. The analog of (131) holds for approximating $\phi_{N-1}(x)$ by $\mathrm{Ai}(s_{N-1})$ with $x = x_{N-1} + s_{N-1}\sigma_{N-1}$ and $s_L \leq s_{N-1} \leq CN^{1/6}$ and $\kappa_N$ replaced by $\kappa_{N-1}$.

To use these local approximations in the integral representation (36), we need a version that applies on the $u$-scale to $\hat{\phi}_N(u)$. Hence, we define

$$(132) \qquad \bar{\phi}_N(u) = \check{\phi}_N(\tanh u), \qquad \bar{\phi}_{N-1}(u) = \check{\phi}_{N-1}(\tanh u).$$

On the $u$-scale, we have $u = u_N + \tau_N t$ with centering $u_N$ and scaling $\tau_N$ defined by

$$(133) \qquad x_N = \tanh u_N, \qquad \tau_N = \omega_N \sigma_N, \qquad \omega_N = (1 - x_N^2)^{-1}.$$

These definitions are suggested by the approximation

$$\tanh(u_N + \tau_N t) \doteq \tanh u_N + \tau_N t \tanh' u_N = x_N + \sigma_N t,$$

where we have used

$$\tanh' u_N = 1/(\tanh^{-1})'(x_N) = 1 - x_N^2 = \omega_N^{-1}.$$

With these definitions, and defining $u_{N-1}, \tau_{N-1}$ and $\omega_{N-1}$ in the corresponding way, it is straightforward (Section 7.2) to show that for $t_L \leq t \leq CN^{1/6}$,

$$\bar{\phi}_N(u_N + \tau_N t) = \mathrm{Ai}(t) + O(N^{-2/3} e^{-t/2}),$$

$$\bar{\phi}_{N-1}(u_{N-1} + \tau_{N-1} t) = \mathrm{Ai}(t) + O(N^{-2/3} e^{-t/2}).$$

We thus obtain good local Airy approximations for both $\bar{\phi}_N$ and $\bar{\phi}_{N-1}$, but with differing centering and scaling values. Our goal is a scaling limit with error term for the kernel $\hat{S}_N$, and so the centering and scaling for $\hat{S}_N$ will need to combine those for $\bar{\phi}_N$ and $\bar{\phi}_{N-1}$ in some fashion. In addition, the integral representation for $\hat{S}_N$ involves global features of $\bar{\phi}_N$ and $\bar{\phi}_{N-1}$, through the transformation $x = \tanh u$. Thus, we use a rescaling $u = \mu + \sigma s$ with the explicit values of $(\mu, \sigma)$ given for real and complex cases in Section 7.3 below, and put

$$(134) \qquad \phi_\tau(s) = \bar{\phi}_N(\mu + \sigma s), \qquad \psi_\tau(s) = \bar{\phi}_{N-1}(\mu + \sigma s).$$

We now convert (36) into a representation in terms of $\phi_\tau$ and $\psi_\tau$. First, observe from (134) and (132) that

$$\phi_\tau(s) = \check{\phi}_N(\tanh(\mu + \sigma s)), \qquad \psi_\tau(s) = \check{\phi}_{N-1}(\tanh(\mu + \sigma s)).$$

From the definition of $\check{\phi}_N$ in (131), we also have

$$(135) \qquad \check{\phi}_N(\tanh u) = \frac{\phi_N(\tanh u)}{\sqrt{\kappa_N \sigma_N} \cosh u} = \frac{\hat{\phi}_N(u)}{\sqrt{\kappa_N \sigma_N}},$$

with a corresponding identity for $\check{\phi}_{N-1}$. Combining the last two displays yields

$$\hat{\phi}_N(\mu + \sigma s) = \sqrt{\kappa_N \sigma_N} \phi_\tau(s), \qquad \hat{\phi}_{N-1}(\mu + \sigma s) = \sqrt{\kappa_{N-1} \sigma_{N-1}} \psi_\tau(s).$$

From (130) and (36) and a change of variables $u = \mu + \sigma s$, $v = \mu + \sigma t$ and $w = \sigma z$,

$$S_\tau(s,t) = \frac{\sigma^2}{2}(\kappa_N - 1)a_N \int_0^\infty \hat{\phi}_N(\mu + \sigma(s+z))\hat{\phi}_{N-1}(\mu + \sigma(t+z)) + \cdots dz$$

$$= \frac{\sigma^2}{2}(\kappa_N - 1)a_N\sqrt{\sigma_N \kappa_N \sigma_{N-1} \kappa_{N-1}} \int_0^\infty \phi_\tau(s+z)\psi_\tau(t+z) + \cdots dz.$$

Thus, we arrive at

(136) $$S_\tau(s,t) = \sigma \hat{S}_N(\mu + \sigma s, \mu + \sigma t) = e_N \bar{S}_\tau(s,t),$$

where

(137) $$\bar{S}_\tau(s,t) = \tfrac{1}{2} \int_0^\infty [\phi_\tau(s+z)\psi_\tau(t+z) + \psi_\tau(s+z)\phi_\tau(t+z)] dz$$

and

(138) $$e_N := \sigma^2(\kappa_N - 1)a_N \sqrt{\sigma_N \sigma_{N-1} \kappa_N \kappa_{N-1}} = 1 + O(N^{-1}).$$

PROOF. Using (122) and (126) combined with $\limsup |x_N| < 1$, we find that both $\sigma_N/\sigma_{N-1}$ and $\omega_N/\omega_{N-1}$ are $1 + O(N^{-1})$, and so

(139) $$e_N = \sigma_N^3 \kappa_N^2 \cdot \omega_N^2 a_N (1 + O(N^{-1})).$$

Combining (82) and (133) of $\omega_N$ with (112), we obtain

(140) $$\frac{1}{4}(x_{N+} - x_{N-}) = k_N = \frac{1}{\kappa_N^2 \sigma_N^3}.$$

From Lemma 1 and (78), $\frac{1}{4}(x_+ - x_-) = a_N(1 + O(N^{-1}))$, and so (140) shows that, indeed, $e_N = 1 + O(N^{-1})$. □

To summarize, we have

(141) $$P\{(\max u_k - \mu)/\sigma \leq s_0\} = \det(I - S_\tau).$$

The Tracy–Widom distribution

$$F_2(s_0) = \det(I - S_A),$$

where the Airy kernel

(142) $$S_A(s,t) = \int_0^\infty \mathrm{Ai}(s+z)\mathrm{Ai}(t+z) dz.$$

To bound the convergence rate of (141) to $F_2(s_0)$, we use the Seiler–Simon bound (32). To bound $S_\tau - S_A$, we use a simple algebraic identity:

(143) $$\begin{aligned}4[\phi\bar{\psi} + \psi\bar{\phi} - 2a\bar{a}] \\ = (\phi + \psi + 2a)(\bar{\phi} + \bar{\psi} - 2\bar{a}) + (\phi + \psi - 2a)(\bar{\phi} + \bar{\psi} + 2\bar{a}) \\ - 2(\phi - \psi)(\bar{\phi} - \bar{\psi}).\end{aligned}$$

Inspection of (143) shows that the essential bounds on simultaneous Airy approximation of $\phi_\tau$ and $\psi_\tau$ are those given in the following lemma. Set

$$\text{(144)} \qquad \mu = \frac{\tau_N^{-1} u_N + \tau_{N-1}^{-1} u_{N-1}}{\tau_N^{-1} + \tau_{N-1}^{-1}}, \qquad \sigma^{-1} = \tfrac{1}{2}(\tau_N^{-1} + \tau_{N-1}^{-1}).$$

LEMMA 5 (Complex case). *There exists $C = C(\cdots)$ such that for $s \geq s_L$,*

$$\text{(145)} \qquad |\phi_\tau(s)| \leq Ce^{-s},$$

$$\text{(146)} \qquad |\psi_\tau(s)| \leq Ce^{-s},$$

$$\text{(147)} \qquad |\phi_\tau(s) - \operatorname{Ai}(s)| \leq CN^{-1/3} e^{-s/4},$$

$$\text{(148)} \qquad |\psi_\tau(s) - \operatorname{Ai}(s)| \leq CN^{-1/3} e^{-s/4},$$

$$\text{(149)} \qquad |\phi_\tau(s) + \psi_\tau(s) - 2\operatorname{Ai}(s)| \leq CN^{-2/3} e^{-s/4}.$$

We remark that the same bounds trivially hold if $\phi_\tau(s)$ and $\psi_\tau(s)$ are replaced by $\sqrt{e_N}\phi_\tau(s)$ and $\sqrt{e_N}\psi_\tau(s)$, respectively.

Of these bounds, the most critical is (149), which provides the $N^{-2/3}$ rate of convergence. We wish here also to acknowledge the influence of El Karoui (2006), whose methods allowed the formulation and proof of Lemma 5, which improved on our earlier, less rigorous approach.

Inspecting (137) and (142), and setting $\phi, \psi$ and $a$ equal to $\sqrt{e_N}\phi_\tau(s+z)$, $\sqrt{e_N}\psi_\tau(s+z)$ and $\operatorname{Ai}(s+z)$, respectively, and $\bar\phi, \bar\psi$ and $\bar a$ to corresponding quantities with $t$ in place of $s$, we are led to an expression for $K_N = S_\tau - S_A$ having the form

$$K_N(s,t) = \int_0^\infty \sum_{i=1}^r a_i(s+z) b_i(t+z)\, dz.$$

Lemma 5 leads to bounds of the form

$$\text{(150)} \qquad |a_i(s)| \leq a_{Ni} e^{-as}, \qquad |b_i(s)| \leq b_{Ni} e^{-as}, \qquad s \geq s_0$$

and hence

$$\text{(151)} \qquad \|K_N\|_1 \leq \frac{e^{-2as_0}}{4a^2} \sum_{i=1}^r a_{Ni} b_{Ni}.$$

To make explicit the role of the rate bounds in Lemma 5, we may write, in the case of $K_N = S_\tau - S_A$, with $a = 1/4$,

$$\|S_\tau - S_A\|_1 \leq C(1 \cdot N^{-2/3} + N^{-2/3} \cdot 1 + N^{-1/3} \cdot N^{-1/3}) e^{-s_0/2} = CN^{-2/3} e^{-s_0/2}.$$

PROOF OF (151). Recall first that the trace class norm of a rank-1 operator $\phi \otimes \psi$ is just $\|\phi\|_2 \|\psi\|_2$. Since the trace norm of an integral is at

most the integral of the trace norms,

$$\|K_N\|_1 \leq \int_0^\infty \sum_i \|a_i(\cdot + z)\|_2 \|b_i(\cdot + z)\|_2 \, dz.$$

From (150), we have $\|a_i(\cdot + z)\|_2^2 \leq a_{Ni}^2 \int_{s_0}^\infty e^{-2a(s+z)} \, ds = a_{Ni}^2 e^{-2a(s_0+z)}/(2a)$, and so, after further integration, we obtain (151). □

REMARK (Summing up: complex case). In conjunction with (141) and (32), this leads to the bound in Theorem 3. Theorem 2 is derived from Theorem 3 in exactly the same manner as Theorem 1 is deduced from Theorem 4, as is detailed in Section 8.4.3. In particular, the logit scale quantities $w_N$ and $\omega_N$ below Theorem 2 are given by $w_N = 2u_N$ and $\omega_N = 2\tau_N$.

### 7.2. *Local bounds—u-scale.*

PROPOSITION 2. *For $t_L \leq t \leq CN^{1/6}$,*

$$\begin{aligned}\bar\phi_N(u_N + \tau_N t) &= \mathrm{Ai}(t) + O(N^{-2/3} e^{-t/2}), \\ \tau_N \bar\phi_N'(u_N + \tau_N t) &= \mathrm{Ai}'(t) + O(N^{-2/3} e^{-t/2})\end{aligned} \quad (152)$$

*and similarly*

$$\begin{aligned}\bar\phi_{N-1}(u_{N-1} + \tau_{N-1} t) &= \mathrm{Ai}(t) + O(N^{-2/3} e^{-t/2}), \\ \tau_{N-1} \bar\phi_{N-1}'(u_{N-1} + \tau_{N-1} t) &= \mathrm{Ai}'(t) + O(N^{-2/3} e^{-t/2}).\end{aligned} \quad (153)$$

PROOF. Consider first $\bar\phi_N(u) = \check\phi_N(\tanh(u_N + \tau_N t))$. By Taylor expansion

$$\begin{aligned}\tanh(u_N + \tau_N t) &= \tanh u_N + \tau_N t \tanh' u_N + \tfrac{1}{2}\tau_N^2 t^2 \tanh''(u^*) \\ &= x_N + \sigma_N(t + \varepsilon_N(t)),\end{aligned} \quad (154)$$

where we use $\tanh' u_N = 1/\omega_N$, and note that for $t_L \leq t \leq CN^{1/6}$,

$$(155) \quad |\varepsilon_N(t)| = |(\omega_N/2)\tau_N t^2 \tanh''(u^*)| \leq C\tau_N t^2 \leq Ct^2 N^{-2/3} \leq CN^{-1/3}.$$

Consequently, from (108),

$$\begin{aligned}\bar\phi_N(u) &= \check\phi_N(x_N + \sigma_N(t + \varepsilon_N(t))) \\ &= \mathrm{Ai}(t + \varepsilon_N(t)) + O(N^{-2/3} e^{-(t+\varepsilon_N(t))/2}) \\ &= \mathrm{Ai}(t) + O(N^{-2/3} e^{-t/2}),\end{aligned}$$

after we appeal to (155) and also observe using the remarks before (118) that

(156)
$$\begin{aligned}|\mathrm{Ai}(t+\varepsilon_N(t)-\mathrm{Ai}(t))|\\ \leq \varepsilon_N(t)\sup\{|\mathrm{Ai}'(u)|:t-\varepsilon_N(t)\leq u\leq t+\varepsilon_N(t)\}\\ \leq CN^{-2/3}t^2\,\mathrm{Ai}'(t-CN^{-1/3})\\ \leq O(N^{-2/3}e^{-t/2}).\end{aligned}$$

Turn now to

$$\tau_N\bar{\phi}'_N(u)=\tau_N(\tanh'u)\check{\phi}'_N(\tanh u)=(\omega_N\tanh'u)\sigma_N\check{\phi}'_N(x_N+\sigma_N(t+\varepsilon_N(t))).$$

Noting that $\omega_N\tanh'u=\omega_N\tanh'u_N+\omega_N\tau_Nt\tanh''u^*=1+O(tN^{-2/3})$, we find

$$\tau_N\bar{\phi}'_N(u)=[1+O(N^{-2/3}t)][\mathrm{Ai}'(t+\varepsilon_N(t))+O(N^{-2/3}e^{-(t+\varepsilon_N(t))/2})].$$

The argument at (156) applies equally with $\mathrm{Ai}'$ in place of $\mathrm{Ai}$, and so

$$\begin{aligned}\tau_N\bar{\phi}'_N(u)&=[1+O(N^{-2/3}t)][\mathrm{Ai}'(t)+O(N^{-2/3}e^{-t/2})]\\ &=\mathrm{Ai}'(t)+O(N^{-2/3}e^{-t/2}).\end{aligned}$$

The arguments for $\bar{\phi}_{N-1}$ and $\tau_{N-1}\bar{\phi}'_{N-1}$ are entirely similar. □

7.3. *Global bounds.* The global bounds that we need for Airy approximation to $\phi_\tau$ and $\psi_\tau$ are very similar in the real and complex cases. The differences in the two statements arise first from the changes in choice of centering $\mu$ and scaling $\sigma$, and second because bounds on convergence of derivatives are also required for the real case.

In the complex case, use (144), and in the real case put

(157) $$\mu=u_N,\qquad \sigma=\tau_N.$$

In either case, set

$$\begin{aligned}\phi_\tau(t)&=\bar{\phi}_N(\mu+\sigma t)=\check{\phi}_N(\tanh(\mu+\sigma t)),\\ \psi_\tau(t)&=\bar{\phi}_{N-1}(\mu+\sigma t)=\check{\phi}_{N-1}(\tanh(\mu+\sigma t)).\end{aligned}$$

The results for the complex case were given in Lemma 5.

LEMMA 6 (Real case). *There exists $C=C(\cdots)$ such that for $s\geq s_L$,*

(158) $$|\phi_\tau(s)|,\qquad |\phi'_\tau(s)|\leq Ce^{-s},$$

(159) $$|\psi_\tau(s)|,\qquad |\psi'_\tau(s)|\leq Ce^{-s},$$

$$|\phi_\tau(s)-\mathrm{Ai}(s)|\leq CN^{-2/3}e^{-s/4},$$

(160)
$$|\phi'_\tau(s) - \text{Ai}'(s)| \leq CN^{-2/3}e^{-s/4},$$

(161)
$$|\psi_\tau(s) - \text{Ai}(s) - \Delta_N \text{Ai}'(s)| \leq CN^{-2/3}e^{-s/4},$$
$$|\psi'_\tau(s) - \text{Ai}'(s) - \Delta_N \text{Ai}''(s)| \leq CN^{-2/3}e^{-s/4}.$$

A trivial corollary of Lemma 6 that will be also needed in the real case is that right-tail integrals $(\tilde{\varepsilon}\psi)(s) = \int_s^\infty \psi(s)\,ds$ satisfy the bounds of type (158)–(161) whenever $\psi$ does.

We shall give proofs of both complex and real cases, indicating the parts that are in common and that are divergent. Proofs for the bounds involving $\phi'_\tau$ and $\psi'_\tau$ are deferred to Section A.8.

*Bounds for $\phi_\tau, \psi_\tau$.* Combine (106) with the bound $|\text{Ai}(x)| \leq M(x)E^{-1}(x)$ and with (94). Since $\mathcal{V}(\zeta)$ is bounded (Section A.3), we have, for $x = \tanh(\mu + \sigma s)$,

$$|\phi_\tau(s)| \leq C\bar{e}_N r_N(x) M(\kappa_N^{2/3}\zeta_N) E^{-1}(\kappa_N^{2/3}\zeta_N),$$

where $\bar{e}_N$ and $r_N(x)$ are defined near (107). Here and in the argument below, analogous bounds hold for $\psi_\tau$ with $N$ replaced by $N-1$.

We consider two cases: $s \in [s_L, s_1]$ and $s \in (s_1, \infty)$.

  (i) Large $s$. First we argue as for the $E_{N3}$ term in the local bound, that

(162)
$$r_N(x) M(\kappa_N^{2/3}\zeta_N) \leq c_0/\sqrt{r}.$$

(This argument is valid for all $s \geq s_1$ and for the $N-1$ case.) For the exponential term, it is convenient to modify slightly the argument used for Lemma 2 and Proposition 3. We suppose that $s_1$ is chosen so that $x \geq \tau(s_1) \geq \max(x_N + r^2\sigma_N, x_{N-1} + r^2\sigma_{N-1})$. Then we have (using the $N$-case as lead example)

$$\sqrt{f(x)} \geq \frac{r\sqrt{\sigma_N(x_{N+} - x_{N-})}}{2(1-x^2)}$$

so as to exploit the inverse hyperbolic tangent integral

$$\int_{\tau(s_1)}^{\tau(s)} \frac{dx}{1-x^2} = \tanh^{-1}\tau(s) - \tanh^{-1}\tau(s_1) = \sigma(s - s_1)$$

to conclude using (140) that

$$\frac{2}{3}\kappa_N\zeta_N^{3/2} = \kappa_N \int_{x_N}^x \sqrt{f_N} \geq \frac{r}{2}\kappa_N\sqrt{\sigma_N(x_{N+} - x_{N-})}\sigma(s-s_1) = r\frac{\sigma}{\sigma_N\omega_N}(s-s_1).$$

Of course, for $\psi_\tau$, the denominator of the last expression is $\sigma_{N-1}\omega_{N-1}$. If $N > N_0$ is chosen large enough so that

$$\frac{3}{4} \leq \frac{\sigma_N \omega_N}{\sigma_{N-1} \omega_{N-1}} \leq \frac{4}{3},$$

then in both real and complex cases

$$\frac{\sigma}{\sigma_N \omega_N}, \qquad \frac{\sigma}{\sigma_{N-1}\omega_{N-1}} \geq \frac{3}{4},$$

so that if we take, say, $r = 4/3$, then in all cases

$$r\frac{\sigma}{\sigma_N \omega_N}(s - s_1) \geq s - s_1.$$

Since $E^{-1}(x) \leq C \exp(-\frac{2}{3}x^{3/2})$, we have

(163) $\quad E^{-1}(\kappa_N^{2/3}\zeta_N) \leq C\exp(-\frac{2}{3}\kappa_N \zeta_N^{3/2}) \leq C\exp(-(s-s_1)) \leq Ce^{-s}.$

Combining (162) and (163), we get, for $s \geq s_1$, that

$$|\phi_\tau(s)| \quad \text{and} \quad |\psi_\tau(s)| \leq Ce^{-s}.$$

(ii) Small $s \in [s_L, s_1]$. We use the bounds $M \leq 1$, $E \geq 1$ and use (107) together with (113) to conclude that

$$|\phi_\tau(s)| \leq Cr_N(x) \leq C \leq Ce^{-s}.$$

For the Airy approximation bounds, we must paste together the local bounds of Proposition 2, derived separately for indices $N$ and $N-1$, into the single scaling $\mu + \sigma t$; and then second, develop adequate bounds for $t \geq t_1 N^{1/6}$.

PROOF OF (149). We establish this first, as it indicates the reason for the choices (144). In order to use (152) and (153), we set

(164) $\qquad \mu + \sigma t = u_N + \tau_N t_N = u_{N-1} + \tau_{N-1} t_{N-1},$

which yields

(165)
$$\phi_\tau(t) + \psi_\tau(t) = \bar\phi_N(u_N + \tau_N t_N) + \bar\phi_{N-1}(u_{N-1} + \tau_{N-1}t_{N-1})$$
$$= \mathrm{Ai}(t_N) + \mathrm{Ai}(t_{N-1}) + O(N^{-2/3}(e^{-t_N/2} + e^{-t_{N-1}/2})).$$

Rewriting (164), we have $t_j = t + c_j + d_j t$, where

$$c_j = \tau_j^{-1}(\mu - u_j), \qquad d_j = \sigma\tau_j^{-1} - 1, \qquad j = N-1, N.$$

Consequently, we have both

$$\mathrm{Ai}(t_N) = \mathrm{Ai}(t) + (c_N + d_N t)\,\mathrm{Ai}'(t) + \tfrac{1}{2}(c_N + d_N t)^2 \,\mathrm{Ai}''(t_N^*),$$
$$\mathrm{Ai}(t_{N-1}) = \mathrm{Ai}(t) + (c_{N-1} + d_{N-1}t)\,\mathrm{Ai}'(t) + \tfrac{1}{2}(c_{N-1} + d_{N-1}t)^2 \,\mathrm{Ai}''(t_{N-1}^*).$$

In the approximation for $\phi_\tau + \psi_\tau$, the terms in $\mathrm{Ai}'(t)$ drop out if we choose $\mu$ and $\sigma$ so that $c_N + c_{N-1} = d_N + d_{N-1} = 0$, and this leads immediately to expressions (144).

With these choices of $\mu$ and $\sigma$, we find from (125) and (126) that

$$c_N = \frac{u_{N-1} - u_N}{\tau_{N-1} + \tau_N} = O(N^{-1/3}), \qquad d_N = \frac{\tau_{N-1} - \tau_N}{\tau_{N-1} + \tau_N} = O(N^{-1}),$$

so that

(166) $$\sup\{|c_N + d_N t| : t_L \leq t \leq t_1 N^{1/6}\} \leq C N^{-1/3}$$

and so

$$\mathrm{Ai}(t_N) + \mathrm{Ai}(t_{N-1}) = 2\,\mathrm{Ai}(t) + O(N^{-2/3}(\mathrm{Ai}''(t_N^*) + \mathrm{Ai}''(t_{N-1}^*))).$$

We conclude for $j = N-1, N$, that

(167) $$t_j \geq t - CN^{-1/3} \quad \text{and} \quad |\mathrm{Ai}''(t_j^*)| \leq C e^{-t/2}.$$

Combining these bounds with (165) establishes (149) for $t \in [t_L, t_1 N^{1/6}]$.

For $t \geq t_1 N^{1/6}$, crude arguments suffice. Indeed, from (145) [and with a similar argument for $\psi_\tau(t)$],

(168) $$|\phi_\tau(t) - \mathrm{Ai}(t)| \leq |\phi_\tau(t)| + |\mathrm{Ai}(t)| \leq C e^{-t} + C e^{-t}$$

(169) $$\leq C N^{-2/3} e^{-t/4}.$$

*Bound for* $\phi_\tau(s) - \mathrm{Ai}(s)$. For $t \geq t_1 N^{1/6}$, we may reuse the bounds for $\phi_\tau$ (and $\psi_\tau$) at (168). Consider, then, the interval $t \in [t_L, t_1 N^{1/6}]$. In the real case, since $\mu = u_N, \sigma = \tau_N$, the bound needed is already established at (152). In the complex case, combining elements from the argument above, we have

$$\phi_\tau(t) = \mathrm{Ai}(t) + (c_N + d_N t)\,\mathrm{Ai}'(t^*) + O(N^{-2/3} e^{-t_N/2}),$$

and (147) follows from (166) and (167). The proof of (148) is analogous.

*Real case bound for* $\psi_\tau(s) - \mathrm{Ai}(s) - \Delta_N \mathrm{Ai}'(s)$. As with earlier cases, the real work lies for $s \in [s_L, s_1 N^{1/6}]$. From the definitions, $\psi_\tau(t) = \bar{\phi}_{N-1}(\mu + \sigma t) = \bar{\phi}_{N-1}(u_N + \tau_N t)$. In order to use (153), we write $u_N + \tau_N t = u_{N-1} + \tau_{N-1} t'$, so that

$$\psi_\tau(t) = \mathrm{Ai}(t') + O(N^{-2/3} e^{-t'/2}),$$

where, using (126),

(170) $$t' - t = \Delta_N + (\tau_N \tau_{N-1}^{-1} - 1)t = \Delta_N + O(tN^{-1}).$$

Consequently,

(171) $$\mathrm{Ai}(t') = \mathrm{Ai}(t) + [\Delta_N + O(tN^{-1})]\,\mathrm{Ai}'(t) + \tfrac{1}{2}[\Delta_N + O(tN^{-1})]^2\,\mathrm{Ai}''(t^*)$$

and since (126) shows that $\Delta_N = O(N^{-1/3})$, we conclude that

(172) $\qquad \text{Ai}(t') = \text{Ai}(t) + \Delta_N \text{Ai}'(t) + O(N^{-2/3} e^{-t/2}).$

Since (170) shows that $e^{-t'/2} = e^{-t/2 + O(N^{-1/3})} \leq Ce^{-t/2}$ (for $t \in t_1 N^{1/6}$), we obtain the desired bound in (161) for $\psi_\tau(t) - \text{Ai}(t) - \Delta_N \text{Ai}'(t)$.

## 8. Orthogonal case: Theorems 1 and 4.

8.1. *Derivation of (48)–(50).* Tracy and Widom (1998) provide a direct derivation of Fredholm determinant representations for eigenvalue probabilities that avoids the Introduction of quaternion determinants. We first review this, with the aim of then making the connection to the results of Adler et al. (2000) (abbreviated as AFNM below), so as to obtain the explicit representation (50) of $S_{N+1,1}$ as a rank-1 modification of a multiple of the unitary kernel $S_{N,2}$.

*Tracy–Widom derivation.* Accordingly, we now fix notation, and specify precisely our use of Tracy and Widom (1998). In the case of the Jacobi orthogonal ensemble (26), the weight function $w(x) = (1-x)^{(\alpha-1)/2}(1+x)^{(\beta-1)/2}$.

Setting $f(x) = -I\{x \geq x_0\} = -\chi_0(x)$, we write the exceedance probability in the form used in Tracy and Widom [(1998), Section 9]:

$$P\left\{\max_{1 \leq k \leq N+1} x_k \leq x_0\right\}$$
$$= E \prod_{k=1}^{N+1} (1 - I\{x_k \geq x_0\})$$
$$= \int \cdots \int \prod_{j<k} |x_j - x_k| \prod_j w(x_j) \prod_j (1 + f(x_j)) \, dx_1 \cdots dx_{N+1}.$$

The argument of Tracy and Widom [(1998), Section 9] establishes that $K_{N+1}$ satisfies (49) with $S_{N+1}$ expressed in the form

(173) $\qquad S_{N+1,1}(x, y) = -\sum_{j,k=0}^{N} \psi_j(x) \mu_{jk} (\varepsilon \psi_k)(y).$

Here $\psi_j(x) = p_j(x) w(x)$ and $\{p_j(x), j = 0, \ldots, N\}$ is an arbitrary sequence of polynomials of exact degree $j$. The coefficients $\{\mu_{jk}\}$ are the entries of $M^{-1}$, where

$$M_{jk} = \iint \varepsilon(x-y) \psi_j(x) \psi_k(y) \, dx \, dy.$$

The function $\varepsilon(x) = \frac{1}{2}\text{sgn}(x)$, and, as usual,

$$(\varepsilon\psi_k)(y) = \int \varepsilon(y-z)\psi_k(z)\,dz.$$

The next step is to make a specific choice of polynomials $p_j$ and hence $\psi_j$.

*Connecting to AFNM.* The key to AFNM's summation formula linking orthogonal and unitary ensembles is the observation that if the unitary ensemble has weight function

$$w_2(x) = e^{-2V(x)} = (1-x)^\alpha(1+x)^\beta,$$

then the corresponding orthogonal ensemble should have a modified weight function

$$\tilde{w}_1(x) = e^{-\tilde{V}(x)} = (1-x)^{(\alpha-1)/2}(1+x)^{(\beta-1)/2}.$$

The derivation of AFNM exploits a sequence of polynomials $\{\tilde{q}_k(x)\}$ that are skew-orthogonal with respect to the skew inner product:

$$\langle f, g\rangle_1 = \iint \varepsilon(y-x) f(x) g(y) \tilde{w}_1(x)\tilde{w}_1(y)\,dx\,dy.$$

Skew-orthogonality means that

$$\langle \tilde{q}_{2j}, \tilde{q}_{2k+1}\rangle_1 = -\langle \tilde{q}_{2k+1}, \tilde{q}_{2j}\rangle_1 = \tilde{r}_j \delta_{jk},$$
$$\langle \tilde{q}_{2j}, \tilde{q}_{2k}\rangle_1 = \langle \tilde{q}_{2j+1}, \tilde{q}_{2k+1}\rangle_1 = 0.$$

Given such a skew-orthogonal sequence, we may fix the functions $\psi_j$ appearing in (173) by setting $p_k = \tilde{q}_k/\sqrt{\tilde{r}_{[k/2]}}$. Since $M_{jk} = -\langle \tilde{q}_j, \tilde{q}_k\rangle_1/\tilde{r}_{[j/2]}$, it follows that $M^{-1}$ is a direct sum of $L = (N+1)/2$ copies of $\begin{pmatrix} 0 & 1 \\ -1 & 0 \end{pmatrix}$, so that (49) takes the form

$$(174) \quad S_{N+1,1}(x,y) = \sum_{k=0}^{L-1}[-\psi_{2k}(x)\varepsilon\psi_{2k+1}(y) + \psi_{2k+1}(x)\varepsilon\psi_{2k}(y)].$$

With the following notational dictionary:

| AFNM Sec. 2 | TW |
|---|---|
| $e^{-\tilde{V}(x)}$ | $w(x)$ |
| $\tilde{q}_k(x)/\sqrt{\tilde{r}_{[k/2]}}$ | $p_k(x)$ |
| $\tilde{q}_k e^{-\tilde{V}}/\sqrt{\tilde{r}_{[k/2]}}$ | $\psi_k(x)$ |
| $\tilde{\Phi}_k(x)/\sqrt{\tilde{r}_{[k/2]}}$ | $\varepsilon\psi_k(x)$ |

we may relate $S_{N+1,1}$ to the function $\tilde{S}_1$ given in Adler et al. [(2000), (2.9)] by

$$S_{N+1,1}(x,y) = \tilde{S}_1(y,x). \tag{175}$$

Actually $\tilde{S}_1$ is defined in AFNM, Section 4 by modifying their formula (2.9) to replace $V$ and $q_k$ by $\tilde{V}$ and $\tilde{q}_k$. This modification has already been incorporated in our discussion.

*Rewriting the AFNM summation formula.* In this subsection only, for consistency with the notation of AFNM, let $p_j(x)$ be the monic orthogonal polynomial associated with $e^{-2V(x)} = w_2(x)$, and define $\gamma_{N-1}$ via

$$\gamma_{N-1}\|p_{N-1}\|_2^2 = \tilde{\kappa}_N = \tfrac{1}{2}(2N + \alpha + \beta).$$

Then, in the Jacobi case, noting that $e^{-(V(x)-\tilde{V}(x))} = \sqrt{1-x^2}$, AFNM's Proposition 4.2 gives the formula, for $N+1$ even,

$$\tilde{S}_1(x,y) = \sqrt{\frac{1-x^2}{1-y^2}} S_{N,2}(x,y) \\ + \gamma_{N-1} e^{-\tilde{V}(y)} p_N(y) \int \varepsilon(x-t) e^{-\tilde{V}(t)} p_{N-1}(t)\,dt. \tag{176}$$

Using the Jacobi polynomial notation of Section 3, we have $p_N = P_N/l_N$. From the definitions, we have $\|p_N\| = \sqrt{h_N}/l_N$ and $a_N = \|p_N\|/\|p_{N-1}\|$. So we may write, using (15) and then (52),

$$e^{-\tilde{V}(y)} p_N(y) = \frac{w_2^{1/2}(y)}{\sqrt{1-y^2}} \frac{P_N(y)}{l_N} = \frac{\phi_N(y)}{\sqrt{1-y^2}} \frac{\sqrt{h_N}}{l_N} = \|p_N\| \tilde{\phi}_N(y)$$

and

$$\int \varepsilon(x-t) e^{-\tilde{V}(t)} p_{N-1}(t)\,dt = \|p_{N-1}\|(\varepsilon\tilde{\phi}_{N-1})(x).$$

The second term in (176) becomes

$$a_N \tilde{\kappa}_N \tilde{\phi}_N(y)(\varepsilon\tilde{\phi}_{N-1})(x).$$

Interchanging the roles of $x$ and $y$ as directed by (175), we obtain (50). [We remark that the possibility of expressing the orthogonal kernel in terms of the unitary kernel plus a finite rank term was shown already by Widom (1999).]

8.2. *Transformation and scaling.* As in the unitary case, to describe the scaling limit it is convenient to use a nonlinear mapping and rescaling $\tau(s) = \tanh(\mu + \sigma s)$ with $\mu = \mu(N)$ and $\sigma = \sigma(N)$ being modified from the unitary setting and to be specified as in (157).

For the matrix kernel appearing in (48) we have, in parallel with (128) and its proof,

$$\det(I - K_{N+1}\chi_0) = \det(I - K_\tau),$$

where $K_\tau$ is an operator with matrix kernel

$$(177) \qquad K_\tau(s,t) = \sqrt{\tau'(s)\tau'(t)} K_{N+1}(\tau(s), \tau(t)).$$

Introducing again $s_k = \tau^{-1}(x_k)$, our aim is to study the convergence of

$$(178) \quad \begin{aligned} F_{N+1}(s_0) &= P\Big\{ \max_{1 \leq k \leq N+1} s_k \leq s_0 \Big\} = P\Big\{ \max_{1 \leq k \leq N+1} x_k \leq x_0 \Big\} \\ &= \sqrt{\det(I - K_\tau)} \end{aligned}$$

to

$$F_1(s_0) = \sqrt{\det(I - K_{GOE})},$$

where, following [Tracy and Widom (2005)]

$$(179) \qquad K_{GOE}(s,t) = \begin{bmatrix} S(s,t) & SD(s,t) \\ IS(s,t) - \varepsilon(s-t) & S(t,s) \end{bmatrix}$$

and the entries of $K_{GOE}$ are given by

$$(180) \quad \begin{aligned} S(s,t) &= S_A(s,t) + \tfrac{1}{2}\operatorname{Ai}(s)\Big(1 - \int_t^\infty \operatorname{Ai}(u)\,du\Big), \\ SD(s,t) &= -\partial_t S_A(s,t) - \tfrac{1}{2}\operatorname{Ai}(s)\operatorname{Ai}(t), \\ IS(s,t) &= -\int_s^\infty S_A(u,t)\,du \\ &\quad + \tfrac{1}{2}\Big(\int_t^s \operatorname{Ai}(u)\,du + \int_s^\infty \operatorname{Ai}(u)\,du \int_t^\infty \operatorname{Ai}(u)\,du\Big), \end{aligned}$$

where $S_A$ is the Airy kernel defined at (142).

Tracy and Widom (2005) describe with some care the nature of the operator convergence of $K_{N+1}\chi_0$ to $K_{GOE}$ for the Gaussian finite $N$ ensemble. We adapt and extend their approach to the Jacobi finite $N$ ensemble focusing on the associated $N^{-2/3}$ rate of convergence. We therefore repeat, for reader convenience, their remarks on weighted Hilbert spaces and regularized 2-determinants in the current setting.

Let $\rho$ be any weight function for which

(181) $$\int_\infty^\infty \frac{1+s^2}{\rho(s)} ds < \infty \quad \text{and}$$

(182) $$\rho(s) \leq C(1+|s|^r) \quad \text{for some positive integer } r.$$

As in Tracy and Widom (2005), we consider $2 \times 2$ Hilbert–Schmidt operator matrices $T$ with trace class diagonal entries. Write $L^2(\rho)$ and $L^2(\rho^{-1})$ for the spaces $L^2((s_0, \infty), \rho(s)\,ds)$ and $L^2((s_0, \infty), \rho^{-1}(s)\,ds)$, respectively. We regard $K_\tau$ as a $2 \times 2$ matrix Hilbert–Schmidt operator on $L^2(\rho) \oplus L^2(\rho^{-1})$ and note that $\varepsilon: L^2(\rho) \to L^2(\rho^{-1})$ as a consequence of the assumption that $\rho^{-1} \in L^1$.

More specifically, if $K_\tau = \begin{pmatrix} K_{11} & K_{12} \\ K_{21} & K_{22} \end{pmatrix}$, we regard $K_{11}$ and $K_{22}$ as trace class operators on $L^2(\rho)$ and $L^2(\rho^{-1})$, respectively, and the off-diagonal elements as Hilbert–Schmidt operators

$$K_{12}: L^2(\rho^{-1}) \to L^2(\rho) \quad \text{and} \quad K_{21}: L^2(\rho) \to L^2(\rho^{-1}).$$

Thus $\operatorname{tr} T$ denotes the sum of the traces of the diagonal elements of $T$. The regularized 2-determinant of a Hilbert–Schmidt operator $T$ with eigenvalues $\mu_k$ is defined by $\det_2(I-T) = \prod(1-\mu_k)e^{\mu_k}$ [cf. Gohberg and Krein (1969), Section IV.2]. Using this, one extends the operator definition of determinant to Hilbert–Schmidt operator matrices $T$ by setting

(183) $$\det(I-T) = \det_2(I-T) e^{-\operatorname{tr} T}.$$

As remarked in Tracy and Widom (2005), the resulting notion of $\det(I - K_\tau)$ is independent of the choice of $\rho$, and allows the derivation of Tracy and Widom (1998) that yields (48)–(50).

To analyze the convergence of $p_{N+1} = F_{N+1}(s_0)$ to $p_\infty = F_1(s_0)$, we note that

$$|p_{N+1} - p_\infty| \leq |p_{N+1}^2 - p_\infty^2|/p_\infty = C(s_0)|p_{N+1}^2 - p_\infty^2|,$$

so that we are led to the difference of determinants

(184) $$|F_{N+1}(s_0) - F(s_0)| \leq C(s_0)|\det(I - K_\tau) - \det(I - K_{GOE})|.$$

Our basic tool will be a Lipschitz bound on the matrix operator determinant for operators in the class $\mathcal{A}$ of $2 \times 2$ Hilbert–Schmidt operator matrices $A = (A_{ij}, i, j = 1, 2)$ on $L^2(\rho) \oplus L^2(\rho^{-1})$ whose diagonal entries are trace class.

PROPOSITION 3. *For $A, B \in \mathcal{A}$, we have*

$$|\det(I-A) - \det(I-B)| \leq C(A,B) \left\{ \sum_{i=1}^2 \|A_{ii} - B_{ii}\|_1 + \sum_{i \neq j} \|A_{ij} - B_{ij}\|_2 \right\}.$$

*The coefficient has the form* $C(A,B) = \sum_{j=1}^{2} c_{1j}(\operatorname{tr} A, \operatorname{tr} B) c_{2j}(A,B)$, *where* $c_{1j}$ *and* $c_{2j}$ *are continuous functions, the latter with respect to the strong (Hilbert–Schmidt norm) topology.*

PROOF. From determinant definition (183), we have

$$\det(I - A) - \det(I - B) = [\det_2(I - A) - \det_2(I - B)]e^{-\operatorname{tr} A} \qquad (185)$$
$$+ \det_2(I - B)e^{-\operatorname{tr} B}[e^{\operatorname{tr} B - \operatorname{tr} A} - 1].$$

A Lipschitz bound for the 2-determinant [Gohberg, Goldberg and Krupnik (2000), page 196 or Simon (1977), Theorem 6.5] gives

$$(186) \quad |\det_2(I - A) - \det_2(I - B)| \leq \|A - B\|_2 \exp\{\tfrac{1}{2}(1 + \|A\|_2 + \|B\|_2)^2\}.$$

The first term of (185) is thus bounded by $C_1(A,B)\|A - B\|_2$, where

$$C_1(A,B) = e^{-\operatorname{tr} A} \exp\{\tfrac{1}{2}(1 + \|A\|_2 + \|B\|_2)^2\}$$

has the requisite form. Since $\|A\|_2 \leq \sum_{i,j} \|A_{ij}\|_2$ and $\|A_{ii}\|_2 \leq \|A_{ii}\|_1$, the first term satisfies the stated bound.

For the second term, we have

$$|\operatorname{tr} A - \operatorname{tr} B| \leq \sum_i |\operatorname{tr}(A_{ii} - B_{ii})| \leq \sum_i \|A_{ii} - B_{ii}\|_1,$$

where we have used the fact that $\operatorname{tr} A = \operatorname{tr} A_{11} + \operatorname{tr} A_{22}$. Thus the second term has the required form $C_2(A,B) = c_{12}(\operatorname{tr} A, \operatorname{tr} B) c_{22}(A,B)$ with $c_{22}(A,B) = \det_2(I - B)$ and $c_1(x,y) = (e^{-x} - e^{-y})/(x - y)$. Bound (186) shows that $c_{22}$ has the necessary continuity. □

8.3. *Representation.* The next step is to establish a representation for $K_\tau(s,t)$ that facilitates the convergence argument. Our starting point is (49), which with the matrix definitions

$$L = \begin{pmatrix} I & -\partial_2 \\ \varepsilon_1 & T \end{pmatrix}, \qquad K^\varepsilon(x,y) = \begin{pmatrix} 0 & 0 \\ -\varepsilon(x-y) & 0 \end{pmatrix},$$

may be written in the form

$$K_{N+1,1}(x,y) = (LS_{N+1,1})(x,y) + K^\varepsilon(x,y).$$

In the unitary case, $S_{N,2}(x,y)$ transformed to $S_\tau(s,t) = e_N \bar{S}_\tau(s,t)$. In the orthogonal setting, we show that $S_{N+1,1}(x,y)$ transforms according to

$$(187) \qquad \tau'(s) S_{N+1}(\tau(s), \tau(t)) = e_N S_\tau^R(s,t),$$

where

$$(188) \qquad S_\tau^R(s,t) = \bar{S}_\tau(s,t) + \tfrac{1}{2} \phi_\tau(s)(\varepsilon \psi_\tau)(t).$$

To establish this, we begin by using the relation

(189) $$\tau'(s) = \sigma \cosh^{-2}(\mu + \sigma s) = \sigma[1 - \tau^2(s)],$$

in combination with the definition $\phi_\tau(s) = \check{\phi}_N(\tau(s))$ and (135) to obtain

(190) $$\begin{aligned}\phi_\tau(s) &= (\kappa_N \sigma_N)^{-1/2}(1 - \tau^2(s))^{1/2}\phi_N(\tau(s)) \\ &= (\sigma\kappa_N\sigma_N)^{-1/2}\sqrt{\tau'(s)}\phi_N(\tau(s)).\end{aligned}$$

Using (189), we may rewrite (50) in the form

$$\tau'(s)S_{N+1}(\tau(s),\tau(t)) \\ = \sqrt{\tau'(s)\tau'(t)}S_{N,2}(\tau(s),\tau(t)) + a_N\tilde\kappa_N\tau'(s)\tilde{\check\phi}_N(\tau(s))(\varepsilon\tilde{\check\phi}_{N-1})(\tau(t)).$$

The first term on the right-hand side equals $e_N \bar S(s,t)$, as may be seen from (129) and (136) in the unitary case. Turning to the components of the second right-hand side term, we use (189) and (190) to write

$$\tau'(s)\tilde{\check\phi}_N(\tau(s)) = \sqrt{\sigma}\sqrt{\tau'(s)}\phi_N(\tau(s)) = \sigma\sqrt{\kappa_N\sigma_N}\phi_\tau(s).$$

From the analog of (190) for $\psi_\tau$ and $x = \tau(t)$, we have

(191) $$\begin{aligned}\int_x^1 \tilde{\check\phi}_{N-1}(y)\,dy &= \int_t^\infty \tilde{\check\phi}_{N-1}(\tau(u))\tau'(u)\,du \\ &= \sigma\sqrt{\sigma_{N-1}\kappa_{N-1}}\int_t^\infty \psi_\tau(u)\,du.\end{aligned}$$

Consequently

$$(\varepsilon\tilde{\check\phi}_{N-1})(\tau(t)) = \sigma\sqrt{\sigma_{N-1}\kappa_{N-1}}(\varepsilon\psi_\tau)(t).$$

Comparing (51) and (138), we find that $e_N/2 = a_N\tilde\kappa_N\sigma^2\sqrt{\sigma_N\sigma_{N-1}\kappa_N\kappa_{N-1}}$. Gathering all this together, we obtain the result promised at (187)–(188):

$$\tau'(s)S_{N+1}(\tau(s),\tau(t)) = e_N\bar S(s,t) + (e_N/2)\phi_\tau(s)(\varepsilon\psi_\tau)(t) = e_N S_\tau^R(s,t).$$

We turn now to the elements of $LS_{N+1,1}(\tau(s),\tau(t))$. Temporarily write $\tilde S(s,t) = \tau'(s)S_{N+1,1}(\tau(s),\tau(t))$. Observing that $\varepsilon(\tau(s) - \tau(t)) = \varepsilon(s-t)$, we have

$$\begin{aligned}(\varepsilon_1 S_{N+1})(\tau(s),\tau(t)) &= \int \varepsilon(\tau(s) - \tau(u))S_{N+1}(\tau(u),\tau(t))\tau'(u)\,du \\ &= \int \varepsilon(s-u)\tilde S(u,t)\,du = \varepsilon_1\tilde S.\end{aligned}$$

Using such change of variables formulae at each matrix entry as needed, we obtain

$$\begin{pmatrix} I & -\partial_2 \\ \varepsilon_1 & T \end{pmatrix} S_{N+1,1}(\tau(s), \tau(t)) = \begin{pmatrix} \frac{1}{\tau'(s)} I & \frac{-1}{\tau'(s)\tau'(t)} \partial_2 \\ \varepsilon_1 & \frac{1}{\tau'(t)} T \end{pmatrix} \tilde{S}(s,t).$$

Note that the operators act with respect to variables $(x, y)$ on the left-hand side, and with respect to variables $(s, t)$ on the right. In terms of $L$, we write this as

$$(192) \quad (LS_{N+1,1})(\tau(s), \tau(t)) = \begin{pmatrix} 1/\tau'(s) & 0 \\ 0 & 1 \end{pmatrix} (L\tilde{S})(s,t) \begin{pmatrix} 1 & 0 \\ 0 & 1/\tau'(t) \end{pmatrix}.$$

We now make use of unimodular matrices $U(\tau) = \begin{pmatrix} \tau & 0 \\ 0 & 1/\tau \end{pmatrix}$. We have, for example, $U(a) K^\varepsilon U(b) = \begin{pmatrix} 0 & 0 \\ (b/a)\varepsilon & 0 \end{pmatrix}$. Setting $a = 1/\sqrt{\tau'(s)}$ and $b = \sqrt{\tau'(t)}$, and combining with (192), we obtain

$$K_\tau(s,t) = \sqrt{\tau'(s)\tau'(t)} (LS_{N+1,1} + K^\varepsilon)(\tau(s), \tau(t))$$
$$= U^{-1/2}(\tau'(s))(L\tilde{S} + K^\varepsilon)(s,t) U^{1/2}(\tau'(t)).$$

We now remark that the eigenvalues of $UKU^{-1}$ are the same as those of $K$, and so $\det(I - UKU^{-1}) = \det(I - K)$. Introduce

$$(193) \qquad q_N(s) = \sqrt{\frac{\tau'(s_0)}{\tau'(s)}} = \frac{\cosh(\mu + \sigma s)}{\cosh(\mu + \sigma s_0)},$$

and abbreviate $U(q_N(s))$ by $U_{q_N}(s)$. It follows that in place of $K_\tau$ we may use

$$\bar{K}_\tau(s,t) = U^{1/2}(\tau'(s_0)) K_\tau(s,t) U^{-1/2}(\tau'(s_0)) = U_{q_N}(s)(L\tilde{S} + K^\varepsilon)(s,t) U_{q_N}^{-1}(t).$$

Recalling (187), we may summarize by saying that $F_{N+1}(s_0) = \sqrt{\det(I - \bar{K}_\tau)}$, with

$$\bar{K}_\tau(s,t) = U_{q_N}(s)(e_N LS_\tau^R + K^\varepsilon)(s,t) U_{q_N}^{-1}(t).$$

REMARK. For later use, we define

$$(194) \qquad \beta_{N-1} = \tfrac{1}{2} \int_{-\infty}^\infty \psi_\tau = [\sigma \sqrt{\sigma_{N-1} \kappa_{N-1}}]^{-1} \tfrac{1}{2} \int_{-1}^1 \tilde{\phi}_{N-1},$$

where the second equality uses (191). Since, here, $N+1$ is even, Lemma 9 both gives an evaluation of $\beta_{N-1}$ and also shows that $\int_{-1}^1 \tilde{\phi}_N = 0$. The analog of the second equality above for $\phi_\tau$ then shows that $\int_{-\infty}^{\infty} \phi_\tau = 0$. Summarizing these remarks, we have

$$(195) \qquad (\varepsilon\phi_\tau)(t) = -\int_t^{\infty} \phi_\tau \quad \text{and} \quad (\varepsilon\psi_\tau)(t) = \beta_{N-1} - \int_t^{\infty} \psi_\tau$$

with

$$(196) \qquad \beta_{N-1}^{-2} = \sigma^2 \sigma_{N-1} \kappa_{N-1}^2 a_{N-1}(1 + \varepsilon_{N-1}).$$

REMARK. The need for care in the left tail is dramatized, for example, by

$$0 = \lim_{N\to\infty} \lim_{t\to-\infty} \int_t^{\infty} \phi_\tau \neq \lim_{t\to-\infty} \lim_{N\to\infty} \int_t^{\infty} \phi_\tau = \int_{-\infty}^{\infty} \mathrm{Ai} = 1.$$

Formula (194) and limit (220) show that there is a similar problem for $\int_t^{\infty} \psi_\tau$.

For the convergence argument, due to the oscillatory behavior of the Airy function in the *left* tail, it is helpful to rewrite expressions involving $\varepsilon$ in terms of the right-tail integration operator

$$(\tilde{\varepsilon}g)(s) = \int_s^{\infty} g(u)\,du.$$

Thus $\varepsilon g = \frac{1}{2}\int_{-\infty}^{\infty} g - \tilde{\varepsilon}g$, and we may rewrite (195) as

$$(197) \qquad \varepsilon\phi_\tau = -\tilde{\varepsilon}\phi_\tau, \qquad \varepsilon\psi_\tau = \beta_{N-1} - \tilde{\varepsilon}\psi_\tau.$$

For kernels $A(s,t)$, we have

$$(\varepsilon A)(s,t) = \tfrac{1}{2}\int_{-\infty}^{\infty} A(u,t)\,du - (\tilde{\varepsilon}_1 A)(s,t),$$

where, of course,

$$(\tilde{\varepsilon}_1 A)(s,t) = \int_s^{\infty} A(u,t)\,du.$$

Using integral representation (137) for $\bar{S}_\tau$ along with the (194) and $\int \phi_\tau = 0$,

$$\int_{-\infty}^{\infty} \bar{S}_\tau(u,t)\,du = \beta_{N-1}\int_0^{\infty} \phi_\tau(t+z)\,dz = \beta_{N-1}(\tilde{\varepsilon}\phi_\tau)(t).$$

As a result, we have

$$\varepsilon_1 \bar{S}_\tau = \tfrac{1}{2}\beta_{N-1} \otimes \tilde{\varepsilon}\phi_\tau - \tilde{\varepsilon}_1 \bar{S}_\tau.$$

From (197), we have

$$S_\tau^R = \bar{S}_\tau - \tfrac{1}{2}\phi_\tau \otimes \tilde{\varepsilon}\psi_\tau + \tfrac{1}{2}\phi_\tau \otimes \beta_{N-1},$$

and, combining the last two displays,

$$\varepsilon_1 S^R_\tau = -\tilde{\varepsilon}_1(\bar{S}_\tau - \tfrac{1}{2}\phi_\tau \otimes \tilde{\varepsilon}\psi_\tau) + \tfrac{1}{2}\beta_{N-1}(1 \otimes \tilde{\varepsilon}\phi_\tau - \tilde{\varepsilon}\phi_\tau \otimes 1).$$

Defining operator matrices

$$\tilde{L} = \begin{bmatrix} I & -\partial_2 \\ -\tilde{\varepsilon}_1 & T \end{bmatrix}, \qquad L_1 = \begin{bmatrix} I & 0 \\ -\tilde{\varepsilon} & 0 \end{bmatrix}, \qquad L_2 = \begin{bmatrix} 0 & 0 \\ \tilde{\varepsilon} & I \end{bmatrix},$$

we arrive at an expression for $LS^R_\tau$ that involves only right-tail integrations:

$$LS^R_\tau = \tilde{L}\left(\bar{S}_\tau - \frac{1}{2}\phi_\tau \otimes \tilde{\varepsilon}\psi_\tau\right) + \frac{\beta_{N-1}}{2}L_1\phi_\tau(s) + \frac{\beta_{N-1}}{2}L_2\phi_\tau(t).$$

Let us rewrite the limiting kernel $K_{GOE}$ in corresponding terms. Until the end of this Section 8, we will write $A(s)$ for the Airy function $\text{Ai}(s)$ to ease notation. For example,

$$IS(s,t) = -\tilde{\varepsilon}_1(S_A(s,t) - \tfrac{1}{2}A(s)(\tilde{\varepsilon}A)(t)) - \tfrac{1}{2}\tilde{\varepsilon}A(s) + \tfrac{1}{2}\tilde{\varepsilon}A(t)$$

and, assembling the other matrix entries correspondingly, we find

$$K_{GOE} = \tilde{L}(S_A - \tfrac{1}{2}A \otimes \tilde{\varepsilon}A) + \tfrac{1}{2}L_1 A(s) + \tfrac{1}{2}L_2 A(t) + K^\varepsilon.$$

SUMMARY. In summary, we may represent

$$K_\tau = U_{q_N}(s)[e_N(K^R_\tau + K^F_{\tau,1} + K^F_{\tau,2}) + K^\varepsilon]U^{-1}_{q_N}(t),$$
$$K_{GOE} = K^R + K^F_1 + K^F_2 + K^\varepsilon,$$

where we have

$$K^R_\tau = \tilde{L}[\bar{S}_\tau - \tfrac{1}{2}\phi_\tau \otimes \tilde{\varepsilon}\psi_\tau], \qquad K^R = \tilde{L}[S_A - \tfrac{1}{2}A \otimes \tilde{\varepsilon}A],$$
$$K^F_{\tau,1} = \tfrac{1}{2}\beta_{N-1}L_1[\phi_\tau(s)], \qquad K^F_1 = \tfrac{1}{2}L_1[A(s)],$$
$$K^F_{\tau,2} = \tfrac{1}{2}\beta_{N-1}L_2[\phi_\tau(t)], \qquad K^F_2 = \tfrac{1}{2}L_2[A(t)].$$

Our goal is to use inequalities (160)–(161) to obtain an $N^{-2/3}$ rate of convergence. Note in particular that $\psi_\tau = A + \Delta_N A' + O(N^{-2/3})$, and so define $A_N(s) = A(s) + \Delta_N A'(s)$. Expression (137) may be rewritten as

$$2\bar{S}_\tau = \phi_\tau \diamond \psi_\tau + \psi_\tau \diamond \phi_\tau,$$

where the convolution like operator $\diamond$ is defined in the obvious way. Replace $\phi_\tau$ by $A$ and $\psi_\tau$ by $A_N$ to define

$$S_{A_N} = \tfrac{1}{2}(A \diamond A_N + A_N \diamond A)$$
$$= A \diamond A + \tfrac{1}{2}\Delta_N(A \diamond A' + A' \diamond A).$$

We have
$$(A \diamond A' + A' \diamond A)(s,t) = \int_0^\infty \frac{d}{dz}[A(s+z)A(t+z)]\,dz = -A(s)A(t),$$
so that
$$S_{A_N} = S_A - \tfrac{1}{2}\Delta_N A \otimes A.$$

Since $\psi_\tau = A + \Delta_N A' + O(N^{-2/3})$, we will see below that it is convenient to set $A_N = A + \Delta_N A'$ and to write the difference
$$K_\tau^R - K^R = \tilde{\mathbf{L}}[\bar{S}_\tau - S_A + \tfrac{1}{2}\Delta_N A \otimes A] - \tfrac{1}{2}\tilde{\mathbf{L}}[\phi_\tau \otimes \tilde{\varepsilon}\psi_\tau - A \otimes \tilde{\varepsilon}A_N]$$
$$= \delta^{R,I} + \delta_0^F.$$

To organize the convergence argument, then, we describe the components of $K_\tau - K_{GOE}$ as

(198) $$K_\tau - K_{GOE} = \delta^{R,D} + \delta^{R,I} + \delta_0^F + \delta_1^F + \delta_2^F + \delta^\varepsilon$$

where, in addition to $\delta^{R,I}$ and $\delta_0^F$ previously defined,
$$\delta^{R,D} = e_N U_{q_N}(s) K_\tau^R U_{q_N}^{-1}(t) - K_\tau^R,$$
$$\delta^{R,I} = \tilde{L}[\bar{S}_\tau - S_A + \tfrac{1}{2}\Delta_N A \otimes A],$$
$$\delta_0^F = -\tfrac{1}{2}\tilde{L}[\phi_\tau \otimes \tilde{\varepsilon}\psi_\tau - A \otimes \tilde{\varepsilon}A_N],$$
$$\delta_i^F = e_N U_{q_N}(s) K_{\tau,i}^F U_{q_N}^{-1}(t) - K_i^F, \qquad i = 1,2,$$
$$\delta^\varepsilon = e_N U_{q_N}(s) K^\varepsilon U_{q_N}^{-1}(t) - K^\varepsilon.$$

### 8.4. Convergence.

#### 8.4.1. Operator bounds.
As a preliminary, we need some bounds on Hilbert–Schmidt and trace norms for repeated use. First, a remark taken verbatim from Tracy and Widom (2005): the norm of a rank-1 kernel $u(x)v(y)$, when regarded as an operator $u \otimes v$ taking a space $L^2(\rho_1)$ to a space $L^2(\rho_2)$ is given by

(199) $$\|u \otimes v\| = \|u\|_{2,\rho_2} \|v\|_{2,\rho_1^{-1}}.$$

(Here norm can be trace, Hilbert–Schmidt or operator norm, since all agree for a rank-1 operator.) Indeed the operator takes a function $h \in L^2(\rho_1)$ to $u(v,h)$, and so its norm is the $L^2(\rho_2)$ norm of $u$ times the norm of $v$ in the space dual to $L^2(\rho_1)$, which is $L^2(\rho_1^{-1})$.

Second, an operator $T: L^2(M,\mu) \to L^2(M,\mu')$ defined by
$$(Tf)(s) = \int K(s,t)f(t)\,d\mu(t)$$

has Hilbert–Schmidt norm given by

$$\|T\|_{HS}^2 = \int \int K^2(s,t)\,d\mu'(s)\,d\mu(t). \tag{200}$$

See, for example, Aubin (1979), Chapter 12.1, Proposition 1.

We use the following notation for a Laplace-type transform:

$$\mathcal{L}(\rho)[t] = \int_{s_0}^{\infty} e^{-tz}\rho(z)\,dz.$$

LEMMA 7. *Let $D$ be an operator taking $L^2(\rho_2)$ to $L^2(\rho_1)$ have kernel*

$$D(s,t) = \alpha(s)\beta(t)(a \diamond b)(s,t),$$

*where we assume, for $s \geq s_0$, that*

$$|\alpha(s)| \leq \alpha_0 e^{\alpha_1 s}, \qquad |\beta(s)| \leq \beta_0 e^{\beta_1 s}, \tag{201}$$

$$|a(s)| \leq a_0 e^{-a_1 s}, \qquad |b(s)| \leq b_0 e^{-b_1 s}. \tag{202}$$

*Assume that $\mathcal{L}(\rho_1)$ and $\mathcal{L}(\rho_2)$ both converge for $t > 0$, and that $a_1 > \alpha_1, b_1 > \beta_1$. Then*

$$\|D\|_{HS} \leq \frac{\alpha_0\beta_0 a_0 b_0}{a_1 + b_1}\{\mathcal{L}(\rho_1)[2(a_1 - \alpha_1)]\mathcal{L}(\rho_2)[2(b_1 - \beta_1)]\}^{1/2}.$$

*If $\rho_1 = \rho_2$, then the trace norm $\|D\|_1$ satisfies the same bound.*

PROOF. Substituting the bounds for $a$ and $b$, one finds

$$|(a \diamond b)(s,t)| \leq \frac{a_0 b_0}{a_1 + b_1} e^{-a_1 s - b_1 t}.$$

The Hilbert–Schmidt bound is a direct consequence of this, (200) and (201):

$$\|D\|_{HS}^2 = \int \int D^2(s,t)\rho_1(s)\rho_2(t)\,ds\,dt$$

$$\leq \left(\frac{\alpha_0\beta_0 a_0 b_0}{a_1 + b_1}\right)^2 \int_{s_0}^{\infty} e^{2\alpha_1 s - 2a_1 s}\rho_1(s)\,ds \int_{s_0}^{\infty} e^{2\beta_1 t - 2b_1 t}\rho_2(t)\,dt.$$

For the trace norm bound, we note that $D$ is an integral over $z$ of rank-1 kernels, and the norm of an integral is at most the integral of the norms. Thus, inserting (199),

$$\|D\| \leq \int_0^{\infty} \|\alpha(s)a(s+z)\beta(t)b(t+z)\|_1\,dz$$

$$\leq \int_0^{\infty} \|\alpha(\cdot)a(\cdot+z)\|_{2,\rho_1}\|\beta(\cdot)b(\cdot+z)\|_{2,\rho_1^{-1}}\,dz.$$

Now insert the bounds assumed at (201) and (202), so that, for example,

$$\|\alpha(\cdot)a(\cdot+z)\|_{2,\rho_1}^2 \le \alpha_0^2 a_0^2 e^{-2a_1 z}\mathcal{L}(\rho_1)[2(a_1-\alpha_1)].$$

The claimed bound for $\|D\|$ now follows after integration over $z$. □

We now make a particular choice of $\rho$ in order to facilitate the operator convergence arguments. For a $\gamma > 0$ to be specified later, let

(203) $$\rho(s) = 1 + e^{\gamma|s|}.$$

For notational convenience, we will let $\rho^+$ and $\rho^-$ be alternate symbols for $\rho$ and $\rho^{-1}$, respectively. With this choice of $\rho$, then for $\gamma < \tau$

(204) $$\mathcal{L}(\rho^\pm)[\tau] \le 2\int_{s_0}^\infty e^{-\tau z \pm \gamma|z|}\,dz \le \frac{4}{\tau-\gamma}e^{-\tau s_0 \pm \gamma|s_0|}.$$

Indeed, if $s_0 \ge 0$, our bound is immediate. If $s_0 < 0$, then split the integral into

$$2\int_{s_0}^0 e^{-\tau z \mp \gamma z}\,dz + \frac{2}{\tau \mp \gamma} \le \frac{4}{\tau-\gamma}e^{-(\tau \pm \gamma)s_0}.$$

We shall also use a related bound, proved similarly. For $|\tau| < \gamma$,

(205) $$\mathcal{L}(\rho^-)[\tau] \le \int_{s_0}^\infty e^{\tau z - \gamma|z|}\,dz \le \frac{2}{\gamma-\tau}e^{-(\gamma-\tau)s_{0+}},$$

where $s_{0+} = \max\{s_0, 0\}$.

COROLLARY 1. *Under the assumptions of Lemma 7, and if $\rho_1$ and $\rho_2$ therein are selected from $\{\rho, \rho^{-1}\}$, and if $\gamma < 2(a_1-\alpha_1), 2(b_1-\beta_1)$, then*

$$\|D\|_{HS}, \|D\|_1 \le C\frac{\alpha_0\beta_0 a_0 b_0}{a_1+b_1}e^{-(a_1+b_1-\alpha_1-\beta_1)s_0+\gamma|s_0|},$$

*where $C = C(a_1, \alpha_1, b_1, \beta_1, \gamma)$.*

CONSEQUENCE. We will make repeated use of Lemma 7 and Corollary 1 in the following way. If any one of $\alpha_0, \beta_0, a_0$ or $b_0$ is $O(N^{-2/3})$ while the others are uniformly bounded in $N$, and the bounds (201) and (202) apply, then the Hilbert–Schmidt (resp., trace) norms $\|D\|$ are $O(N^{-2/3})$. The convergence conditions for $\mathcal{L}(\rho_2)$ and $\mathcal{L}(\rho_1^{-1})$ follow from (181)–(182).

8.4.2. *Convergence details.*

PROOF OF THEOREM 4. Insert the conclusion of Proposition 3 into (184) to obtain

$$|F_{N+1}(s_0) - F_1(s_0)|$$

(206)
$$\leq C(s_0)C(K_\tau, K_{GOE})\bigg\{\sum_i \|K_{\tau,ii} - K_{GOE,ii}\|_1 + \sum_{i\neq j} \|K_{\tau,ij} - K_{GOE,ij}\|_2\bigg\}.$$

We exploit decomposition (198): the convergence of the matrix entries of $K_\tau - K_{GOE}$ is reduced to establishing the entrywise convergence of (i) terms involving integral kernels, $\delta^{R,D}$ and $\delta^{R,I}$, (ii) finite rank terms $\delta_0^F, \delta_1^F$ and $\delta_2^F$, and (iii) a term $\delta^\varepsilon$ involving versions of the convolution operator $\varepsilon$. We establish both Hilbert–Schmidt and trace norm bounds for the diagonal elements and Hilbert–Schmidt bounds for the off-diagonal entries. The distinction is moot for the finite rank terms $\delta_i^F$, and $\delta^\varepsilon$ involves only the $(2,1)$ entry, so the trace bounds are actually also needed only for the $\delta^R$ term.

For each term, we show $\|\delta_{ij}\| \leq N^{-2/3}$, so that the term $\{\cdot\}$ in (206) is $\leq CN^{-2/3}$. We have both $\|K_\tau - K_{GOE}\|_2$ and $\operatorname{tr} K_\tau - \operatorname{tr} K_{GOE}$ converging to 0 at $N^{-2/3}$ rate, so that $C(K_\tau, K_{GOE})$ remains bounded as $N \to \infty$.

$\delta^R$ *terms.* For both $\delta^{R,D}$ and $\delta^{R,I}$, we use Corollary 1 to establish the needed Hilbert–Schmidt and trace norm bounds for each entry in the $2 \times 2$ matrices comprising $\delta^{R,D}$ and $\delta^{R,I}$.

$\delta^{R,I}$ *term.* We have $\delta^{R,I} = \tilde{L}[\bar{S}_\tau - S_{A_N}]$, and

$$\bar{S}_\tau - S_{A_N} = (\phi_\tau - A)\diamond\psi_\tau + A\diamond(\psi_\tau - A_N) + (\psi_\tau - A_N)\diamond\phi_\tau + A_N\diamond(\phi_\tau - A).$$

In turn, for $\partial_2(\bar{S}_\tau - S_{A_N})$ we replace the second slot arguments $\psi_\tau, (\psi_\tau - A_N), \phi_\tau$ and $(\phi_\tau - A)$ by their derivatives, and for $\tilde{\varepsilon}(\bar{S}_\tau - S_{A_N})$, we replace the first slot arguments $(\phi_\tau - A), A, (\psi_\tau - A_N)$ and $A_N$ by their right-tail integrals.

Consider, for example, the first term $(\phi_\tau - A)\diamond\psi_\tau$. Use the abbreviation $D^{(k)}\psi$ to denote any of $\psi', \psi$ or $\tilde{\varepsilon}\psi$. Then we have the bounds

(207) $\quad |D^{(k)}(\phi_\tau - A)| \leq CN^{-2/3}e^{-s/4}, \qquad |D^{(k)}\psi_\tau| \leq Ce^{-s}.$

We apply Lemma 7 and Corollary 1 with $\alpha(s) = \beta(s) \equiv 1$ and with

$$a_0 = CN^{-2/3}, \qquad b_0 = C, \qquad a_1 = \tfrac{1}{4}, \qquad b_1 = 1.$$

The argument is entirely parallel for each of the second through fourth terms. Thus, if $D_{ij}$ denotes any matrix entry in any component of $\delta^{R,I}$, we obtain

(208) $\qquad \|D_{ij}\| \leq C^2 N^{-2/3} e^{-5s_0/4 + \gamma|s_0|}.$

Table 4

|  | $\alpha_0$ | $\alpha_1$ |
|---|---|---|
| 1 | 1 | 0 |
| $q_N^\pm(s)$ | $2e^{-\sigma s_0}$ | $\sigma$ |
| $(e_N - 1)q_N^\pm(s)$ | $CN^{-1}e^{-\sigma s_0}$ | $\sigma$ |
| $(q_N^\pm(s) - 1)$ | $CN^{-2/3}e^{-(\sigma+\delta)s_0}$ | $\sigma + \delta$ |

$\delta^{R,D}$ *term.* Decompose $K_\tau^R = K_\tau^{R,c} + K_\tau^{R,1}$ into "convolution" and "rank-1" terms

$$K_\tau^{R,c} = \tilde{L}\bar{S}_\tau, \qquad K_\tau^{R,1} = \tfrac{1}{2}\tilde{L}(\phi_\tau \otimes \tilde{\varepsilon}\psi_\tau),$$

respectively. Correspondingly, in the following telescoping decomposition, we have

$$\delta^{R,D} = (e_N - 1)Q_N(s)K_\tau^R Q_N^{-1}(t) + (Q_N(s) - I)K_\tau^R Q_N^{-1}(t) + K_\tau^R(Q_N^{-1}(t) - I)$$
$$= \delta^{R,D,c} + \delta^{R,D,1}.$$

The elements of the component terms of $\delta^{R,D,c}$ are all of the form $\alpha(s)\beta(t)(a \diamond b)(s,t)$. In order to apply Lemma 7, we verify conditions (201) and (202). Since $\tilde{L}\bar{S}_\tau = \begin{pmatrix} \bar{S}_\tau & -\partial_2 \bar{S}_\tau \\ -\tilde{\varepsilon}_1 \bar{S}_\tau & T\bar{S}_\tau \end{pmatrix}$, the terms $a \diamond b$ are all of the form

$$D^{(k)}\phi_\tau \diamond D^{(l)}\psi_\tau \quad \text{or} \quad D^{(k)}\psi_\tau \diamond D^{(l)}\phi_\tau, \qquad k = 0, -1; \ l = 0, 1.$$

All functions $D^{(k)}\phi_\tau$ and $D^{(k)}\psi_\tau$ satisfy (202) with $a_0 = b_0 = C$ and $a_1 = b_1 = 1$.

Inspecting the decomposition above, we see that the multipliers $\alpha(s)$ and $\beta(t)$ are chosen from the list $q_N^\pm(s), (e_N - 1)q_N^\pm(s)$ or $(q_N^\pm(s) - 1)$. To develop bounds for $q_N(s)$ and $q_N^{-1}(s)$, note that $c(a,b) = \cosh(a+b)/\cosh a$ satisfies, for all $a$ and $b$,

(209)
$$c(a,b) \quad \text{and} \quad 1/c(a,b) \leq 2e^{|b|},$$
$$|c(a,b) - 1| \quad \text{and} \quad |(1/c(a,b)) - 1| \leq 2be^{|b|}.$$

These inequalities are applied to $q_N(s) = c(\mu + \sigma s_0, \sigma(s - s_0))$, yielding

(210) $\qquad |q_N(s)| \quad \text{and} \quad |q_N^{-1}(s)| \leq 2e^{\sigma(s-s_0)},$

(211) $\qquad |q_N(s) - 1| \quad \text{and} \quad |q_N^{-1}(s) - 1| \leq 2\sigma(s - s_0)e^{\sigma(s-s_0)}.$

As a result, we may collect the bounds for $\alpha_0$ and $\alpha_1$ (resp., $\beta_0$ and $\beta_1$) in (201) in Table 4.

In the last line, we have used the bound $s - s_0 \leq \delta^{-1}e^{\delta(s-s_0)}$ for all $s \geq s_0$, where $\delta$ can be chosen arbitrarily small. Consequently, $C$ depends on $\delta$.

Denoting by $D_{ij}$ any matrix entry in any one of the terms of $\delta^{R,D,c}$, we therefore have from Corollary 1

$$\|D_{ij}\| \leq CN^{-2/3}e^{-2s_0+\gamma|s_0|}. \tag{212}$$

(Note the cancellation of terms of form $\sigma s_0$ or $(\sigma+\delta)s_0$ in the exponent.)

Many of the remaining steps will involve matrices of rank-1 operators on $L^2(\rho) \oplus L^2(\rho^{-1})$. Henceforth, we abbreviate the $L^2$ norms on $L^2(\rho)$ and $L^2(\rho^{-1})$ by $\|\cdot\|_+$ and $\|\cdot\|_-$, respectively. Let us record, using the remark leading to (199), the bound

$$\begin{pmatrix} \|a_{11} \otimes b_{11}\|_1 & \|a_{12} \otimes b_{12}\|_2 \\ \|a_{21} \otimes b_{21}\|_2 & \|a_{22} \otimes b_{22}\|_1 \end{pmatrix} \leq \begin{pmatrix} \|a_{11}\|_+ \|b_{11}\|_- & \|a_{12}\|_+ \|b_{12}\|_+ \\ \|a_{21}\|_- \|b_{12}\|_- & \|a_{22}\|_- \|b_{22}\|_+ \end{pmatrix}. \tag{213}$$

On the left, $\|\cdot\|_1$ and $\|\cdot\|_2$ denote trace and Hilbert–Schmidt norms. Indeed, apply (199) to $a_{ij} \otimes b_{ij} : L^2(\rho_j) \to L^2(\rho_i)$, where $\rho_1 = \rho$ and $\rho_2 = \rho^{-1}$.

Turning now to the $\delta^{R,D,1}$ term, and observing that

$$\tilde{L}(\phi_\tau \otimes \tilde{\varepsilon}\psi_\tau) = \begin{pmatrix} \phi_\tau \otimes \tilde{\varepsilon}\psi_\tau & -\phi_\tau \otimes \psi_\tau \\ -\tilde{\varepsilon}\phi_\tau \otimes \tilde{\varepsilon}\psi_\tau & \tilde{\varepsilon}\psi_\tau \otimes \phi_\tau \end{pmatrix},$$

we see that every term in $\delta^{R,D,1}$ is of the form $a \otimes b$, where

$$a(s)b(t) = \ell_j(s)\zeta_j(s)\zeta_k(t)\ell_k(t),$$

where $\zeta_j(s)$ and $\zeta_k(t)$ are chosen from the list $\{\phi_\tau, \tilde{\varepsilon}\phi_\tau, \psi_\tau, \tilde{\varepsilon}\psi_\tau\}$ and $\ell_j(s), \ell_k(s)$ are chosen from the rows of Table 4, with the conventions that $j$ and $k$ indicate rows and that one of $j, k$ equals 1 or 2 and the other equals 3 or 4. If we abbreviate the bounds summarized in the table by $|\ell_j(s)| \leq C_{jN}e^{-l_j s_0}e^{l_j s}$, and then use Corollary 1, we obtain

$$\|\ell_j \zeta_j\|^2_{\rho^\pm} \leq C^2_{jN} e^{-2l_j s_0} C^2 \int_{s_0}^{\infty} e^{2l_j s - 2s} \rho^\pm(s)\, ds$$

$$\leq CC^2_{jN} e^{-2l_j s_0} e^{-2(1-l_j)s_0 + \gamma|s_0|},$$

so that

$$\|a \otimes b\| \leq \|\ell_j \zeta_j\|_{\rho^\pm} \|\ell_k \zeta_k\|_{\rho^\pm}$$
$$\leq CN^{-2/3}e^{-2s_0+\gamma|s_0|}. \tag{214}$$

Applying these bounds to $\tilde{L}(a \otimes \tilde{\varepsilon}b)$, we obtain

$$\|\tilde{L}(a \otimes \tilde{\varepsilon}b)\| \leq \begin{pmatrix} A_+ B_- & A_+ B_+ \\ A_- B_- & A_+ B_- \end{pmatrix}, \tag{215}$$

where

$$A_+ = \|a\|_+, \qquad B_+ = \|b\|_+,$$
$$A_- = \|\tilde{\varepsilon}a\|_-, \qquad B_- = \|\tilde{\varepsilon}b\|_-.$$

For $\delta_0^F$ write
$$-2\delta_0^F = \tilde{L}[\phi_\tau \otimes \tilde{\varepsilon}(\psi_\tau - A_N) + (\phi_\tau - A) \otimes \tilde{\varepsilon} A_N]$$
so that we may apply (215), first with $a = \phi_\tau$, $b = \psi_\tau - A_N$ and then with $a = \phi_\tau - A$, $b = A_N$.

In the first case, we have from Corollary 1
$$A_+^2 = \|\phi_\tau\|_+^2 = \int_{s_0}^\infty \phi_\tau^2 \rho \leq C^2 \int_{s_0}^\infty e^{-2s+\gamma|s|}\,ds$$
$$\leq \frac{4}{2-\gamma} C^2 e^{-2s_0+\gamma|s_0|},$$
with the same bound applying also to $A_-^2$. In a similar vein,
$$B_-^2 = \|\tilde{\varepsilon}(\psi_\tau - A_N)\|_\pm^2 \leq C^2 N^{-4/3} \int_{s_0}^\infty e^{-s/2 \pm \gamma|s|}\,ds$$
$$\leq \frac{8}{1-2\gamma} C^2 N^{-4/3} e^{-s_0/2+\gamma|s_0|}$$
with the same bound also for $B_+^2$. Hence

(216) $$A_\pm B_\pm \leq C(\gamma) N^{-2/3} e^{-5s_0/4 + \gamma|s_0|}.$$

The same bound works for the second case, with $a = \phi_\tau - A, b = A_N$, as well.

$\delta_i^F$ *term.* We have
$$2\delta_1^F = \begin{pmatrix} u_{N1} \otimes q_N^{-1} - A \otimes 1 & 0 \\ -u_{N2} \otimes q_N^{-1} + \tilde{\varepsilon} A \otimes 1 & 0 \end{pmatrix},$$
$$2(\delta_2^F)^t = \begin{pmatrix} 0 & q_N^{-1} \otimes u_{N2} - 1 \otimes \tilde{\varepsilon} A \\ 0 & q_N^{-1} \otimes u_{N1} - 1 \otimes A \end{pmatrix}$$
with
$$u_{N1} = \gamma_N q_N \phi_\tau, \qquad u_{N2} = \gamma_N q_N^{-1} \tilde{\varepsilon} \phi_\tau, \qquad \gamma_N = e_N \beta_{N-1}.$$

Using (213), we find that the norms of the first column of $\delta_1^F$ are bounded by
$$\begin{pmatrix} \|u_{N1} - A\|_+ \|q_N^{-1}\|_- + \|A\|_+ \|q_N^{-1} - 1\|_- \\ \|u_{N2} - \tilde{\varepsilon} A\|_- \|q_N^{-1}\|_- + \|\tilde{\varepsilon} A\|_- \|q_N^{-1} - 1\|_- \end{pmatrix}$$
while the norms of the second column of $(\delta_2^F)^t$ are bounded by the same quantities, with the rows interchanged.

From (210),
$$\|q_N^{-1}\|_-^2 \leq 4 \int_{s_0}^\infty e^{2\sigma(s-s_0)-\gamma|s|}\,ds \leq C e^{-2\sigma s_0 - (\gamma - 2\sigma)s_{0+}},$$

so that

$$\|q_N^{-1}\|_- \le \begin{cases} Ce^{-\gamma s_0/2}, & s_0 \ge 0, \\ Ce^{-\sigma s_0}, & s_0 < 0. \end{cases} \tag{217}$$

Using (211),

$$\|q_N^{-1} - 1\|_-^2 \le C(\delta)\sigma^2 e^{-2(\sigma+\delta)s_0} \int_{s_0}^\infty e^{2(\sigma+\delta)s - \gamma|s|}\,ds,$$

so that

$$\|q_N^{-1} - 1\|_- \le \begin{cases} CN^{-2/3}e^{-\gamma s_0/2}, & s_0 \ge 0, \\ CN^{-2/3}e^{-(\sigma+\delta)s_0}, & s_0 < 0. \end{cases} \tag{218}$$

Using (158) for $\phi_\tau$ and (211) for $q_N - 1$,

$$\|(q_N - 1)\phi_\tau\|_+^2 \le C\sigma^2 e^{-2(\sigma+\delta)s_0} \int_{s_0}^\infty e^{-2(1-\sigma-\delta)s + \gamma|s|}\,ds,$$

so that, using Corollary 1,

$$\|(q_N - 1)\phi_\tau\|_+ \le C\sigma e^{-s_0 + (\gamma/2)|s_0|}.$$

A similar argument gives

$$\|q_N \phi_\tau\|_+ \le C e^{-s_0 + (\gamma/2)|s_0|}.$$

We have

$$\|u_{N1} - A\|_+ \le |\gamma_N - 1|\|q_N \phi_\tau\|_+ + \|(q_N - 1)\phi_\tau\|_+ + \|\phi_\tau - A\|_+. \tag{219}$$

First, we show that $|\gamma_N - 1| = O(N^{-1})$. For $e_N$, refer to (138). For $\beta_{N-1}$, we exploit (196) and Lemma 9 (noting that $N + 1$ is even) to write

$$\begin{aligned}\beta_{N-1}^{-2} &= \sigma^2 \sigma_{N-1} \kappa_{N-1}^2 a_{N-1}(1 + \varepsilon_N) = \sigma_{N-1}^3 \omega_{N-1}^2 \kappa_{N-1}^2 a_{N-1}(1 + \varepsilon_N) \\ &= 1 + \varepsilon_N,\end{aligned} \tag{220}$$

where we have used $\sigma = \sigma_N \omega_N = \sigma_{N-1}\omega_{N-1}(1 + \varepsilon_N)$ and then (139).

Assembling the bounds developed just above decomposition (219) along with (207) yields

$$\|u_{N1} - A\|_+ \le (CN^{-1}e^{-s_0} + CN^{-2/3}e^{-s_0} + CN^{-2/3}e^{-s_0/4})e^{\gamma|s_0|/2}$$

and

$$\|u_{N1} - A\|_+ \|q_N^{-1}\|_- \le CN^{-2/3}e^{-(1/2+\gamma)s_0/2}. \tag{221}$$

With a similar decomposition,

$$\|u_{N2} - \tilde\varepsilon A\|_- \le |\gamma_N - 1|\|q_N^{-1}\tilde\varepsilon\phi_\tau\|_- + \|(q_N^{-1} - 1)\tilde\varepsilon\phi_\tau\|_- + \|\tilde\varepsilon(\phi_\tau - A)\|_-$$

and these terms are bounded exactly as are the corresponding terms in $\|u_{N1} - A\|_+$.

Finally, observe that

$$\|\tilde{\varepsilon}A\|_{\pm}^2 \leq C^2 \int_{s_0}^{\infty} e^{-2s+\gamma|s|} \, ds \leq C(\gamma) e^{-2s_0+\gamma|s_0|},$$

so that

(222) $\qquad \|\tilde{\varepsilon}A\|_{-}\|q_N^{-1} - 1\|_{-} \leq \begin{cases} CN^{-2/3} e^{-s_0}, & s_0 > 0, \\ CN^{-2/3} e^{-(1+\gamma)s_0}, & s_0 < 0. \end{cases}$

$\delta^{\varepsilon}$ *term.* The only nonzero entry in the $\delta^{\varepsilon}$ term is

$$\delta_{21}^{\varepsilon} = [e_N q_N^{-1}(s) q_N^{-1}(s) - 1]\varepsilon(s-t)$$
$$= [(e_N - 1)q_N^{-1}(s)q_N^{-1}(t) + (q_N^{-1}(s) - 1)q_N^{-1}(t) + (q_N^{-1}(t) - 1)]\varepsilon(s-t)$$
$$= (\tilde{\varepsilon}_1 + \tilde{\varepsilon}_2 + \tilde{\varepsilon}_3)(s,t).$$

Each of $\tilde{\varepsilon}_j(s,t)$ has the form $\varepsilon_{ab}(s,t) = a(s)b(t)\varepsilon(s-t)$. We regard $\varepsilon_{ab}$ as an operator mapping $L^2(\mathbb{R}, d\mu(t) = \rho(t)\,dt) \to L^2(\mathbb{R}, d\mu'(s) = \rho^{-1}(s)\,ds)$, thus

$$(\varepsilon_{ab}f)(s) = \int [a(s)\varepsilon(s-t)b(t)/\rho(t)]f(t)\,d\mu(t)$$

so that $K(s,t)$ in (200) has the form $\varepsilon_{ab}(s,t)/\rho(t)$, and

$$\|\varepsilon_{ab}\|_{HS}^2 = \tfrac{1}{4} \iint a^2(s)b^2(t)\rho^{-2}(t)\,d\mu'(s)\,d\mu(t) = \tfrac{1}{4}\|a\|_{-}^2 \|b\|_{-}^2.$$

Hence

$$\|\delta_{21}^{\varepsilon}\|_{HS} \leq |e_N - 1|\|q_N^{-1}\|_{-}^2 + \|q_N^{-1}\|_{-}\|q_N^{-1} - 1\|_{-} + \|1\|_{-}\|q_N^{-1} - 1\|_{-}$$

and from (138), (217) and (218), it is bounded by

(223) $\qquad \begin{cases} C(N^{-1} + N^{-2/3} + N^{-2/3})e^{-\gamma s_0}, & s_0 \geq 0, \\ C(N^{-1} + N^{-2/3} + N^{-2/3})e^{-(2\sigma+\gamma)s_0}, & s_0 < 0. \end{cases}$

We remark that the exponential right-tail bound here is possible due to the assumption (203) on the weight function $\rho$.

At last we can assemble the bounds obtained in (208), (212), (214), (216), (221), (222) and (223). Each term has a component $CN^{-2/3}$ where $C$ depends on $\gamma, \delta$ and of course $\alpha(N)/N$ and $\beta(N)/N$. We only track the tail dependence on $s_0$ for $s_0 > 0$. With regard to that dependence, (214) is dominated by (208), and so (206) is bounded by

$$CN^{-2/3}(e^{-5s_0/4+\gamma s_0} + e^{-(1+2\gamma)s_0/4} + e^{-s_0} + e^{-\gamma s_0}).$$

It remains to choose a suitable value of $\gamma$; it is clear that $\gamma = \tfrac{1}{2}$ yields a bound $CN^{-2/3}e^{-s_0/2}$. [The choice of $\gamma$ could be further optimized, but this is perhaps not worthwhile until best bounds are found on the exponential rate in (160)–(161).]

8.4.3. *Summing up.*

THEOREM 4. *The previous subsection established that*

$$|F_{N+1}(s_0) - F_1(s_0)| \leq CN^{-2/3} e^{-s_0/2},$$

with the constant $C$ depending on $s_0$ when $s_0 < 0$, and, referring to (178) and recalling that $x = \tau(s) = \tanh(\mu + \sigma s)$, we have

$$F_{N+1}(s_0) = P\{\tau^{-1}(x_{(1)}) \leq s_0\} = P\{(\tanh^{-1} x_{(1)} - \mu)/\sigma \leq s_0\}.$$

The $JOE(N+1, \alpha, \beta)$ setting is linked to the $JUE(N, \alpha, \beta)$ via (50), and through equating $\mu = u_N$ and $\sigma = \tau_N$. The $u$-scale centering and scaling values are related to their $x$-scale versions via (133):

$$u_N = \tanh^{-1} x_N, \qquad \tau_N = \sigma_N/(1 - x_N^2).$$

Finally, on the $x$-scale, we have from (76) and (112):

$$x_N = -\cos(\varphi + \gamma), \qquad \sigma_N^3 = \frac{2\sin^4(\varphi + \gamma)}{\kappa_N^2 \sin\varphi \sin\gamma}.$$

Hence, for all $\mu$ we arrive at

$$\mu = u_N = \tanh^{-1} x_N = \frac{1}{2}\log\frac{1+x_N}{1-x_N} = \frac{1}{2}\log\frac{1-\cos(\varphi+\gamma)}{1+\cos(\varphi+\gamma)}$$

(224)
$$= \log\tan(\varphi+\gamma)/2$$

and

(225) $$\sigma^3 = \tau_N^3 = \frac{\sigma_N^3}{(1-x_N^2)^3} = \frac{2}{\kappa_N^2}\frac{1}{\sin^2(\varphi+\gamma)\sin\varphi\sin\gamma}.$$

THEOREM 1. While Theorem 1 is just a relabeling of Theorem 4, it may be useful to collect the parameterizations and formula leading to the centering and scaling expressions (5) and (6).

First we identify the double Wishart setting of Definition 1 with the appropriate JOE. Identification (27) was made under the additional assumption that $n \geq p$. If $n < p$, we use identity (2) and density (7) with parameters $(p', m', n') = (n, m+n-p, p)$. In either case, then, we use $JOE(N+1, \alpha, \beta)$ with the identification

$$\begin{pmatrix} N+1 \\ \alpha \\ \beta \end{pmatrix} = \begin{pmatrix} p \wedge n \\ m-p \\ |n-p| \end{pmatrix}.$$

Noting that $|n-p| = p \vee n - p \wedge n$, we have

$$\kappa_N = 2(N+1) + \alpha + \beta - 1 = m+n-1,$$
$$\alpha + \beta = m+n-2(p \wedge n), \qquad \alpha - \beta = m+n-2(p \vee n).$$

Thus, the (20) defining $\gamma$ and $\varphi$ become

$$\cos\gamma = \frac{\alpha+\beta}{\kappa} = 1 - \frac{2(p\wedge n - 1/2)}{m+n-1},$$

$$\cos\varphi = \frac{\alpha-\beta}{\kappa} = 1 - \frac{2(p\vee n - 1/2)}{m+n-1},$$

which yield the half-angle forms (6).

Recall that $\theta_i = (1+x_i)/2$ so that

$$w_i = \log\frac{\theta_i}{1-\theta_i} = \log\frac{1+x_i}{1-x_i} = 2\tanh^{-1}x_i.$$

Thus $\mu_p = 2\mu$ and $\sigma_p = 2\sigma$, and so we recover (5) from (224) and (225).

## APPENDIX

**A.1. Proof of Lemma 1.** From (57) and (58), respectively,

$$\frac{h_N}{h_{N-1}} = \frac{(N+\alpha)(N+\beta)}{N(N+\alpha+\beta)}\frac{2N+\alpha+\beta-1}{2N+\alpha+\beta+1}$$

and

$$\frac{l_{N-1}}{l_N} = \frac{2N(N+\alpha+\beta)}{(2N+\alpha+\beta)(2N+\alpha+\beta-1)},$$

whence

$$(226) \qquad a_N = 2\left[\frac{N(N+\alpha)(N+\beta)(N+\alpha+\beta)}{\kappa(\kappa-1)^2(\kappa-2)}\right]^{1/2}.$$

Now use the $(a,b)$ parameters defined at (63)–(65). We have, for example, $(N+\alpha)/\kappa = (1+a-1/(2N_+))/(2+a+b)$, so that

$$a_N = 2\left[\frac{(1+a)(1+b)(1+a+b)}{(2+a+b)^4}\right]^{1/2}[1+O(N^{-1})].$$

From (68) and (69) follow

$$\frac{\sin^2\varphi}{4} = \frac{(1+a)(1+b)}{(2+a+b)^2}, \qquad \frac{\sin^2\gamma}{4} = \frac{1+a+b}{(2+a+b)^2},$$

and now (86) is immediate.

**A.2. Choice of $f$ and $g$ in (74) for LG approximation.** We elaborate on consequences of the key remark that the $O(1/\kappa)$ error bound (94) is available only if the integral $\mathcal{V}(\zeta) = \int_\zeta^{\zeta(b)} |\psi(v)v^{-1/2}|\,dv < \infty$. In view of Remark A, we consider convergence at both endpoints, corresponding to $a = -1$ as well as $b = 1$.

We refer to arguments in Olver (1974) to show that these convergence requirements lead to the specific choices of $f$ and $g$ made in (74). Indeed, remarks in Olver [(1974), Section 11.4.1] show that it suffices to show finite total variation as $x \to \pm 1$ of the error control function

$$F(x) = \int_{x_*}^{x} f^{-1/4} \frac{d^2}{dx^2}(f^{-1/4}) - gf^{-1/2} \, dx$$

associated with the LG approximations of Olver [(1974), Chapter 6.2] at both endpoints. In turn, the discussion of Olver [(1974), Section 6.4.3] shows that this is obtained if

$$f_0 = \lim_{x \to c}(c-x)^2 f(x) > 0 \quad \text{and} \quad g_0 = \lim_{x \to c}(c-x)^2 g(x) = -1/4,$$

for both endpoints $c = \pm 1$.

Return to (70)–(72) and write

$$\frac{\mathfrak{n}(x)}{4(1-x^2)^2} = \frac{u^2 \mathsf{F}(x) + \mathsf{G}(x)}{4(1-x^2)^2} = u^2 f(x) + g(x),$$

where $\mathsf{F}(x)$ and $\mathsf{G}(x)$ are quadratic polynomials which clearly determine $f$ and $g$. The requirements on $g_0$ imply that $\mathsf{G}(\pm 1) = -4$. The coefficient of $x^2$ in $\mathfrak{n}(x)$ is $(2N + \alpha + \beta + 1)^2 - 1 = \kappa^2 - 1$. If, for convenience, we take $\mathsf{F}(x)$ to be monic, then it is natural to set the large parameter $u = \kappa$. The condition on $f_0$ follows from the fact that the zeros $x_\pm$ of $\mathfrak{n}(x)$ lie in the interior of $[-1, 1]$. Further, $\mathsf{G}''(x) = -2$, which with the two previous constraints implies $\mathsf{G}(x) = -3 - x^2$. Since $\mathfrak{n}(1) = 4\alpha^2 - 4$ and $\mathfrak{n}(-1) = 4\beta^2 - 4$, we conclude easily that $\mathsf{F}(1) = 4\alpha^2/\kappa^2 = 4\lambda^2$ and $\mathsf{F}(-1) = 4\beta^2/\kappa^2 = 4\mu^2$, thus arriving at the expressions (74) for $f(x)$ and $g(x)$.

**A.3. Error control function.** Clearly, for $x_0 \le x \le 1$ we have $\mathcal{V}(\zeta(x)) \le \mathcal{V}(\zeta(x_0))$, and it is our goal here to show that $\mathcal{V}(\zeta(x); \lambda, \mu)$ is uniformly bounded for $(\lambda, \mu) \in D_\delta$.

We first observe that

$$\mathcal{V}(\zeta(x_0)) = \mathcal{V}_{[x_0, x_1]}(H) + \mathcal{V}_{[x_1, 1]}(H),$$

where $x_0 = \frac{1}{2}(x_+ + x_-)$ is defined before (96) and $x_1 > x_+$ will be specified below. We then note that since $H'(x) = -\dot\zeta |\zeta|^{-1/2} \psi(\zeta)$, we have

$$\mathcal{V}_{[x_0, x_1]}(H) = \int_{\zeta(x_0)}^{\zeta(x_1)} \frac{|\psi(\zeta)|}{|\zeta|^{1/2}} \, d\zeta.$$

Our approach is to use the fact that $(\zeta, \lambda, \mu) \to \psi(\zeta; \lambda, \mu)$ is continuous, and hence is bounded, by $M(\delta)$ say, on the compact set of $(\zeta, \lambda, \mu)$ for which both $\zeta \in [\zeta(x_0(\lambda, \mu)), \zeta(x_1)]$ and $(\lambda, \mu) \in D_\delta$—we use here the continuity of $x_0$ in $(\lambda, \mu)$ and of $\zeta$ in $x$. For $(\lambda, \mu)$ in $D_\delta$, there exist finite bounds $\zeta_-(\delta) \leq \zeta(x_0(\lambda, \mu))$ and $\zeta_+(\delta) \geq x_1$ so that

$$\mathcal{V}_{[x_0, x_1]}(H) \leq M(\delta) \int_{\zeta_-(\delta)}^{\zeta_+(\delta)} |\zeta|^{-1/2}\, d\zeta \leq M_1(\delta).$$

Continuity of $\psi(\zeta; \lambda, \mu)$ is a consequence of Olver (1974), Lemma 11.3.1, and the continuous dependence of $f$ and hence $\zeta$ on $(\lambda, \mu)$. Indeed, the lemma uses the decomposition

$$\hat{f}(x; \lambda, \mu) = \{p(x; \lambda, \mu)\}^2 \{\tfrac{3}{2} q(x; \lambda, \mu)\}^{-2/3},$$

where

$$p(x) = (x - x_-)^{1/2}/(2(1 - x^2)) \quad \text{and}$$
$$q(x) = (x - x_+)^{-3/2} \int_{x_+}^{x} (t - x_+)^{1/2} p(t)\, dt$$

and as $x \to x_+$, $q(x) \to \tfrac{2}{3} p(x_+)$. For $(\lambda, \mu) \in D_\delta$, we have $p(x_+) \geq \delta_5(\delta) > 0$, and so the continuous dependence of $x_-$ and $x_+$ on $(\lambda, \mu)$ carries through to $\psi$.

Turning to $\mathcal{V}_{[x_1, 1]}(H)$, note from Olver (1974), (11.4.01), that $H(x) = F(x) + (5/24)\zeta^{-3/2}$, so that

$$\mathcal{V}_{[x_1, 1]}(H) = \mathcal{V}_{[x_1, 1]}(F) + \tfrac{5}{24} \zeta(x_1)^{-3/2}.$$

Since $\zeta(x_1)$ is bounded below on $D_\delta$, it remains to bound

$$\mathcal{V}_{[x_1, 1]}(F) = \int_{x_1}^{1} |\mathcal{L}(x)|\, dx,$$

where $\mathcal{L} = f^{-1/4}(f^{-1/4})'' - gf^{-1/2}$.

To organize the calculation write

$$f(x) = (1-x)^{-2}\mathsf{f}(x), \qquad \mathsf{f}(x) = \frac{\lambda^2}{4} + (1-x)\mathsf{f}_1(x),$$
$$\mathsf{f}_1(x) = \frac{(1+x)(\lambda^2 - 1) + 2\mu^2}{4(1+x)^2},$$
$$g(x) = (1-x)^{-2}\mathsf{g}(x), \qquad \mathsf{g}(x) = -\tfrac{1}{4} + (1-x)\mathsf{g}_1(x),$$
$$\mathsf{g}_1(x) = \frac{-1}{2(1+x)^2},$$

from which one obtains

$$f^{-1/4}(f^{-1/4})'' = -\tfrac{1}{4}(1-x)^{-1}\mathsf{f}^{-1/2} + \mathcal{B},$$
$$gf^{-1/2} = (1-x)^{-1}\mathsf{f}^{-1/2}\mathsf{g},$$

where

$$\mathcal{B} = \frac{1}{4}\mathsf{f}^{-1/2}\left[\frac{\mathsf{f}'}{\mathsf{f}} - (1-x)\left\{\frac{\mathsf{f}''}{\mathsf{f}} - \frac{5}{4}\left(\frac{\mathsf{f}'}{\mathsf{f}}\right)^2\right\}\right].$$

Since we have arranged the decomposition $\kappa^2 f + g$ precisely so that $\mathsf{g} = -\frac{1}{4} + (1-x)\mathsf{g}_1$, the $(1-x)^{-1}$ term cancels and

$$\mathcal{L} = f^{-1/4}(f^{-1/4})'' - gf^{-1/2} = -\mathsf{f}^{-1/2}\mathsf{g}_1 + \mathcal{B}.$$

Consequently, to show that $\mathcal{L}$ is bounded uniformly over $D_\delta$ on $[x_1, 1]$, it is enough to choose $x_1$ close enough to 1 so that

$$\inf \mathsf{f}(x) \geq \delta_4(\delta) > 0,$$

where the infimum is taken over $x \in [x_1, 1]$ and $(\lambda, \mu) \in D_\delta$. And indeed, then $|\mathsf{g}_1|, |\mathsf{f}'|$ and $|\mathsf{f}''|$ are uniformly bounded for such $(x, \lambda, \mu)$.

**A.4. Behavior of LG transform as $x \to 1$.** Write the leading term $f(x)$ of (74) in terms of

(227) $$\sqrt{f(x)} = \frac{R(x)}{2(1-x^2)}, \qquad R(x) = \sqrt{(x-x_+)(x-x_-)}.$$

Thus

(228) $$I(x) := (2/3)\zeta^{3/2} = \frac{1}{2}\int_{x_+}^{x} \frac{R(x')}{1-x'^2}\,dx'.$$

PROPOSITION 4. *Let $N$ be fixed. As $x \to 1$,*

(229) $$4I(x) = 4\int_{x_+}^{x} \sqrt{f(x')}\,dx' = \frac{2a}{2+a+b}\log(1-x)^{-1} + c_{0N} + o(1),$$

*where*

(230) $$c_{0N} = \frac{2}{2+a+b}\log\left[\frac{(2a^2)^a(1+b)^{1+b}}{(1+a)^{1+a}(1+a+b)^{1+a+b}}\right].$$

[Recall that the turning points $x_\pm$ of (75) and (76) are related by $\lambda, \mu$ of (63) to $a, b$ defined at (64).]

PROOF. To ease notation, introduce new variables $s = 1+x$ and $t = 1-x$, which are both positive for $|x| \leq 1$. With slight abuse, we set

(231) $$R(s) = \sqrt{(s-s_+)(s-s_-)}, \qquad s_\pm = 1+x_\pm,$$

(232) $$R(t) = \sqrt{(t-t_-)(t-t_+)}, \qquad t_\pm = 1-x_\mp,$$

and note that $R(s) = R(t) = R(x)$ for $|x| \leq 1$. Consequently,

$$(233) \qquad 4I(x) = \int_{s_+}^{s} \frac{R(s')}{s'}\,ds' + \int_{t}^{t_-} \frac{R(t')}{t'}\,dt'.$$

The following "elementary" indefinite integral formula may be derived from Gradshteyn and Ryzhik (1980), 2.267, 2.261 and 2.266. Let $u$ denote any of the variables $x, s$ or $t$, so that $R(u) = \sqrt{(u-u_+)(u-u_-)}$ and set

$$(234) \qquad \bar{u} = (u_+ + u_-)/2, \qquad \hat{u}^2 = u_+ u_-.$$

Then

$$\int \frac{R(u)}{u}\,du = R(u) - \hat{u}\log\{u^{-1}|\bar{u}u - \hat{u}^2 - \hat{u}R(u)|\} - \bar{u}\log|u - \bar{u} + R(u)|.$$

Before using this to evaluate the definite integrals in (233) note that

$$(235) \qquad \operatorname{sgn}[\bar{u}u - \hat{u}^2 - \hat{u}R(u)] = \operatorname{sgn}[\bar{u}u - \hat{u}^2],$$

$$(236) \qquad \operatorname{sgn}[u - \bar{u} + R(u)] = \operatorname{sgn}[u - \bar{u}]$$

and all four quantities are positive if $u \geq u_+$ and negative if $u \leq u_-$. [For (236), let $\Delta = (u_+ - u_-)/2$ and observe that $(u - \bar{u})^2 - R^2(u) = \Delta^2 > 0$, while for (235), note that

$$(\bar{u}u - \hat{u}^2)^2 - \hat{u}^2 R^2(u) = u^2\Delta^2 > 0.]$$

As a result, we have

$$\int_{s_+}^{s} \frac{R(s')}{s'}\,ds' = R(s) - \hat{s}\log\left[\frac{s_+}{s} \cdot \frac{\bar{s}s - \hat{s}^2 - \hat{s}R(s)}{\bar{s}s_+ - \hat{s}^2}\right] - \bar{s}\log\left[\frac{s - \bar{s} + R(s)}{s_+ - \bar{s}}\right]$$

with an analogous expression for the integral in $t$ in (233), namely

$$\int_{t}^{t_-} \frac{R(t')}{t'}\,dt' = -R(t) + \hat{t}\log\left[\frac{t_-}{t} \cdot \frac{\hat{t}R(t) + \hat{t}^2 - \bar{t}t}{\hat{t}^2 - \bar{t}t_-}\right] + \bar{t}\log\left[\frac{\bar{t} - t - R(t)}{\bar{t} - t_-}\right].$$

Define $\delta$ through the equations $\delta^{-1} = \bar{s} - s_- = t_+ - \bar{t} = x_+ - \bar{x}$. Adding the two previous displays yields

$$(237) \quad 4I(x) = -\hat{s}\log\{\delta s^{-1}(\bar{s}s - \hat{s}^2 - \hat{s}R(s))\} - \bar{s}\log\{\delta(s - \bar{s} + R(s))\}$$

$$(238) \qquad\qquad + \hat{t}\log\{\delta t^{-1}(\hat{t}R_t + \hat{t}^2 - \bar{t}t)\} + \bar{t}\log\{\delta(\bar{t} - t - R(t))\}.$$

As $x \nearrow 1$, we have $s \nearrow 2$ and $t \searrow 0$, and so, noting that $R(1) = \sqrt{t_+ t_-} = \hat{t}$, we obtain $4I(x) = \hat{t}\log t^{-1} + c_{0N} + o(1)$. This is the desired approximation (229), with

$$(239) \qquad c_{0N} = -\hat{s}\log T_1 - \bar{s}\log T_2 + +\hat{t}\log T_3 + \bar{t}\log T_4,$$

where

$$T_1 = \tfrac{1}{2}\delta(2\bar{s} - \hat{s}^2 - \hat{s}\hat{t}), \qquad T_2 = \delta(\bar{t} + \hat{t}), \qquad T_3 = 2\delta\hat{t}^2, \qquad T_4 = \delta(\bar{t} - \hat{t}).$$

To convert the previous expression for $c_{0N}$ to that given in terms of $a, b$ in (230), we proceed via the angle parameters $\gamma, \varphi$ of (66). In preparation, we set out, in parallel for the $s$ and $t$ variables, some relations that follow from the definitions, the fact that $x_\pm = -\cos\varphi\cos\gamma \pm \sin\varphi\sin\gamma$ and some algebra:

$$\text{(240)} \quad \hat{s}^2 = s_+ s_- = (1+x_+)(1+x_-) \qquad \hat{t}^2 = (1-x_+)(1-x_-)$$

$$\text{(241)} \quad = (1-\cos\varphi\cos\gamma)^2 \qquad\qquad = (1+\cos\varphi\cos\gamma)^2$$

$$\quad\quad - \sin^2\varphi\sin^2\gamma \qquad\qquad\qquad - \sin^2\varphi\sin^2\gamma$$

$$\text{(242)} \quad = (\cos\gamma - \cos\varphi)^2, \qquad\qquad = (\cos\gamma + \cos\varphi)^2,$$

so that

$$\text{(243)} \quad \hat{s} = \cos\gamma - \cos\varphi = 2\mu, \qquad \hat{t} = \cos\gamma + \cos\varphi = 2\lambda,$$

$$\text{(244)} \quad \bar{s} = 1 - \cos\varphi\cos\gamma = 1+\bar{x}, \qquad \bar{t} = 1 + \cos\varphi\cos\gamma = 1 - \bar{x}$$

and, since $\hat{s} + \hat{t} = 2\cos\gamma$,

$$\text{(245)} \quad \bar{s} - \hat{s}(\hat{s}+\hat{t})/2 = 1 - \cos\varphi\cos\gamma - (\cos\gamma - \cos\varphi)\cos\gamma = \sin^2\gamma,$$

while

$$\bar{t} + \hat{t} = (1+\cos\gamma)(1+\cos\varphi), \qquad \bar{t} - \hat{t} = (1-\cos\gamma)(1-\cos\varphi).$$

From (68) and (69) for $\cos\gamma$ and $\cos\varphi$,

$$\frac{1 \pm \cos\gamma}{\sin\gamma} = (1+a+b)^{\pm 1/2}, \qquad \frac{1 \pm \cos\varphi}{\sin\varphi} = \left(\frac{1+a}{1+b}\right)^{\pm 1/2}$$

and since (78) shows that $\delta^{-1} = \sin\varphi\sin\gamma$, we find that

$$\delta(\bar{t} \pm \hat{t}) = \frac{1 \pm \cos\gamma}{\sin\gamma}\frac{1 \pm \cos\varphi}{\sin\varphi} = \left[(1+a+b)\frac{1+a}{1+b}\right]^{\pm 1/2}.$$

Thus $T_4 = 1/T_2$, and so

$$-\bar{s}\log T_2 + \bar{t}\log T_4 = \log[(1+a)^{-1}(1+b)^{-1}(1+a+b)^{-1}].$$

Using now (245) and (240) and similar trigonometric manipulations,

$$\text{(246)} \quad T_1 = \frac{\sin\gamma}{\sin\varphi} = \frac{\sqrt{1+a+b}}{\sqrt{1+a}\sqrt{1+b}},$$

$$\text{(247)} \quad T_3 = \frac{2(\cos\gamma + \cos\varphi)^2}{\sin\gamma\sin\varphi} = \frac{2a^2}{\sqrt{1+a+b}\sqrt{1+a}\sqrt{1+b}},$$

so that

$$-\hat{s}\log T_1 + \hat{t}\log T_3 = -\mu\log[(1+a)^{-1}(1+b)^{-1}(1+a+b)]$$
$$+ \lambda\log[4a^4(1+a)^{-1}(1+b)^{-1}(1+a+b)^{-1}].$$

We obtain expression (230) for $c_{0N}$ from (239) by combining the previous displays and computing $1 \pm \lambda \pm \mu$ using (65). $\square$

LARGEST EIGENVALUE IN MULTIVARIATE ANALYSIS    75**A.5. Identification of $c_N$.** We first remark that as $\zeta \to \infty$ when $x \to 1$, we may substitute the large $x$ behavior of $\mathrm{Ai}(x)$ given by $\mathrm{Ai}(x) \sim [2\sqrt{\pi} x^{1/4}]^{-1} \times \exp\{-(2/3)x^{3/2}\}$ to obtain

$$(248) \qquad w_2(x,\kappa) \sim [2\sqrt{\pi}]^{-1} \kappa^{-1/6} f^{-1/4}(x) \exp\{-(2/3)\kappa\zeta^{3/2}\}.$$

Consequently, using (228), we may express $c_N$ in terms of the limit

$$c_N = \lim_{x \to 1^-} w_N(x) \cdot 2\sqrt{\pi} \kappa^{1/6} f^{1/4}(x) \exp\{\kappa I(x)\}.$$

Consider first the dependence on $x$. From (70) and (59), we have as $x \to 1$,

$$(249) \qquad w_N(x) \sim w_N (1-x)^{(\alpha+1)/2}, \qquad w_N = 2^{(\beta+1)/2} \binom{N+\alpha}{N}.$$

Since $R(x) \to R(1) = \hat{t} = 2\lambda$ [compare (240) and (234)], we have from (227) that

$$f^{1/4}(x) \sim \sqrt{\lambda/2}(1-x)^{-1/2},$$

while Proposition 4 along with (65) and (63) implies

$$\exp\{\kappa I(x)\} \sim e^{\kappa c_{0N}/4}(1-x)^{-\alpha/2}.$$

Multiply the last three displays: the resulting exponent of $(1-x)$ is identically 0. Consequently

$$(250) \qquad c_N = w_N \cdot 2\sqrt{\pi} \kappa^{1/6} \sqrt{\lambda/2} e^{\kappa c_{0N}/4}.$$

Using [ ] to denote the quantity in brackets in (230), we have

$$e^{\kappa c_{0N}/4} = [\,]^{N_+/2} = 2^{\alpha/2} \left[\frac{a^a}{(1+a)^{1+a}}\right]^{N_+} \left[\frac{(1+a)^{1+a}(1+b)^{1+b}}{(1+a+b)^{1+a+b}}\right]^{N_+/2}.$$

LEMMA 8. *Let $N_+ = N+1/2, \alpha = N_+ a$ and $\beta = N_+ b$. There exist bounded remainders $\theta_1(a)$ and $\theta_2(a,b)$ such that*

$$(251) \qquad \binom{N+\alpha}{N} = \frac{1}{\sqrt{2\pi N_+ a}} \left[\frac{(1+a)^{1+a}}{a^a}\right]^{N_+} \exp\left\{\frac{\theta_1}{N}\right\}$$

*and*

$$(252) \qquad \frac{(N+\alpha)!(N+\beta)!}{N!(N+\alpha+\beta)!} = \left[\frac{(1+a)^{1+a}(1+b)^{1+b}}{(1+a+b)^{1+a+b}}\right]^{N_+} \exp\left\{\frac{\theta_2}{N}\right\}.$$

PROOF. Stirling's approximation $x! = \sqrt{2\pi} e^{-x} x^{x+1/2} e^{\theta/x}$ has $0 \leq \theta \leq 1/12$. Consequently

$$\frac{(N+\alpha)!}{N!\alpha!}$$
$$= \frac{1}{\sqrt{2\pi}} \frac{[N_+(1+a) - 1/2]^{N_+(1+a)}}{[N_+ - 1/2]^{N_+}[N_+ a]^{N_+ a+1/2}} \exp\{\theta(a)/N\}$$
$$= \frac{1}{\sqrt{2\pi N_+ a}} \left[\frac{(1+a)^{1+a}}{a^a}\right]^{N_+} \frac{(1 - 1/[2N_+(1+a)])^{N_+(1+a)}}{(1 - 1/[2N_+])^{N_+}} \exp\{\theta(a)/N\}.$$

The result (251) follows from the relation $(1 - v/N)^N = \exp\{-v + \theta v^2 N^{-1}\}$ where $0 \leq \theta \leq 1$ for $N > 2v$. The argument for (252) is completely analogous. □

Combining Lemma 8 with (249) of $w_N$,

$$e^{\kappa c_{0N}/4} = \frac{e^{\theta''/N}}{w_N} \left(\frac{\kappa_N h_N}{2\pi\alpha}\right)^{1/2}.$$

Substituting this into (250), we finally obtain (99).

LEMMA 9.

$$(253) \qquad \int_{-1}^{1} \tilde{\phi}_N = 2(\kappa_N a_N)^{-1/2}(1 + \varepsilon_N) \qquad \text{for } N \text{ even}$$

and 0 for $N$ odd.

PROOF. From Nagao and Forrester (1995), (A.7), we have

$$(254) \quad \int_{-1}^{1} (1-x)^{(\alpha-1)/2}(1+x)^{(\beta-1)/2} P_N^{\alpha,\beta}(x) \, dx$$
$$= 2^{(\alpha+\beta)/2} \frac{\Gamma((N+\alpha+1)/2)\Gamma((N+\beta+1)/2)}{\Gamma((N+\alpha+\beta+1)/2)\Gamma((N+2)/2)}$$

if $N$ is even, and zero if $N$ is odd. [Identify our parameters $(N, \alpha, \beta)$ with NF's $(n, 2b+1, 2a+1)$ after noting that they use the opposite convention for Jacobi polynomial indices: our $P_N^{\alpha,\beta}$ is their $P_n^{(2a+1, 2b+1)}$.]

The function $\tilde{\phi}_N$ equals $h_N^{-1/2}$ times the integrand of (254), and so after combining this integral with expression (57) for $h_N$, we obtain

$$\int_{-1}^{1} \tilde{\phi}_N = \left[\frac{\kappa}{2(N+1)(N+\alpha+\beta+1)}\right]^{1/2} \frac{r(N+\alpha)r(N+\beta)}{r(N+1)r(N+\alpha+\beta+1)},$$

where

(255) $$r(N) = \frac{\Gamma((N+1)/2)}{\sqrt{\Gamma(N+1)}} = \left(\frac{2\pi}{N}\right)^{1/4} 2^{-N/2} e^{\theta/N},$$

with the last equality following from Stirling's formula as in the proof of Lemma 8. Substituting (255), we further find

$$\int_{-1}^{1} \tilde{\phi}_N = \sqrt{2\kappa}[N(N+\alpha)(N+\beta)(N+\alpha+\beta+1)]^{-1/4} e^{\theta/N}$$

$$= 2(\kappa_N a_N)^{-1/2}(1+\varepsilon_N)$$

after exploiting (226). □

**A.6. Proof of (128).** We show that $S_\tau$ has the same eigenvalues as $S_N$. Indeed, suppose that $g \in L_2[x_0, 1)$ satisfies $S_N g = \lambda g$. Set

(256) $$h(s) = \sqrt{\tau'(s)} g(\tau(s)).$$

First observe that

$$\int_{s_0}^{\infty} h^2(s) \, ds = \int_{s_0}^{\infty} g^2(\tau(s)) \tau'(s) \, ds = \int_{x_0}^{1} g^2(x) \, dx,$$

so that $h \in L_2[s_0, \infty)$ if and only if $g \in L_2[x_0, 1)$. In addition

$$(S_\tau h)(s) = \int_{s_0}^{\infty} S_\tau(s,t) h(t) \, dt = \sqrt{\tau'(s)} \int_{s_0}^{\infty} S_N(\tau(s), \tau(t)) g(\tau(t)) \tau'(t) \, dt$$

$$= \sqrt{\tau'(s)} \int_{t_0}^{1} S_N(\tau(s), y) g(y) \, dy = \sqrt{\tau'(s)} \lambda g(\tau(s)) = \lambda h(s).$$

**A.7. Proof of bounds for $\sigma_N \check{\phi}'_N(x)$.**

PROOF OF (109). Differentiate (106) to obtain

$$\sigma_N \check{\phi}'_N(x) = \sigma_N \bar{e}_N r'_N(x)[\text{Ai}(\kappa^{2/3}\zeta) + \varepsilon_2(x, \kappa)]$$
(257) $$\qquad + \bar{e}_N r_N(x)[\text{Ai}'(\kappa^{2/3}\zeta) \sigma_N \kappa^{2/3} \dot{\zeta}(x) + \sigma_N \partial_x \varepsilon_2(x, \kappa)]$$
$$= D_{N1} + D'_{N1}.$$

Using (104) to rewrite $\sigma_N \kappa^{2/3} \zeta$ as $\dot\zeta/\dot\zeta_N$, we further decompose the difference $\sigma_N \check{\phi}'_N(x) - \text{Ai}(s_N)$ as $\sum_{i=1}^{5} D_{Ni}$, with the new terms given by

$$D_{N2} = [\bar{e}_N r_N(x) - 1][\dot\zeta/\dot\zeta_N(x)] \text{Ai}'(\kappa^{2/3}\zeta),$$
$$D_{N3} = [\dot\zeta/\dot\zeta_N(x) - 1] \text{Ai}'(\kappa^{2/3}\zeta),$$
$$D_{N4} = \text{Ai}'(\kappa^{2/3}\zeta) - \text{Ai}'(s_N),$$
$$D_{N5} = \bar{e}_N r_N(x) \sigma_N \partial_x \varepsilon_2(x, \kappa).$$

We first observe from (94) and the uniform bound on $\mathcal{V}$ that

$$|\operatorname{Ai}(\kappa^{2/3}\zeta) + \varepsilon_2(x,\kappa)| \leq CM(\kappa^{2/3}\zeta)E^{-1}(\kappa^{2/3}\zeta) \tag{258}$$

and that, using (117) and (113),

$$\left|\frac{r'_N}{r_N}(x)\right| = \frac{1}{2}\left(\frac{\dot\zeta(x)}{\dot\zeta_N}\right)^{-1}\left|\frac{\ddot\zeta(x)}{\dot\zeta_N}\right|^{-1} \leq C.$$

As a result, combining the two previous bounds with the argument used for $E_{N3}$, we obtain

$$\begin{aligned}|D_{N1}| &\leq \sigma_N |(r'_N/r_N)(x)| \cdot r_N(x) M(\kappa^{2/3}\zeta) E^{-1}(\kappa^{2/3}\zeta) \\ &\leq C\sigma_N e^{-s_N} \leq CN^{-2/3} e^{-s_N/2}.\end{aligned} \tag{259}$$

Before turning to $D_{N2}$ and $D_{N3}$, we first remark, using $|\operatorname{Ai}'(x)| \leq N(x)E^{-1}(x)$, that on $[s_L, s_1 N^{1/6}]$,

$$|\operatorname{Ai}'(\kappa^{2/3}\zeta)| \leq N(\kappa^{2/3}\zeta) E^{-1}(\kappa^{2/3}\zeta) \leq C s_N^{1/4} e^{-s_N}. \tag{260}$$

Indeed, we bound $N(\kappa^{2/3}\zeta)$ by using $N(x) \leq C|x|^{1/4}$ and (114) to conclude that $|\kappa^{2/3}\zeta| \leq 2s_N$. The bound for $E^{-1}(\kappa^{2/3}\zeta)$ uses Proposition 3.

Combining (116), (117), (260) and (115), we find

$$|D_{N2}| \leq C(1+s_N)\sigma_N \cdot 2 \cdot C s_N^{1/4} e^{-s_N} \leq CN^{-2/3} e^{-s_N/2},$$
$$|D_{N3}| \leq C s_N \sigma_N \cdot C s_N^{1/4} e^{-s_N} \leq CN^{-2/3} e^{-s_N/2}.$$

$D_{N4}$ is treated in exactly the same manner as the $E_{N2}$ term above, additionally using the equation $\operatorname{Ai}''(x) = x\operatorname{Ai}(x)$.

Using (95), we can rewrite $D_{N5}$ as

$$|D_{N5}| \leq C\bar{e}_N \kappa_N^{-1} \cdot r_N(x) \sigma_N \kappa_N^{2/3} \hat{f}^{1/2}(x) \cdot N(\kappa^{2/3}\zeta) E^{-1}(\kappa^{2/3}\zeta).$$

From (110) and (104), we note that $\sigma_N \kappa_N^{2/3} \hat{f}^{1/2}(x) = \dot\zeta(x)/\dot\zeta_N$ and in combination with (117)

$$r_N(x)\sigma_N \kappa_N^{2/3}\hat{f}^{1/2}(x) = [\dot\zeta(x)/\dot\zeta_N]^{1/2} \leq \sqrt{2} \tag{261}$$

on $[s_L, s_1 N^{1/6}]$. Bringing in (260), we conclude

$$|D_{N5}| \leq C\kappa_N^{-1} \cdot \sqrt{2} \cdot C s_N^{1/4} e^{-s_N} \leq CN^{-2/3} e^{-s_N/2}.$$

### A.8. Proofs of bounds for $\phi'_\tau, \psi'_\tau$.

*Preliminaries on $r_N$ and $r'_N$.* Starting from (107), we find that

$$r'_N(x) = -\frac{1}{2}\left(\frac{\dot\zeta(x)}{\dot\zeta_N}\right)^{-1/2}\frac{\ddot\zeta(x)}{\dot\zeta(x)} \quad \text{and} \quad \frac{r'_N(x)}{r_N(x)} = -\frac{1}{2}\frac{\ddot\zeta(x)}{\dot\zeta(x)}.$$

Writing $I(\sqrt{f})$ for $\int_{x_+}^x \sqrt{f}$, taking logarithms in (96) and differentiating yields

$$(\log\zeta)' = \frac{2}{3}\frac{\sqrt{f}}{I(\sqrt{f})} \quad \text{and} \quad \frac{(\log\zeta)''}{(\log\zeta)'} = \frac{1}{2}\frac{f'}{f} - \frac{\sqrt{f}}{I(\sqrt{f})}.$$

From this one readily finds that

$$\frac{\ddot\zeta}{\dot\zeta} = \frac{(\log\zeta)''}{(\log\zeta)'} + (\log\zeta)' = \frac{1}{2}\frac{f'}{f} - \frac{1}{3}\frac{\sqrt{f}}{I(\sqrt{f})}.$$

By straightforward algebra and bounds on both $f'/f$ and $\sqrt{f}/I(\sqrt{f})$, one can check that for $x > x_+$,

$$\left|\frac{\ddot\zeta}{\dot\zeta}(x)\right| \leq \frac{C}{(x-x_+)(1-x^2)},$$

where $C$ depends on $x_+$ and $x_-$. Recalling that $\tau'(s) = \sigma(1 - \tau^2(s))$, we have, for $x = \tau(s) = x_N + \sigma_N s_N(s)$,

(262) $$\tau'(s)\left|\frac{r'_N(x)}{r_N(x)}\right| \leq \frac{C\sigma}{\sigma_N s_N(s)} \leq CN^{-1/6}$$

since for $s \geq s_1 N^{1/6}$, we have $s_N(s) \geq CN^{1/6}$.

*Bound for $\phi'_\tau(s)$.* The differentiated function $\phi'_\tau(s) = \breve\phi'_N(\tau(s))\tau'(s) = \tilde D_{N1}(s) + \tilde D'_{N1}(s)$ may be written in the form (257) with $\sigma_N$ replaced by $\tau'(s)$. The analog of (259) is

$$|\tilde D_{N1}(s)| \leq \tau'(s)\left|\frac{r'_N(x)}{r_N(x)}\right| r_N(x) M(\kappa^{2/3}\zeta) E^{-1}(\kappa^{2/3}\zeta).$$

On $[s_L, s_1 N^{1/6}]$, this is bounded by $CN^{-2/3}e^{-s_N/2}$ exactly as in the local bound case. For $s > s_1 N^{1/6}$, we use (262) together with (162) and (163) as above to get

$$|\tilde D_{N1}(s)| \leq CN^{-1/6} \cdot c_0/\sqrt{r} \cdot Ce^{-s} \leq Ce^{-s}.$$

Using (95), (110) and the uniform bound on $\mathcal{V}$,

$$|\tilde D'_{N1}(s)| \leq \bar e_N r_N(x)\tau'(s)[\text{Ai}'(\kappa^{2/3}\zeta)\kappa^{2/3}\dot\zeta + |\partial_x\varepsilon_2(x,\kappa)|]$$
$$\leq \bar e_N r_N(x)\tau'(s)\kappa^{2/3}\hat f^{1/2}(x) N(\kappa^{2/3}\zeta) E^{-1}(\kappa^{2/3}\zeta).$$

From (261) and (107),

$$r_N(x)\tau'(s)\kappa^{2/3}\hat{f}^{1/2}(x) = \sigma_N^{-1}\tau'(s)r_N^{-1}(x).$$

Using the $N$-asymptotics from (93) and then (107), we have

$$r_N^{-1}(x)N(\kappa^{2/3}\zeta) \leq Cr_N^{-1}(x)\kappa_N^{1/6}\zeta^{1/4} = C[\kappa_N\sigma_N\sqrt{f(x)}]^{1/2}.$$

Using $x = \tanh(u_N + \tau_N s)$ and (154) along with (155) and (104), we find

$$\kappa_N\sigma_N\sqrt{f(x)} \leq \kappa_N\sigma_N^{3/2}[s + \varepsilon_N(s)]^{1/2}/(1-x^2) \leq Cs^{1/2}/(1-x^2).$$

Since $\tau'(s) = \sigma[1 - \tau^2(s)]$, we conclude that

$$\sigma_N^{-1}\tau'(s)r_N^{-1}(x)N(\kappa^{2/3}\zeta) \leq C(\sigma/\sigma_N)s^{1/4}$$

and hence that

$$|\tilde{D}'_{N1}(s)| \leq Cs^{1/4}e^{-s} \leq Ce^{-s/2}.$$

*Bounds for* $\phi'_\tau(s) - \mathrm{Ai}'(s), \psi'_\tau(s) - \mathrm{Ai}'(s)$. Again, the real work is on $[s_L, s_1 N^{1/6}]$. For $\phi'_\tau(s)$, since $\mu = u_N$, $\sigma = \tau_N$, the bound needed is already established at (153). For $\psi'_\tau(s)$, we follow the approach taken for $\psi_\tau(s)$, differentiating $\psi_\tau(t) = \phi_{N-1}(u_{N-1} + \tau_{N-1}t)$ to yield

$$\psi'_\tau(t) = \tau_{N-1}\phi'_{N-1}(u_{N-1} + \tau_{N-1}t')(dt'/dt)$$
$$= [\mathrm{Ai}'(t') + O(N^{-2/3}e^{-t/2})][1 + O(N^{-1})],$$

using (153) and $dt'/dt = \tau_N\tau_{N-1}^{-1} = 1 + O(N^{-1})$. We now argue exactly as at (171) and (172), increasing each order of derivative by one. Since $\mathrm{Ai}''(t) = t\mathrm{Ai}(t)$, we nevertheless obtain the same bounds as before, so the proof of (161) follows.

**Acknowledgments.** Special thanks to Peter Forrester for an introduction to the literature relating orthogonal to unitary ensembles, and especially for pointing out (50), and more generally for conversations during visits he hosted to Melbourne, and also for providing early drafts of his book manuscript. Thanks also to Noureddine El Karoui, Plamen Koev, Zongming Ma, Bala Rajaratnam, Sasha Soshnikov and Craig Tracy for conversations and suggestions. An early version of the main results of this paper was presented at an Oberwolfach meeting in March 2002.

Department of Statistics
Sequoia Hall
390 Serra Mall
Stanford University
Stanford, California 94305-4065
E-mail: imj@stanford.edu